\documentclass[11pt]{article}
\pdfoutput=1
\usepackage[T1]{fontenc}
\usepackage[utf8]{inputenc}
\usepackage{mathpazo,microtype,amsmath,amssymb,amsthm}
\usepackage[margin=1in]{geometry}
\usepackage{booktabs,tabularx,tikz,needspace,etoolbox,xurl}
\AtBeginEnvironment{thebibliography}{\small}
\usetikzlibrary{arrows.meta,decorations.markings}
\tikzset{mid arrow/.style={postaction={decorate},
 decoration={markings,mark=at position 0.5 with {\arrow{#1}}}}}
\usepackage[colorlinks=true,linkcolor=blue,citecolor=blue,urlcolor=blue]{hyperref}
\newtheorem{theorem}{Theorem}[section]
\newtheorem{proposition}[theorem]{Proposition}
\newtheorem{lemma}[theorem]{Lemma}
\newtheorem{corollary}[theorem]{Corollary}
\newtheorem{maintheorem}{Theorem}

\theoremstyle{definition}
\newtheorem{definition}[theorem]{Definition}
\newtheorem{example}[theorem]{Example}
\newtheorem{question}[theorem]{Question}
\newtheorem{conjecture}[theorem]{Conjecture}
\theoremstyle{remark}
\newtheorem{remark}[theorem]{Remark}
\newcommand{\kk}{\Bbbk}
\newcommand{\m}{\mathfrak m}
\newcommand{\Tor}{\operatorname{Tor}}
\newcommand{\pd}{\operatorname{pd}}
\newcommand{\depth}{\operatorname{depth}}
\newcommand{\A}{\mathcal A}
\newcommand{\Podd}{\mathcal P}
\title{Additive Relations and Betti Growth in Goldbach Rings}
\author{Michael Nelson\thanks{Email: \href{mailto:mnelso32@gmail.com}{mnelso32@gmail.com}.}}
\date{}
\begin{document}
\maketitle
\begin{abstract}
The binary Goldbach ring over a field $\kk$ is the quotient of
$R=\kk[x_p:p\text{ an odd prime}]$
by the relations $x_px_q=x_{p'}x_{q'}$ whenever $p+q=p'+q'$, with
$\deg x_p=(p,1)$. We prove that if $g$ is a map from the odd primes to an
abelian group and $g(p)+g(q)$ depends only on $p+q$, then
$g(p)=g(3)+\frac{p-3}{2}(g(5)-g(3))$ for every odd prime $p$.
This settles K\'atai and Phong's conjecture on shifted prime arguments
and proves, for odd primes, the affine classification they conjectured
for binary sums.
The universal group of binary prime relations is free of rank two.
We obtain a minimal bigraded free resolution from compatible resolutions
of finite truncations and determine the growth of its linear Betti numbers
in two regimes. For fixed monomial degree $d\ge3$, as the Goldbach degree
$n\to\infty$ through $n\equiv d\pmod2$, we prove
$\beta_{d-1,n,d}\sim(d-1)T_d(n)/d!$, where $T_d(n)$ counts ordered
representations by $d$ odd primes. Explicit Koszul cycles from common
prime shifts span a positive fraction of this strand.
For finite supports $\A_m=\{p\le m:p\text{ an odd prime}\}$, fix the
homological degree $i\ge1$ and sum over all Goldbach degrees.
As $m\to\infty$, the linear-strand totals are asymptotic to
$i|\A_m|^{i+1}/(i+1)!$ when all binary relations on $\A_m$ are retained,
and to $(1-2^{-i})|\A_m|^{i+1}/(i+1)!$ when only relations of common
weight at most $m$ are retained. These leading estimates are uniform
over fields. The interpolation theorem gives eventual depth at most two
for the variable truncations; both truncation families have
asymptotically maximal projective dimension.
\end{abstract}
\noindent{\small\textit{2020 Mathematics Subject Classification.}
Primary 13D02; Secondary 11P32, 13F65, 11B13, 11A25.\par
\par}
\section{Introduction}
\label{sec:introduction}

The equality $3+17=7+13$ gives two decompositions of $20$ into odd
primes. Multiplying the corresponding monomials by $x_3$ gives a move
between decompositions of $23$ into three odd primes. A different move,
$3+7=5+5$, now gives
\[
 x_3^2x_{17}=x_3x_7x_{13}=x_5^2x_{13}.
\]
These equalities hold in the binary Goldbach ring over a field $\kk$
\[
 R=\kk[x_p:p\text{ an odd prime}],\qquad
 G=R/(x_px_q-x_{p'}x_{q'}\mid p+q=p'+q').
\]
Each move preserves the sum of the labels and the number of factors.
Accordingly, we give $x_p$ bidegree $(p,1)$ and call the two entries
the Goldbach degree and the monomial degree.
Goldbach's conjecture says that $G_{n,2}\ne0$ for every even $n\ge6$;
here we study the relations among these representations.
The monomials of $R_{n,d}$ correspond to unordered representations
of $n$ as a sum of $d$ odd primes. In $G_{n,d}$, two such monomials
give the same class precisely when a sequence of binary moves joins them.

Our main results concern interpolation and the growth of relations
among these representations. Equal pair sums force an affine
interpolation law for maps from the odd primes to any abelian group
(Theorem~\ref{main:prime-interpolation}), settling K\'atai and Phong's
conjecture on shifted prime arguments \cite[Conj.~6]{KP15Problems}.
Explicit Koszul cycles obtained by shifting primes span a positive
fraction of the linear syzygies in each sufficiently large admissible
Goldbach degree, and we determine the leading Betti asymptotic
(Theorem~\ref{main:weightwise}). Finally, the two natural finite
truncations have different leading constants for their total linear
strands (Theorem~\ref{main:cutoff-growth}).

The moves above already occur as edges of a strand of the minimal free
resolution of $G$ over $R$. In $F_{23,3}$ the vertices $x_3x_7x_{13}$,
$x_5^2x_{13}$, and $x_5x_7x_{11}$ form a triangle. Its oriented boundary
is a relation among the three edge generators, giving a minimal syzygy
(Example~\ref{ex:original23}). How does the number of such generators
grow with the Goldbach degree at a fixed monomial degree? Write
\[
 \beta_{i,n,d}=\dim_\kk\Tor_i^R(G,\kk)_{n,d},\qquad
 T_d(n)=\#\{(p_1,\ldots,p_d):p_j\text{ odd prime},\
                                  p_1+\cdots+p_d=n\}.
\]
A star in an internal degree means that we sum over that degree. Thus
$\beta_{i,*,d}=\sum_n\beta_{i,n,d}$ and
$\beta_i=\beta_{i,*,*}=\sum_{n,d}\beta_{i,n,d}$.
The count $T_d(n)$ includes repetitions. Write $\A$ for the odd primes
and $\A_t=\A\cap[1,t]$. We distinguish two growth questions. In the
full ring, we fix $d$ and let the Goldbach degree $n\to\infty$,
counting one bidegree at a time. For finite truncations on $\A_m$,
we fix $i$, sum over all Goldbach degrees in monomial degree $i+1$,
and let the cutoff $m\to\infty$.

For the first question, common shifts of primes give explicit Koszul
cycles in the linear strand. For example, adding two to each of $3,11,17$ gives
the primes $5,13,19$. In the Koszul complex on the variables this
common shift gives the cycle
\[
 x_5\varepsilon_{11}\varepsilon_{17}
 -x_{13}\varepsilon_3\varepsilon_{17}
 +x_{19}\varepsilon_3\varepsilon_{11}.
\]
Its differential vanishes by the binary relations in $G$. Up to sign,
this cycle represents the upper quadrilateral in the drawing of
$F_{33,3}$ in Example~\ref{ex:weight33}.
More generally, let $W^+_{n,d}$ be the span of the classes obtained
from $d$ distinct primes $p_1,\ldots,p_d$ and a shift $h>0$ for which
all $p_j+h$ are prime and $h+\sum_jp_j=n$. These classes have
homological degree $d-1$. The following estimates compare their span
with the entire linear strand.

\begin{maintheorem}[Growth in Goldbach degree and common shifts]
\label{main:weightwise}
Let $G$ be the binary Goldbach ring on all odd primes. Fix $d\ge3$.
Over every field,
\[
 \beta_{d-1,n,d}
 =\frac{d-1}{d!}\,T_d(n)+O_d\bigl((1+|\A_n|)^{d-2}\bigr).
\]
There are constants
$c_d,C_d>0$ independent of the field such that, for every sufficiently
large $n\equiv d\pmod2$,
\[
 c_d\,\frac{n^{d-1}}{(\log n)^d}
 \ \le\ \dim_\kk W^+_{n,d}
 \ \le\ \beta_{d-1,n,d}
 \ \le\ C_d\,\frac{n^{d-1}}{(\log n)^d}.
\]
\end{maintheorem}

In particular, for fixed $d\ge3$, as $n\to\infty$ through
$n\equiv d\pmod2$,
\[
 \beta_{d-1,n,d}\sim\frac{d-1}{d!}T_d(n)
 \sim\frac{\mathfrak S_d(n)}{d!(d-2)!}
             \frac{n^{d-1}}{(\log n)^d},
\]
where $\mathfrak S_d(n)$ is the prime singular series defined in
Corollary~\ref{cor:prime-betti-singular-series}. The estimates are
uniform over coefficient fields, and the translated cycles account
for a positive fraction of the linear strand in each such weight.
The coefficient $(d-1)/d!$ has a homological explanation.
For a partition into $d$ distinct primes there are $d$ ways to leave
one prime in the coefficient and put the other $d-1$ in the exterior
product. One incoming Koszul boundary relates these $d$ terms, leaving
$d-1$ dimensions before the outgoing differential is applied. Each
such partition has $d!$ orderings in $T_d(n)$. The error estimate
accounts for partitions with repeated primes and for the rank of the
outgoing differential. For $2\le k\le d$, if the defining moves replace
$k$ summands at a time, the same argument gives the
coefficient $\binom{d-1}{k-1}/d!$ in homological degree $d-k+1$
(Theorem~\ref{thm:fixed-weight-betti-euler}).

Theorem~\ref{thm:fixed-weight-betti-euler} proves the algebraic estimate
for arbitrary positive supports; Corollary~\ref{cor:prime-betti-singular-series}
combines it with known prime asymptotics. To account for a positive
fraction by explicit classes, Theorem~\ref{thm:prime-shift-proportion}
selects independent translated cycles by their exterior coordinates
and counts simultaneous prime completions using the linear forms
theorem of Tao and Ter\"av\"ainen \cite[Thm.~1.6]{TT}.
Related
translated cycles for lattice polygons appear in the work of
Castryck--Cools--Demeyer--Lemmens \cite[Secs.~4.1--4.2]{CCDL19}.
For a single shift, Proposition~\ref{thm:scroll-shift} in
Appendix~\ref{app:scroll-shift} uses the full scroll resolution to count
the resulting linear syzygies by Goldbach degree.

For the second growth question, $m$ bounds the available prime labels.
The variable cutoff $G^{<m}$ retains all binary relations on $\A_m$;
the relation-weight cutoff $G^{\le m}$ retains only those of common
weight at most $m$ on the same support. The superscripts distinguish
the ideals, while both support cutoffs are inclusive. Write
$\beta^{<m}$ and $\beta^{\le m}$ for their Betti numbers over
$R^{<m}=\kk[x_p:p\in\A_m]$.

\begin{maintheorem}[Linear-strand growth under truncation]
\label{main:cutoff-growth}
For every fixed $i\ge1$, as $m\to\infty$,
\[
 \beta^{<m}_{i,*,i+1}
       \sim\frac{i}{(i+1)!}|\A_m|^{i+1},\qquad
 \beta^{\le m}_{i,*,i+1}
       \sim\frac{1-2^{-i}}{(i+1)!}|\A_m|^{i+1},
\]
uniformly over coefficient fields.
\end{maintheorem}

Here $|\A_m|\sim m/\log m$, and each total sums over all Goldbach
degrees $n$ while keeping monomial degree $i+1$ fixed. For cubic
second syzygies the leading constants are $1/3$ and $1/8$.
Theorem~\ref{grow:both-cutoffs} proves these formulas and their
analogues for every fixed arity. Their proof counts across the finite
support and accounts for which relation weights are retained.
In contrast, Theorem~\ref{main:weightwise} concerns one bidegree of
the full ring. Stabilization identifies each fixed bidegree with its
finite-cutoff value once $m$ is large enough
(Proposition~\ref{setup:stabilization}).

The same binary relations also determine which functions on the odd
primes preserve equal pair sums.

\begin{maintheorem}[Interpolation on the odd primes]\label{main:prime-interpolation}
Let $\Podd$ be the set of odd primes, let $H$ be an abelian group,
and let $g:\Podd\to H$ satisfy
$g(p)+g(q)=g(p')+g(q')$ whenever $p+q=p'+q'$. Then for every odd prime $p$,
\[
 g(p)=g(3)+\frac{p-3}{2}\bigl(g(5)-g(3)\bigr).
\]
\end{maintheorem}

We prove this by computing the universal Freiman group generated by the
odd primes modulo their equal-pair-sum relations
(Theorem~\ref{depth:prime-pair-group}). It is free of rank two.
The proof adjoins auxiliary prime labels to realize a prescribed
equal-sum identity by binary moves, and then cancels those auxiliary
classes in the group. This converts the problem into a system of
linear forms to which the prime-counting theorem applies.
To our knowledge, this explicit integral identification has not
previously been recorded. It determines the group itself, including
the absence of torsion, and gives interpolation in every abelian target.

The complex-valued case proves the odd-prime affine-classification
assertion in Conjecture~1 of K\'atai--Phong \cite[p.~469]{KP15}, who prove affine classification
for functions preserving sums of at least three primes and propose
a binary analogue. Park's survey \cite[p.~6, after Thm.~3.7]{Park25}
still records the binary analogue as a conjecture.
Theorem~\ref{main:prime-interpolation} extends this odd-prime
classification to arbitrary abelian targets. The restriction to odd
primes matters: the printed all-primes formulation leaves an
independent value at $2$ under its pair-sum hypothesis, as explained
after Theorem~\ref{depth:prime-pair-group}.
Applying our interpolation formula to shifted arguments also proves
the distinct Conjecture~6 of K\'atai--Phong
\cite[p.~235]{KP15Problems}, as shown in Corollary~\ref{depth:KP-interpolation}.

Translation classes also occur in arbitrarily high homological
degrees. Unbounded multiplicities of positive prime differences imply
$\pd_R G=\infty$ (Corollary~\ref{cor:density}); arbitrarily long prime
progressions give another route through Green--Tao \cite{GreenTao08}.

Goldbach's conjecture has entered commutative algebra before:
Coykendall and Spicer \cite{CS12} show that a variant of it yields
Cohen--Kaplansky domains with any prescribed number of nonprime
irreducible elements. More broadly, the connection between additive
structure and commutative algebra has
been developed through Hilbert functions of iterated sumsets
\cite{EliahouMazumdar22}, projective monomial curves \cite{Elias22},
and the interpretation of Freiman's doubling inequality through
monomial ideals and their fiber cones \cite{HerzogHibiZhu}.
The connectivity of monomials under chosen moves also fits the
binomial-ideal framework of Diaconis--Eisenbud--Sturmfels \cite{DES98}.
To distinguish our ring from the full toric algebra, let
$J=\ker(R\to\kk[u,v])$, where $x_p\mapsto u^pv$.
The full toric quotient $R/J$ identifies all monomials of the same
bidegree, whereas the binary Goldbach ring is
$G=R/\langle J_{*,2}\rangle$, with $J_{*,2}=\bigoplus_n J_{n,2}$.
Thus $G$ imposes the quadratic equal-sum relations and their ideal
consequences, retaining the components under pair replacements.
The natural surjection $G\to R/J$ is an isomorphism precisely when
every nonempty binary-move fiber is connected.
For the full odd-prime support, we ask whether this is always so, or
equivalently whether $G$ is a domain
(Question~\ref{question:full-support-connectivity}).
For general supports connectivity can fail. We study how the fiber
components and the prime representation counts $T_d(n)$ enter the
minimal resolution.

The group computation also controls the finite variable-cutoff
resolutions. For all sufficiently large $m$, the rank-two group gives
\[
 \depth G^{<m}\le2,\qquad
 |\A_m|-2\le\pd_{R^{<m}}G^{<m}\le|\A_m|.
\]
An exhaustive prime-pair verification also gives
$\depth G^{<m}\le2$ at every integer cutoff
$3\le m\le20\,000\,000$, over every field
(Corollary~\ref{depth:verified-range}). This finite statement uses
triangular relations and localization, without computing the minimal
resolution; it does not determine the exact depth.
Theorem~\ref{depth:relation-cutoff-pd} gives an asymptotically maximal
projective dimension for the relation-weight cutoff as well.
Group completion identifies monomials of the same bidegree after
inverting the variables; it does not itself allow cancellation in $G$
or prove connectivity of every binary-move fiber.

The Betti estimates also have a Hilbert-series consequence. Put
\[
 \gamma_{n,d}=\sum_i(-1)^i\beta_{i,n,d},\qquad
 E_G(u,v)=\sum_{n,d}\gamma_{n,d}u^nv^d.
\]
If $H_G$ and $H_R$ denote the bigraded Hilbert series, the Euler
identity gives $E_G=H_G/H_R$. Corollary~\ref{cor:hilbert-euler-comparison}
shows that $(-1)^{d-1}\gamma_{n,d}$ differs from the dominant binary
Betti number by a lower-order term for fixed $d\ge3$ as
$n\to\infty$ through $n\equiv d\pmod2$. Thus the prime
representation asymptotic also governs coefficients determined by the
Hilbert series. The component basis makes $H_G$, and hence $E_G$,
independent of the field, although individual Betti numbers can vary.

The companion paper \cite{CompanionTopology} studies residue-field
Betti numbers, integral torsion, and higher operations in these rings.
Here the emphasis is on the arithmetic estimates and their consequences
for interpolation, depth, and projective dimension.
Throughout, $\kk$ denotes a field of characteristic $\ell$, where
$\ell=0$ or $\ell$ is prime. We reserve $k$ for the arity of a
defining additive relation. The variable $n$ is the Goldbach degree, $d$ the monomial degree, $i$ the
homological degree, and $m$ the truncation parameter.

Sections~\ref{sec:setup}--\ref{sec:koszul} establish the algebraic and
geometric framework. Section~\ref{sec:betti-euler-shifts} proves the
fixed-degree Betti estimates and the positive proportion of shifted
cycles as $n\to\infty$, then records their Hilbert-series consequence.
Section~\ref{sec:depth} treats interpolation and depth, and
Section~\ref{sec:growth} proves the truncation estimates as $m\to\infty$.
Section~\ref{sec:questions} concludes with further questions.
Appendix~\ref{app:scroll-shift} gives the scroll bound for a fixed
shift, and Appendix~\ref{app:saturation} studies saturation, quotients
by a prime variable, and dual Koszul cocycles.

\section{Goldbach rings}
\label{sec:setup}

Let $\A$ be a nonempty subset of the integers bounded below, and set
\[
 R=R_\A=\kk[x_a\mid a\in\A].
\]
We write $R$, $I$, and $G$ when the support $\A$ is understood.
Our main example is the set of all positive odd primes. As throughout
the paper, $\kk$ is a field of characteristic $\ell$, with $\ell=0$
or $\ell$ prime. The estimates below hold over every field; individual
Betti numbers can nevertheless depend on the characteristic. For example,
for the full binary odd-prime ring,
$(\beta_{5,135,7},\beta_{6,135,7})$ equals $(81,9)$ outside
characteristic two and $(82,10)$ in characteristic two
\cite[Theorem~6.3]{CompanionTopology}; this example is not used in the
proofs below. Integral forms over $\mathbb Z$ carry an
explicit subscript $\mathbb Z$.

\begin{definition}
For $k\ge2$, the \emph{$k$-ary Goldbach ring} supported on $\A$ is
the $\kk$-algebra $G=G_\A=R/I$, where
\begin{equation}\label{setup:relation-ideal}
 I=I_\A=
 \bigl\langle x_{a_1}\cdots x_{a_k}-x_{b_1}\cdots x_{b_k}
 \mid a_1+\cdots+a_k=b_1+\cdots+b_k\bigr\rangle .
\end{equation}
All the labels $a_i,b_i$ belong to $\A$.
\end{definition}
When $k=2$ we call $G$ a binary Goldbach ring, and when $k=3$ a
ternary Goldbach ring. The arity $k$ specifies how many factors a
single relation may replace.
Multiplying such a relation by a monomial keeps all the other factors
fixed. Thus, in monomial degree $d\ge k$, the quotient remembers which
representations can be connected by replacements of $k$ factors at a
time. For $d<k$ no defining relation is present, while in degree $k$
all monomials with the same sum are identified. This distinguishes a
fixed-arity Goldbach ring from the toric quotient, which identifies
equal sums in every monomial degree.

We keep $k$ for the arity, $n$ for the Goldbach degree, and $d$ for the
monomial degree throughout. Statements involving binary moves will
explicitly concern $k=2$; the definitions and grading constructions
below apply to every fixed $k\ge2$.

The labels need only be bounded below. The next observation lets us
pass to positive labels when a proof needs a positive grading, without
changing the ring or its homological invariants. It also records how
the original additive weights change under that passage.

\begin{proposition}[Affine invariance]
\label{setup:affine-invariance}
Fix integers $\nu>0$ and $N$, and put
\[
 \A'=\nu\A+N=\{\nu a+N\mid a\in\A\}.
\]
For each fixed $k\ge2$, the bijection $a\mapsto \nu a+N$ induces
isomorphisms
\[
 R_\A\cong R_{\A'},\qquad G_\A\cong G_{\A'} .
\]
Under the bigrading $\deg x_a=(a,1)$ introduced below, bidegree
$(n,d)$ becomes $(\nu n+Nd,d)$.
\end{proposition}
\begin{proof}
Write $a'=\nu a+N$ and send $x_a$ to $x_{a'}$. If
$\sum_{j=1}^k a_j=\sum_{j=1}^k b_j$, then
\[
 \sum_{j=1}^k a'_j
 =\nu\sum_{j=1}^k a_j+kN
 =\nu\sum_{j=1}^k b_j+kN
 =\sum_{j=1}^k b'_j .
\]
Thus the isomorphism $R_\A\to R_{\A'}$ sends $I_\A$ into $I_{\A'}$.
Its inverse, induced by $a'\mapsto(a'-N)/\nu$, sends $I_{\A'}$
into $I_\A$. A monomial with $d$ factors and sum of labels $n$
acquires sum $\nu n+Nd$.
\end{proof}

In particular, every support bounded below is affinely equivalent
to a support contained in the positive integers. Whenever we form
support cutoffs, we first make this translation and work with positive
labels. For every integer $t$, write
\[
 \A_t=\A\cap[1,t].
\]
Thus $\A_t$ includes its upper endpoint and is empty if $t<1$.
We use $n$ for Goldbach degree and $m$ for a truncation parameter:
$\A_n$ and $\A_m$ apply the same cutoff operation, with the index
indicating its role in the argument. In particular, they are equal
when $n=m$.

\begin{remark}[Monomial notation]
For a finite multiset $\boldsymbol a=\{a_1,\ldots,a_d\}$ in $\A$, write
$x_{\boldsymbol a}=x_{a_1}\cdots x_{a_d}$; repetitions are allowed and
the empty multiset represents $1$. We use parentheses
$\boldsymbol a=(a_1,\ldots,a_d)$ when the order matters. For a finite
ordered subset $\sigma=\{a_1<\cdots<a_d\}$, put
$x_\sigma=\prod_{a\in\sigma}x_a$.

Equivalently, let $\mathcal F$ consist of the functions
$\alpha:\A\to\mathbb Z_{\ge0}$ with finite support
$\operatorname{supp}\alpha=\{a:\alpha(a)\ne0\}$. Then
\[
 x^\alpha=\prod_{a\in\operatorname{supp}\alpha}x_a^{\alpha(a)}
\]
ranges over the monomials of $R$. For instance, $\alpha(3)=2$ and
$\alpha(11)=4$, with all other values zero, represents $x_3^2x_{11}^4$.
Every polynomial is a finite sum of these monomials with coefficients
in $\kk$, uniquely after like terms are combined. These two notations
let us pass between additive representations and algebraic expressions.
\end{remark}

\subsection{Bigraded algebra structures on \texorpdfstring{$R$}{R} and \texorpdfstring{$G$}{G}}
There is a natural bigraded $\kk$-algebra structure underlying our setup
which organizes the algebraic and homological objects we study. For
$x_{\boldsymbol a}=x_{a_1}\cdots x_{a_d}$, set
\[
 \deg_1(x_{\boldsymbol a})=\sum_{j=1}^d a_j,
 \qquad \deg_2(x_{\boldsymbol a})=d.
\]
We call the first grading the \emph{Goldbach degree}, or additive weight,
and the second the \emph{monomial degree}. Homological degree is
separate from these two internal gradings. For the extension to Koszul
terms, see \hyperref[setup:koszul-grading]{the degree convention in
Section~\ref*{sec:koszul}}.

Let $R_{n,d}$ be the $\kk$-span of the monomials of bidegree $(n,d)$,
and denote their set by $\mathcal M_{n,d}$. Then
\[
 R=\bigoplus_{n\in\mathbb Z,\ d\ge0}R_{n,d},\qquad
 R_{n,d}R_{n',d'}\subseteq R_{n+n',d+d'}.
\]
Each piece is finite dimensional. Indeed, writing $b=\min\A$, every
label in a monomial of degree $d>0$ lies between $b$ and
$n-(d-1)b$.
The dimension of $R_{n,d}$ counts the unordered expressions of $n$ as
a sum of $d$ elements of $\A$, with repetitions permitted. For positive
supports these are the corresponding unordered partitions. In monomial degree
zero, $R_{0,0}=\kk$ and all other pieces vanish.

The generators of $I$ are bihomogeneous, so the quotient
inherits a bigraded algebra structure
\[
 G=\bigoplus_{n,d}G_{n,d},\qquad
 G_{n,d}G_{n',d'}\subseteq G_{n+n',d+d'}.
\]

\begin{example}[Goldbach degree $23$]
\label{setup:weight23}
Take $\A$ to be the positive odd primes and let $G$ be the binary
Goldbach ring. The nonzero pieces of $R$ in weight $23$ are
\begin{align*}
 R_{23,1}&=\kk x_{23},\\
 R_{23,3}&=\kk x_3^2x_{17}\oplus\kk x_3x_7x_{13}
          \oplus\kk x_5^2x_{13}\oplus\kk x_5x_7x_{11},\\
 R_{23,5}&=\kk x_3^4x_{11}\oplus\kk x_3^3x_7^2
          \oplus\kk x_3^2x_5^2x_7\oplus\kk x_3x_5^4,\\
 R_{23,7}&=\kk x_3^6x_5.
\end{align*}
For instance, the monomials in $R_{23,5}$ correspond to the four
unordered partitions
\[
 (3,3,3,3,11),\quad(3,3,3,7,7),\quad
 (3,3,5,5,7),\quad(3,5,5,5,5).
\]
In $G$, the four cubic monomials represent the same class, since
\begin{equation}\label{setup:weight23-cubic-moves}
 x_3^2x_{17}=x_3x_7x_{13}=x_5^2x_{13}=x_5x_7x_{11}.
\end{equation}
Each step applies one binary Goldbach move, using the pair sums
$3+17=7+13$, $3+7=5+5$, and $5+13=7+11$ in order. Similarly,
\[
 x_3^4x_{11}=x_3^3x_7^2=x_3^2x_5^2x_7=x_3x_5^4,
\]
using $3+11=7+7$ and then $3+7=5+5$ twice. Thus these two fibers
are connected by binary moves. Their common classes are nonzero,
since evaluation at $x_p=1$ descends to $G$ and sends every monomial
to $1$. Together with the singleton fibers in monomial degrees one and seven,
this proves
\[
 \dim_\kk G_{23,d}=1\quad(d=1,3,5,7).
\]
All other monomial-degree components in weight $23$ vanish.
\end{example}

For the binary ring on the odd primes, the same method works in larger
fibers. For example, in $G_{36,4}$,
\[
 x_3^2x_{11}x_{19}
 =x_3x_5x_{11}x_{17}
 =x_5x_7x_{11}x_{13}
 =x_5^2x_{13}^2 .
\]
The three moves use $3+19=5+17$, $3+17=7+13$, and $7+11=5+13$.
To count all the classes in a given bidegree we must determine which
partitions can be joined by such paths.

\subsection{The associated toric ideal}
Let $S$ be the $\kk$-algebra
\[
 S=\kk[u^a v\mid a\in\A]\subseteq\kk[u,u^{-1},v].
\]
For a positive support we can take the ambient ring to be $\kk[u,v]$.
Giving $u$ degree $(1,0)$ and $v$ degree $(0,1)$ makes
\[
 \phi:R\twoheadrightarrow S,\qquad x_a\longmapsto u^a v,
\]
a surjective bigraded $\kk$-algebra map. Its kernel $J$ is the
associated toric ideal. For the $k$-ary ideal $I$, each defining generator maps to zero:
its two terms have image $u^n v^k$, where $n$ is their common sum.
Thus $I\subseteq J$ for every fixed arity, and $\phi$ induces a
surjection $\overline\phi:G\to S$. In each represented bidegree
$(n,d)$, the target has the single basis monomial $u^n v^d$.

The map $\overline\phi$ forgets how representations are connected:
it retains only their common sum and number of factors. Its kernel
therefore measures the difference between allowing replacements of
exactly $k$ factors and allowing equal-sum relations of every arity.

To describe that kernel, take the monomials in $\mathcal M_{n,d}$ as
vertices. Join two distinct vertices when they have the form
$x^\alpha x_{a_1}\cdots x_{a_k}$ and
$x^\alpha x_{b_1}\cdots x_{b_k}$ with
$\sum_{j=1}^k a_j=\sum_{j=1}^k b_j$. The monomial $x^\alpha$ is
held fixed. Write $\Gamma_{n,d}$ for this graph, with the arity
understood from the specified ring. In the binary case $k=2$, its
edges are precisely the pair replacements used in the weight-$23$
example.

\begin{proposition}[Component bases]
\label{setup:components}
For a fixed arity $k\ge2$, one monomial from each connected
component of $\Gamma_{n,d}$ gives a basis of $G_{n,d}$ over $\kk$.
\end{proposition}
\begin{proof}
In bidegree $(n,d)$ the ideal is spanned by differences of the
endpoints of edges, the monomial multiples of the defining
binomials. Summing along a path gives the difference of its endpoints.
The map from the free module on vertices to the free module on
components sends each monomial to its component. Its kernel is
spanned by these differences, giving the asserted basis. The same
proof works for any $k$.
\end{proof}

Consequently $I=J$ if and only if each nonempty fiber is connected,
equivalently if and only if $\dim_\kk G_{n,d}\le1$ in every
bidegree. In that case $\overline\phi$ identifies $G$ with the
integral domain $S$. For the binary ring on the full odd-prime support, the converse
also follows from Theorem~\ref{depth:prime-pair-group}. Goldbach degree and
monomial degree identify the relation group with the subgroup generated by
$(p,1)$, so $\overline\phi$ becomes an isomorphism after the variables
are inverted. If $G$ is a domain it embeds in that
localization, since each monomial has nonzero image in $S$.
For general supports the converse can fail. For the binary ring on $\A=\{3,5,11\}$
the pair sums are distinct, hence $I=0$, although
$x_3^3x_{11}-x_5^4\in J$.

\begin{example}
Take the binary ring on the odd primes $\A$ congruent to $1$ modulo $4$. Then
$x_5^2x_{29}$ and $x_{13}^3$ both map to $u^{39}v^3$, but no
binary move connects them. Indeed $13+13$ is the only pair of
labels in $\A$ with sum $26$, so $x_{13}^3$ is an isolated vertex.
The other two triples of weight $39$ are joined by
$x_5^2x_{29}=x_5x_{17}^2$ in $G$. Thus $\dim_\kk G_{39,3}=2$.
\end{example}

Every nonempty binary fiber of monomial degree two is connected, since one
move joins any two decompositions of the pair sum. For a binary ring on a finite
support $\A$, write $G_{*,d}=\bigoplus_n G_{n,d}$, where $*$ denotes
the direct sum over all Goldbach degrees. Then
\begin{equation}\label{setup:low-hilbert}
 \dim_\kk G_{*,1}=|\A|,\qquad \dim_\kk G_{*,2}=|\A+\A|.
\end{equation}
For the toric quotient, write $S_{*,d}=\bigoplus_n S_{n,d}$.
Its basis consists of the monomials $u^n v^d$ with $n\in d\A$, where
$d\A$ denotes the $d$-fold sumset. Thus, for every finite support,
\[
 \dim_\kk S_{*,d}=|d\A|.
\]
This identifies iterated sumset cardinalities with the Hilbert function
of a standard graded algebra, the correspondence used by
Eliahou--Mazumdar \cite{EliahouMazumdar22} and Elias \cite{Elias22}.

\paragraph{The curve defined by the toric ideal.}
The toric quotient provides a geometric comparison for the Goldbach
ring. Its monomials are indexed by attained bidegrees rather than by
move components, so its Hilbert function is simpler. This comparison
does not identify $G$ with $S$: connectivity, cancellation, and depth
must still be checked in the relevant ring. On the odd-prime support
the toric quotient has a particularly explicit curve description,
independent of which arity was used to define $G$.

For the odd primes write $y=u^3v$, $t=u^2$ and $e_p=(p-3)/2$.
Then $S=\kk[yt^{e_p}:p\in\Podd]$. Write $S^{<m}$ for the subalgebra
generated by the $x_p=yt^{e_p}$ with $p\in\A_m$. For $m\ge5$,
the first two generators are $y$ and $yt$, so the fraction field
of $S^{<m}$ is $\kk(y,t)$. Thus $\operatorname{Spec}S^{<m}$ is a
two-dimensional affine cone. With the monomial grading its projective
curve is the closure of $t\mapsto[1:t:\cdots:t^{e_P}]$, where $P$
is the largest prime at most $m$ and the omitted coordinates are the
other $t^{e_p}$ in increasing order. This rational monomial curve is
locally Cohen--Macaulay, since its local rings at closed points are
one-dimensional Noetherian domains.

The finite exponent sets $\{e_p:p\in\A_m\}$ place these curves in the
sumset framework of Elias \cite{Elias22}. Relations between sumset
structure, Castelnuovo--Mumford regularity, and the Cohen--Macaulay
property of their homogeneous coordinate rings are developed by
Gimenez--Gonz\'alez-S\'anchez \cite{GimenezGonzalez23} and
Hoa--Tien \cite{HoaTien26}. The recent preprint
\cite{ChaudharyKhandelwal26} further studies the difference between
sumset regularity and Castelnuovo--Mumford regularity.

The cone can have depth one. At $m=42$, $x_3,x_{41}$ are homogeneous
parameters since $x_p^{19}=x_3^{19-e_p}x_{41}^{e_p}$ in $S^{<42}$.
The class of $x_5x_7$ modulo $x_3$ is nonzero: dividing by $x_3$ would
require a coefficient of monomial degree one and Goldbach degree $9$.
Yet $x_{41}x_5x_7=x_3x_{13}x_{37}$, so this pair of parameters is
not a regular sequence. The domain $S^{<42}$ therefore has depth
one at the vertex over every field.

On the full support $S$ still has dimension two. Any finite strict
chain of prime ideals contracts strictly to a sufficiently large
finite cutoff: include one witness to each strict inclusion. For the
lower bound, the ideal generated by the $x_p$ with $p>3$ has quotient
$\kk[x_3]$ and lies strictly between $(0)$ and $S_+=\bigoplus_{d>0}S_{*,d}$.
For $q>p>3$ we have $x_p^{e_q}=x_3^{e_q-e_p}x_q^{e_p}$, giving
$\sqrt{x_3 S}=S_+$. Thus $D(x_3)$ is the complement of the vertex,
and $S[x_3^{-1}]=\kk[y^{\pm1},t]$. Taking its monomial-degree-zero part gives
$\operatorname{Proj}S=D_+(x_3)\cong\mathbb A^1_\kk$.
At the vertex $S_+/S_+^2$ has one basis element for every odd prime.
The local ring there is non-Noetherian, and the maximum length of a
regular sequence in its maximal ideal is one. Indeed, $x_3$ is regular.
For any two nonzero nonunits $f,g$ in that local ring, some powers
satisfy $f^N=x_3 a$ and $g^M=x_3 b$. Then $a\notin(f^N)$ but
$a g^M=b f^N$, ruling out a regular sequence $f,g$ since powers of
a regular sequence are regular.

Finally, affine invariance supplies a positive grading without changing
the original bidegrees. Translating every label by $1-b$ gives the
auxiliary weight
\begin{equation}\label{setup:positive-weight}
 w(n,d)=n-(b-1)d,\qquad w(x_a)=a-b+1\ge1.
\end{equation}
There are finitely many variables, and consequently finitely many
monomials, in every bounded auxiliary-weight range. This connected,
locally finite positive grading will allow us to construct minimal
resolutions even when $\A$ is infinite.

\section{The minimal free resolution of \texorpdfstring{$G$}{G} over \texorpdfstring{$R$}{R}}
\label{setup:resolution-construction}
Throughout this section $G=R/I$ is the binary Goldbach ring on
$\A$ over the field $\kk$.
Our goal is to construct its bigraded minimal free resolution $(F,\partial)$
over $R$, where by minimal we mean that $F_0=R$ and
$\partial F\subseteq\m F$ for $\m=\langle x_a\mid a\in\A\rangle$.
Since $I$ need not be finitely generated, we first construct bigraded
minimal free resolutions of truncated versions of $G$ over polynomial
subrings of $R$ with finitely many variables. These resolutions will
give rise to a directed system over $R$ whose direct limit is the
bigraded minimal free resolution we seek.

A bigraded differential preserves both degrees. We can therefore consider
one bidegree at a time and write
\begin{equation}\label{setup:strand-decomposition}
 F=\bigoplus_{n\in\mathbb Z,\ d\ge0}F_{n,d},
 \qquad F_{i,n,d}:=(F_{n,d})_i=(F_i)_{n,d}.
\end{equation}
For any bigraded complex $C$ we use the same convention: $C_{n,d}$
is the entire fixed-bidegree complex, and $C_{i,n,d}$ is its term
in homological degree $i$.
Each $F_{n,d}$ is a complex of $\kk$-vector spaces, and multiplication
by $x_a$ sends it into $F_{n+a,d+1}$. We will obtain $F_{n,d}$ as the
eventual value of the corresponding finite-stage complexes. After
the construction, we describe its first differentials by graphs of
binary moves and exhibit polygonal and polyhedral strands. These
examples show concretely what stabilizes at each cutoff.

\subsection{First truncation}
For the cutoff constructions, first apply the affine translation of
Proposition~\ref{setup:affine-invariance} to make the support positive.
We write $\A$ for this translated support and use its Goldbach degrees
below. For each positive integer $m$, put $\A_m=\A\cap[1,m]$. Define
\[
 \begin{aligned}
 R^{<m}&=\kk[x_a\mid a\in\A_m],\\
 I^{\le m}&=\bigl\langle x_{a_1}x_{a_2}-x_{b_1}x_{b_2}\mid
 a_1+a_2=b_1+b_2\le m,\ a_1,a_2,b_1,b_2\in\A_m\bigr\rangle .
 \end{aligned}
\]
Set $G^{\le m}=R^{<m}/I^{\le m}$. The superscript $\le m$ records the
bound on the common sum in the defining relation. The ambient ring
$R^{<m}$ uses the inclusive support $\A_m$; the symbol $<$ distinguishes
the variable-cutoff notation rather than a strict endpoint. Positivity
forces each label in a defining relation of weight at most $m$ to be
less than $m$. When
we introduce the second truncation, we will bound only the individual
labels. Both $R^{<m}$ and $G^{\le m}$ inherit bigraded $\kk$-algebra
structures. Since every variable and every defining relation occurs at
some finite stage, we recover $R$ and $G$ as bigraded direct limits:
\[
 R=\varinjlim_m R^{<m},\qquad G=\varinjlim_m G^{\le m}.
\]
We also write $R_m=R^{<m}$ when several cutoffs occur together.
Let $F^{\le m}$ be the minimal bigraded free resolution of $G^{\le m}$
over $R^{<m}$. Its Betti number
\[
 \beta^{\le m}_{i,n,d}
   =\dim_\kk\Tor_i^{R^{<m}}(G^{\le m},\kk)_{n,d}
\]
counts the basis generators of homological degree $i$, Goldbach
degree $n$, and monomial degree $d$.

\subsubsection{The minimal free resolution as a direct limit}
We now transport everything to $R$. Put
\begin{equation}\label{setup:first-base-change}
 I^{(m)}=I^{\le m}R,\qquad G^{(m)}=R/I^{(m)},\qquad
 F^{(m)}=R\otimes_{R^{<m}} F^{\le m}.
\end{equation}
The inclusion $R^{<m}\hookrightarrow R$ makes $R$ a bigraded free
$R^{<m}$-module, with basis the monomials in the remaining variables.
Consequently $F^{(m)}$ is a finite-length minimal bigraded free
resolution of $G^{(m)}$ over $R$. The canonical maps
$G^{(m)}\to G^{(m+1)}$ lift to bigraded comparison maps of $R$-complexes
\begin{equation}\label{setup:first-comparisons}
 F^{(m)}\longrightarrow F^{(m+1)}\longrightarrow F^{(m+2)}
 \longrightarrow\cdots.
\end{equation}
We claim that we can take $F$ to be the bigraded direct limit of this
system. To see why the direct limit is free, we choose the comparisons
so that in each fixed bidegree $(n,d)$ they are eventually the identity
on the same $\kk$-complex.

Write $b=\min\A$ and let $w(n,d)=n-(b-1)d$ be the positive weight of
\eqref{setup:positive-weight}. We will also use $e(n,d)=n-bd$.
For a monomial $x_{a_1}\cdots x_{a_d}$ of bidegree $(n,d)$ this is
$\sum_j(a_j-b)\ge0$. Thus an omitted label $a\ge m+1$ already accounts for
at least $m+1-b$ of this sum. This observation gives a cutoff that keeps
track of both $n$ and $d$.

\begin{lemma}[Adjoining a relation]\label{setup:mapping-cone}
Let $J\subsetneq R$ be a finitely generated bihomogeneous ideal and
let $P$ be a minimal bigraded free resolution of $R/J$. Suppose
$f\notin J$ is bihomogeneous of positive monomial degree and bidegree $\alpha$.
One can choose a minimal resolution of $R/(J,f)$ and a comparison
from $P$ which retains the old complex and its chosen basis in every
bidegree $\lambda$ satisfying
\[
 w(\lambda)<w(\alpha)\quad\text{or}\quad e(\lambda)<e(\alpha).
\]
\end{lemma}
\begin{proof}
Choose a finite polynomial subring containing $f$ and generators of
$J$. We may form the colon ideal and its resolution there and extend
to $R$ by freeness. The short exact sequence
\[
 0\longrightarrow R/(J:f)(-\alpha)
 \xrightarrow{\ f\ }R/J\longrightarrow R/(J,f)\longrightarrow0
\]
gives a comparison $L(-\alpha)\to P$, where $L$ resolves $R/(J:f)$.
The mapping cone $C$ is a free resolution of $R/(J,f)$, with
\begin{equation}\label{setup:cone-modules}
 C_i=P_i\oplus L_{i-1}(-\alpha).
\end{equation}
All basis degrees of a minimal $L$ lie in the semigroup of monomial
bidegrees: starting with $L_0=R$, each new generator is chosen in a
bidegree already represented by a monomial multiple of a preceding
generator. In particular both $w$ and $e$ are nonnegative on these
degrees. The second summand in \eqref{setup:cone-modules} therefore
has $w\ge w(\alpha)$ and $e\ge e(\alpha)$.

If $C$ is not minimal, we split off contractible pairs
$R(-\lambda)\xrightarrow{1}R(-\lambda)$. Initially a scalar entry
must involve the new summand because $P$ is minimal. Such a pair
has both $w(\lambda)\ge w(\alpha)$ and $e(\lambda)\ge e(\alpha)$.
Cancelling it leaves the indicated low bidegrees unchanged: no
homogeneous map can take a generator of lower $w$ or lower $e$
to a multiple of this pair, since coefficients have nonnegative
$w$ and $e$. In particular no subsequent cancellation starts in
those bidegrees. Projecting the inclusion $P\to C$ onto the remaining
minimal summand gives the desired comparison.
\end{proof}

\begin{theorem}[Stabilization of the first truncation]
\label{setup:first-stabilization}
For each fixed bidegree $(n,d)$, the maps
\[
 F^{(m)}_{n,d}\longrightarrow F^{(m+1)}_{n,d}
              \longrightarrow F^{(m+2)}_{n,d}\longrightarrow\cdots
\]
can be chosen to be the identity as soon as
\[
 w(n,d)<\gamma_m,\qquad
 \gamma_m=\min\{m+1-2(b-1),\,m-b+3\}.
\]
Another sufficient condition is
\begin{equation}\label{setup:first-bidegree-bound}
 m\ge n-b(d-2)\quad\text{and}\quad m\ge n-b(d-1).
\end{equation}
In particular, for positive supports the extended strands stabilize
whenever $m\ge n$.
\end{theorem}
\begin{proof}
Fix $(n,d)$ and consider the construction of $F^{(m+1)}$ starting
from $F^{(m)}$. If $I^{(m)}=I^{(m+1)}$ we take the identity.
Otherwise write $I^{(m+1)}=I^{(m)}+\langle g_1,\ldots,g_t\rangle$
and put $I^{(m),s}=I^{(m)}+\langle g_1,\ldots,g_s\rangle$.
We construct the comparison by lifting each quotient in turn:
\[
 F^{(m)}\longrightarrow F^{(m),1}\longrightarrow\cdots
 \longrightarrow F^{(m),t}=F^{(m+1)}.
\]
At step $s$ we apply Lemma~\ref{setup:mapping-cone} to $I^{(m),s}$
and $g_{s+1}$, skipping any redundant relation. We form and minimize
these cones over $R^{<m+1}$ and then extend scalars, so the resulting
$F^{(m+1)}$ has the form \eqref{setup:first-base-change}. Each new defining
binomial has common sum $m+1$ or involves the newly admitted label
$m+1$. In the former case its $w$-weight is $m+1-2(b-1)$, and in the
latter case its $w$-weight is at least $m-b+3$. This proves the first
bound. For $e$ the corresponding lower bounds are $m+1-2b$ and
$m+1-b$, respectively. Requiring $n-bd$ to be strictly below both
gives \eqref{setup:first-bidegree-bound}. Both bounds increase with
$m$, so the comparison retains that entire strand at every later
stage. If $b\ge1$, all new relations have Goldbach degree at least
$m+1$. The cone argument with Goldbach degree as the positive weight
therefore preserves all degrees $n\le m$.
\end{proof}

\begin{corollary}[The resolution as a direct limit]
\label{setup:minimal}
With these choices the complex
\begin{equation}\label{setup:resolution-limit}
 F=\varinjlim_m F^{(m)}
\end{equation}
is a minimal bigraded free resolution of $G$ over $R$. Its strand
$F_{n,d}$ is a finite $\kk$-complex for every $(n,d)$.
\end{corollary}
\begin{proof}
To see that $F$ is free over $R$, choose for every $w$-weight bound
$W$ a stage after which all basis generators of weights at most $W$
are retained. We may choose these stages increasingly with $W$.
The resulting sets of stable generators are nested. Their union spans
$F$, since each homogeneous element belongs to some stabilized weight
range. Suppose that distinct stable generators $e_1,\ldots,e_s$ and
elements $f_1,\ldots,f_s\in R$ satisfy
\[
 f_1e_1+\cdots+f_s e_s=0 .
\]
We may take the relation to be homogeneous of bidegree $(n,d)$.
It then occurs in a single stabilized range, where all the $e_j$ are
basis elements of one finite-stage free module $F_i^{(m)}$. Thus
$f_1=\cdots=f_s=0$. This proves freeness.

By \eqref{setup:strand-decomposition} we may compute homology on each
strand separately. Choose a sufficiently large $m=m(n,d)$ for each
bidegree. Then
\[
 \mathrm{H}_*(F)=\bigoplus_{n,d}\mathrm{H}_*(F_{n,d})
       =\bigoplus_{n,d}\mathrm{H}_*(F^{(m(n,d))}_{n,d})
       \cong\bigoplus_{n,d}G_{n,d}=G
\]
in homological degree zero, with zero homology in positive degrees.
Here the finite-stage complexes are resolutions and the quotients
$G^{(m)}$ have direct limit $G$. Minimality follows since the
differential of each stable generator lies in $\m F$ at the stage
where it was chosen. Finally, every bounded $w$-weight range is
finite dimensional at its stabilized stage. Since every such stage
uses a finite-length resolution, the stabilized strand is finite.
\end{proof}

The usual comparison argument gives uniqueness up to bigraded chain
isomorphism. The choices above identify the strands literally. Any
other bigraded comparison lifting the same quotient map differs from
ours by $\partial h+h\partial$. By minimality this difference has
coefficients in $\m$. On the stabilized basis in a bounded $w$-weight
range, the difference strictly lowers basis weight, so adding it
to the identity still gives an isomorphism on that range.

\subsection{Second truncation}
In the previous subsection we constructed $F$ as the direct limit of
the minimal free resolutions associated to the first truncation of $G$.
We now use all binary relations on the variables already present:
\[
 I^{<m}=\bigl\langle x_{a_1}x_{a_2}-x_{b_1}x_{b_2}
 \mid a_1+a_2=b_1+b_2,\ a_1,a_2,b_1,b_2\in\A_m\bigr\rangle .
\]
Set $G^{<m}=R^{<m}/I^{<m}$. Thus $I^{\le m}\subseteq I^{<m}$, and the
second truncation can include relations whose Goldbach degrees exceed
$m$. Just as before, $G=\varinjlim_m G^{<m}$. Let $F^{<m}$ be the
minimal bigraded free resolution of $G^{<m}$ over $R^{<m}$. We
transport everything to $R$ by putting
\[
 I^{[m]}=I^{<m}R,\qquad G^{[m]}=R/I^{[m]},\qquad
 F^{[m]}=R\otimes_{R^{<m}}F^{<m}.
\]
We again obtain a directed system of $R$-complexes and bigraded
comparison maps.

\begin{theorem}[Stabilization of the second truncation]
\label{setup:second-stabilization}
The minimal free resolution of $G$ over $R$ is also the direct limit
of the $F^{[m]}$, up to bigraded chain isomorphism. For every fixed
$(n,d)$ the comparisons
\[
 F^{[m]}_{n,d}\longrightarrow F^{[m+1]}_{n,d}
              \longrightarrow F^{[m+2]}_{n,d}\longrightarrow\cdots
\]
can be chosen to be the identity if $w(n,d)<m-b+3$, or if
\begin{equation}\label{setup:second-bidegree-bound}
 m\ge n-b(d-1).
\end{equation}
\end{theorem}
\begin{proof}
Every new defining binomial involves $m+1$. It therefore has
$w$-weight at least $m-b+3$ and $e$-weight at least $m+1-b$.
Adjoin these binomials one at a time as in the first truncation.
The mapping cone preserves $(n,d)$ if either weight is below its
respective bound. For $e$ this says exactly \eqref{setup:second-bidegree-bound}.
The stable basis and strandwise homology arguments in
Corollary~\ref{setup:minimal} give a minimal free resolution of $G$.
By uniqueness it is isomorphic to $F$.
\end{proof}

These are sufficient bounds on the cutoff for a given bidegree. They
also give bounds before extending scalars: an omitted variable has
$e$-weight at least $m+1-b$, so it cannot contribute to a strand satisfying
\eqref{setup:second-bidegree-bound}. Thus \eqref{setup:first-bidegree-bound}
suffices for $F^{\le m}_{n,d}\cong F_{n,d}$, and
\eqref{setup:second-bidegree-bound} suffices for $F^{<m}_{n,d}\cong F_{n,d}$.

In a nonempty fiber there is also an exact comparison of the two
unextended strands. The variable cutoff recovers $F_{n,d}$ as soon
as every monomial of $R_{n,d}$ is present. The relation cutoff
also requires $a+a'\le m$ whenever two distinct vertices
$x^\alpha x_ax_{a'}$ and $x^\alpha x_bx_{b'}$ in the fiber differ by the move
$a+a'=b+b'$.
Thus the variable cutoff always recovers the whole strand at least
as early as the relation cutoff.

Indeed, a missing vertex lowers the dimension in homological degree
zero. Once every vertex is present, a missing quadratic relation of
the indicated kind lowers the dimension in degree one. A fiber
of $h$ quadratic monomials contributes $h-1$ minimal
generators when its relations are included and zero when they are
omitted. Conversely, the cone in \eqref{setup:cone-modules} for a
new relation $f$ can affect $(n,d)$ only if $(n,d)-\deg f$ is a
monomial bidegree, since all shifts of the colon resolution belong
to that semigroup. A monomial $x^\alpha$ of this degree would put both
terms of $x^\alpha f$ in the fiber. Under the stated conditions every such
relation is already included. The cone and its cancellations
therefore preserve the strand, and extending scalars adds no
monomial multiple involving an omitted label in this bidegree.

For example, in the odd-prime fiber $(23,3)$ the largest label is
$17$ and the largest weight of such a quadratic move is $20$, from
$3+17=7+13$. The variable cutoff recovers the strand at $m=17$,
and the relation cutoff at $m=20$ (Example~\ref{ex:original23}).

For either truncation the number of basis generators in a fixed
homological degree and bidegree stabilizes to
\[
 \beta_{i,n,d}=\dim_\kk(F_i/\m F_i)_{n,d}
              =\dim_\kk\Tor_i^R(G,\kk)_{n,d}.
\]
Write $\beta^{<m}_{i,n,d}=\dim_\kk\Tor_i^{R^{<m}}(G^{<m},\kk)_{n,d}$.
Let $\sigma^{<}_{i,n,d}$ be the least positive integer $M$ such that
$\beta^{<m}_{i,n,d}=\beta_{i,n,d}$ for every $m\ge M$.
Define $\sigma^{\le}_{i,n,d}$ in the same way for the first truncation.
The cutoff bounds for $F_{n,d}$ apply to these indices, although an
individual Betti number can stabilize earlier.
The term $F_{i,n,d}$ also contains monomial multiples of basis
generators of lower bidegrees. Finally, define $\delta^{\le m}$ and
$\rho^{\le m}$ to be the depth and projective dimension of $G^{\le m}$
over $R^{<m}$, and similarly define $\delta^{<m}$ and $\rho^{<m}$ for
$G^{<m}$. Depth is taken at the homogeneous maximal ideal. By the
Auslander--Buchsbaum formula,
\begin{equation}\label{setup:AB-cutoffs}
 \rho^{\le m}+\delta^{\le m}=|\A_m|
                       =\rho^{<m}+\delta^{<m}.
\end{equation}
For the odd primes, $|\A_m|=\pi(m)-1$ for integers $m\ge3$,
where $\pi(y)$ is the usual prime-counting function, including $2$.
Section~\ref{sec:depth} uses this identity to estimate projective
dimension.

\begin{remark}
The mapping-cone argument applies to any proper bihomogeneous ideal
in $R$. Local finiteness lets us choose a homogeneous generating set
finite in each bounded $w$-weight range. Adjoining these generators in
order of weight gives the same stabilization. When $\A$ is positive,
the construction also works for an ideal homogeneous in Goldbach degree
alone, for instance
$\langle x_a-x_{a'}x_{a''}x_{a'''}\mid a=a'+a''+a'''\rangle$.
In that case the complexes decompose by $F_n$.
\end{remark}

\subsection{Minimal moves and the geometry of finite strands}
\label{sec:examples}
We now specialize the preceding construction to its combinatorial
description. A monomial is a vertex, a chosen binary relation gives
an edge, and a syzygy records a relation among those edges. The
following order fixes the quadratic generators and their orientations.
Higher cells will be specified by their boundaries in individual
strands; no general cellular resolution is assumed.
\subsubsection{The product order}
In each strand we compare monomials of the same bidegree $(n,d)$.
Write them as $x_{\boldsymbol a}=x_{a_1}\cdots x_{a_d}$ and
$x_{\boldsymbol b}=x_{b_1}\cdots x_{b_d}$, so that
$\sum_j a_j=\sum_j b_j=n$. On the odd primes our order is simply
\[
 x_{\boldsymbol a}\prec x_{\boldsymbol b}
 \quad\Longleftrightarrow\quad
 a_1\cdots a_d<b_1\cdots b_d.
\]
Unique factorization distinguishes all the monomials in the fiber.

The connection with \emph{majorization} explains why balanced pairs
give the largest products. List both partitions in decreasing order.
Then $\boldsymbol a$ majorizes $\boldsymbol b$ if
$\sum_{j=1}^s a_j\ge\sum_{j=1}^s b_j$ for each $1\le s<d$.
For distinct partitions of positive labels this implies
$a_1\cdots a_d<b_1\cdots b_d$, by strict concavity of $\log$ and
the majorization inequality \cite[Ch.~3]{MOA11}.
For $d=2$ the larger of the two parts determines majorization.
For $d=3$ the product also compares partitions that majorization
leaves incomparable. In $R_{31,3}$, for instance, the first two
partial sums for $(13,13,5)$ are $13,26$, whereas those for
$(17,7,7)$ are $17,24$. Their products are $845$ and $833$.

For an arbitrary bounded-below support $\A$, put
$c=\max\{0,2-\min\A\}$ and compare first by
$P_c(x^\alpha)=\prod_{a\in\A}(a+c)^{\alpha(a)}$.
The same majorization argument applies to $\log(a+c)$.
If two products agree, take the largest label $a_0$ at which the
exponents differ and put $x^\alpha\prec x^\beta$ if
$\alpha(a_0)<\beta(a_0)$. For example, $x_2x_6^2$ and $x_3^2x_8$
on $\A=\{2,3,6,8\}$ have the same sum, monomial degree and product.
The comparison by $\log P_c$ is a weight order with weights
$\log(a+c)$ \cite[Definition~1.2 and Section~6]{Vietri02}.

\begin{proposition}[Product monomial order]\label{moves:order}
The same comparison on all monomials is a monomial well-order on $R$.
With $c$ fixed, it restricts to each polynomial subring on a subset
of the variables.
\end{proposition}
\begin{proof}
The product comparison followed by the stated lexicographic tie-break
is a total order, with $1$ its least element. Products multiply, and
adding the same exponent vector preserves the largest differing label
and its exponent comparison. Hence multiplying by any monomial
preserves strict inequalities.

For $N\ge1$, only finitely many monomials satisfy $P_c\le N$.
Every factor satisfies $2\le a+c\le N$, and the number of factors
is at most $\log_2 N$. Every nonempty set of monomials therefore has a
least element, obtained by taking its least product score and then
the least of the finitely many ties. Restriction changes none of
these comparisons.
\end{proof}

\subsubsection{Minimal binary Goldbach moves}
For every represented sum $n$, let
\[
 \mathcal M_{n,2}=\{x_ax_{a'}:a,a'\in\A,\ a\le a',\ a+a'=n\}.
\]
Choose $b\le b'$ in $\A$ with $b+b'=n$ and $b'-b$ as small as
possible. We call $u_n=x_bx_{b'}$ the \emph{anchor} of this fiber.
If $n/2\in\A$, the anchor is the square $x_{n/2}^2$.

\begin{proposition}[Anchors and minimal generators]\label{moves:anchors}
The anchor $u_n$ is the unique largest monomial in $\mathcal M_{n,2}$
for the product order, independently of the shift $c$. The binomials
\[
 u_n-x_ax_{a'},\qquad
 x_ax_{a'}\in\mathcal M_{n,2}\setminus\{u_n\},
\]
as $n$ varies, form a minimal bihomogeneous generating set of $I$
over every field. If a quadratic fiber is nonempty in $R^{<m}$, its
anchor already belongs to $R^{<m}$. Consequently these chosen
generating sets are nested under variable truncation.
\end{proposition}
\begin{proof}
For $a+a'=n$ and $a\le a'$,
\[
 (a+c)(a'+c)=\frac{(n+2c)^2-(a'-a)^2}{4}.
\]
The sum and the gap determine the pair. Thus $u_n$ is the unique
product maximum for every $c$. For any other pair we have
$a<b\le b'<a'$: decreasing the gap at fixed sum raises the smaller
label and lowers the larger one. Hence if $a'\le m$, both $b$ and $b'$
are at most $m$. This proves the assertion about cutoffs.

For any two monomials $v,w\in\mathcal M_{n,2}$, their difference
is $(u_n-w)-(u_n-v)$. Within a fiber, the displayed binomials are
linearly independent because each nonanchor monomial appears in just one.
Different fibers have different Goldbach degrees. Since all generators
have monomial degree two and $\m I$ has monomial degree at least
three, their classes form a basis of $I/\m I$. This proves
minimality. The fixed anchors give the asserted nesting.
\end{proof}

\begin{definition}[Minimal binary Goldbach move]\label{moves:minimal}
A \emph{minimal binary Goldbach move} is a replacement
\[
 x^\alpha x_ax_{a'}\longrightarrow x^\alpha u_n,\qquad
 x_ax_{a'}\in\mathcal M_{n,2}\setminus\{u_n\},
\]
where $x^\alpha$ is the fixed monomial factor. The arrow points
from the nonanchor multiple to the anchor multiple, and hence
increases the product order. We also allow traversal in the opposite
direction when discussing undirected paths.
\end{definition}

Here ``minimal'' refers to the chosen minimal generating set of the
ideal. When we wish to specify a quadratic generator by the monomial
$v=x_ax_{a'}$ it replaces, we write $\left.e\right|_v$. Thus
$\partial(\left.e\right|_v)=u_n-v$ for $n=a+a'$.

\begin{proposition}[The vertices and edges of a strand]\label{moves:graph}
Choose these quadratic generators as the basis of $F_1$.
For every bidegree $(n,d)$,
the map $F_{1,n,d}\to F_{0,n,d}$ is the oriented incidence map
of the graph whose vertices are the monomials $x^\alpha$ of Goldbach degree
$n$ and monomial degree $d$
and whose edges are the minimal binary Goldbach moves in that fiber.
Its connected components are those of $\Gamma_{n,d}$.
\end{proposition}
\begin{proof}
Since $F_{0,n,d}=R_{n,d}$, its basis consists of the stated monomials.
The basis of $F_{1,n,d}$ consists of the multiples
$x^\alpha(\left.e\right|_v)$ of bidegree $(n,d)$. For $v=x_ax_{a'}$,
\[
 \partial\bigl(x^\alpha(\left.e\right|_v)\bigr)
     =x^\alpha u_{a+a'}-x^\alpha x_ax_{a'}.
\]
This is exactly the incidence boundary of the indicated edge.
Any elementary binary move between two multiples of pairs of sum
$a+a'$ with the same fixed factor $x^\alpha$ passes through $x^\alpha u_{a+a'}$
using at most two star edges. Thus the chosen edges connect the
vertices in each component of $\Gamma_{n,d}$.
\end{proof}

Along an arrow $x^\alpha x_ax_{a'}\to x^\alpha u_n$, the pair $(a',a)$ majorizes
$(b',b)$. Appending the labels of $x^\alpha$ to both pairs
preserves majorization: the sum of the $j$ largest labels is the
maximum over all ways to choose $j$ labels from the pair and from $x^\alpha$.
Thus the source partition majorizes the target along every directed path.
Every arrow strictly increases the product order, so the move graph
has no directed cycles. Its reachability relation is a partial order,
but the graph need not be its Hasse diagram: a chosen edge may also
be traversed by a longer directed path. Nor do these quadratic
generators have to be a Gr\"obner basis. The product order selects
and orients the edges; choosing higher syzygies and realizing their
boundaries by cells is a further problem.
Majorization alone does not create an edge: $\Gamma_{n,d}$ allows all
binary replacements, while the chosen graph retains only anchor moves.

\begin{proposition}[Maximal representatives]\label{moves:maximal-representatives}
Every connected component $C$ of $\Gamma_{n,d}$ has a unique greatest
monomial $m_C$ in the product order, called its \emph{component maximum}.
The classes of these monomials form a basis of $G_{n,d}$.
From any vertex, successive product-increasing minimal binary moves
terminate at a \emph{one-move local maximum}: a monomial from which
every nontrivial binary move strictly decreases $P_c$.
\end{proposition}
\begin{proof}
The fiber is finite, so each component has a greatest element in the
total product order, and a strictly increasing sequence of moves
terminates. At a terminal vertex every pair of factors is its sum's
anchor; otherwise replacing that pair by its anchor would increase
the product. Every nontrivial replacement of an anchor decreases
$P_c$. The basis assertion is Proposition~\ref{setup:components}.
\end{proof}

\begin{example}[Local maxima need not be unique]\label{moves:local-maxima}
On the odd primes, the minimal moves in bidegree $(21,3)$ include
\[
 x_5^2x_{11}\;\longleftrightarrow\;
 x_3x_7x_{11}\;\longleftrightarrow\;x_7^3,
 \qquad 275\;>\;231\;<\;343.
\]
The displayed numbers are the products of the labels. Both endpoints
are one-move local maxima: their quadratic factors are the anchors
$u_{10}=x_5^2$, $u_{16}=x_5x_{11}$, and $u_{14}=x_7^2$.
The component maximum is $x_7^3$, since the product of three positive
labels with sum $21$ is at most $7^3$. Thus reaching the component
maximum from a local maximum can require a product-decreasing move.
\end{example}

\subsubsection{Illustrated minimal-resolution strands}
For the examples in this subsection, let $\A=\Podd$ be the positive
odd primes and let $G$ be the full binary quotient. We now describe
the $\kk$-complex $F_{n,d}$ for several choices of
$(n,d)$. The vertices are labeled by the monomials $x^\alpha\in R_{n,d}$,
the edges correspond to the chosen quadratic Goldbach relations,
and the two-dimensional cells record relations among these moves.
We use the maximum-product anchor in each quadratic fiber. A
quadratic basis element has boundary equal to the anchor
minus the nonanchor monomial. Arrows in the following pictures indicate this
cellular boundary orientation; they point toward increasing product.

We label a generator of Goldbach degree $n$ and monomial degree $d$
by $e_{n,d}$, using $e_{n,d}',e_{n,d}'',\ldots$ when several occur.
The two subscripts record its bidegree. In the quadratic fiber we
label the generators $\left.e\right|_v$ in order of increasing
$P_c(u_n)-P_c(v)$, beginning with $e_{n,2}$. These gaps are distinct.
Thus each additional prime indicates a smaller second monomial $v$
in the product order. In higher degrees the primes distinguish the
chosen generators of the same bidegree.
The quadratic and cubic generators in the following pictures have
homological degrees one and two, respectively.

\begin{example}[The strand $F_{18,2}$]\label{ex:original18}
The $\kk$-complex $F_{18,2}$ is supported on the labeled interval below:
\begin{center}
\begin{tikzpicture}[scale=0.8]
\node[circle,fill=black,inner sep=1.5pt,label=below:$x_5x_{13}$]
  (e1) at (-2,0) {};
\node[circle,fill=black,inner sep=1.5pt,label=below:$x_7x_{11}$]
  (e3) at (2,0) {};
\draw[thick,mid arrow={Stealth[reversed,length=1.45mm,width=1mm]}] (e3)--(e1);
\node at (0,0.3) {$e_{18,2}$};
\end{tikzpicture}
\end{center}
In particular, in homological degree zero we have
\[
 F_{0,18,2}=R_{18,2}
 =\kk x_5x_{13}\oplus\kk x_7x_{11},
\]
where the two monomials correspond to the vertices. In homological
degree one, $F_{1,18,2}=\kk e_{18,2}$, and the differential is the
usual boundary of an oriented edge:
\[
 \partial(e_{18,2})=x_7x_{11}-x_5x_{13}.
\]
\end{example}

\Needspace{6cm}
\begin{example}[The strand $F_{42,2}$]\label{ex:original42}
The maximum-product anchor is $x_{19}x_{23}$, and $F_{42,2}$ is
supported on the following star.
\begin{center}
\begin{tikzpicture}[scale=0.8]
\node[circle,fill=black,inner sep=1.5pt,label=below:$x_{19}x_{23}$]
  (e1) at (0,0) {};
\node[circle,fill=black,inner sep=1.5pt,label=above:$x_{11}x_{31}$]
  (e2) at (-4,4) {};
\node[circle,fill=black,inner sep=1.5pt,label=above:$x_{13}x_{29}$]
  (e3) at (0,4) {};
\node[circle,fill=black,inner sep=1.5pt,label=above:$x_5x_{37}$]
  (e4) at (4,4) {};
\draw[thick,mid arrow={Stealth[reversed,length=1.45mm,width=1mm]}] (e1)--(e2);
\draw[thick,mid arrow={Stealth[reversed,length=1.45mm,width=1mm]}] (e1)--(e3);
\draw[thick,mid arrow={Stealth[reversed,length=1.45mm,width=1mm]}] (e1)--(e4);
\node at (-2.8,2) {$e_{42,2}'$};
\node at (-0.5,2.7) {$e_{42,2}$};
\node at (2.8,2) {$e_{42,2}''$};
\end{tikzpicture}
\end{center}
The boundaries are
\[
\begin{aligned}
 \partial(e_{42,2})&=x_{19}x_{23}-x_{13}x_{29},\\
 \partial(e_{42,2}')&=x_{19}x_{23}-x_{11}x_{31},\\
 \partial(e_{42,2}'')&=x_{19}x_{23}-x_5x_{37}.
\end{aligned}
\]
In $\Gamma_{42,2}$ a single edge joins $x_{13}x_{29}$ to $x_{11}x_{31}$
because $13+29=11+31$. In the pictured star each edge is one basis
element. The path between these leaves goes first along $e_{42,2}$ to
the anchor, then backwards along $e_{42,2}'$. Its boundary is
\[
 \partial(e_{42,2}-e_{42,2}')=x_{11}x_{31}-x_{13}x_{29}.
\]
\end{example}

\Needspace{8cm}
\begin{example}[The strand $F_{23,3}$]\label{ex:original23}
The strand $F_{23,3}$ is supported on the labeled simplicial complex:
\nopagebreak
\begin{center}
\begin{tikzpicture}[scale=0.65,every node/.style={font=\small}]
\draw[fill=blue!12] (3,0)--(6,4)--(6,-4)--cycle;
\node[circle,fill=black,inner sep=1.5pt,label=below:$x_3^2x_{17}$]
  (e2) at (-3,0) {};
\node[circle,fill=black,inner sep=1.5pt,label=below left:$x_3x_7x_{13}$]
  (e4) at (3,0) {};
\node[circle,fill=black,inner sep=1.5pt,label=above:$x_5^2x_{13}$]
  (e6) at (6,4) {};
\node[circle,fill=black,inner sep=1.5pt,label=below:$x_5x_7x_{11}$]
  (e8) at (6,-4) {};
\draw[thick,mid arrow={Stealth[reversed,length=1.45mm,width=1mm]}] (e4)--(e2);
\draw[thick,mid arrow={Stealth[reversed,length=1.45mm,width=1mm]}] (e6)--(e4);
\draw[thick,mid arrow={Stealth[reversed,length=1.45mm,width=1mm]}] (e8)--(e6);
\draw[thick,mid arrow={Stealth[reversed,length=1.45mm,width=1mm]}] (e8)--(e4);
\node at (0,0.4) {$x_3e_{20,2}$};
\node at (3.7,2.3) {$x_{13}e_{10,2}$};
\node[right] at (6.2,0) {$x_5e_{18,2}$};
\node at (3.8,-2.2) {$x_7e_{16,2}$};
\node at (4.7,0) {$e_{23,3}$};
\end{tikzpicture}
\end{center}
The quadratic boundaries are
\[
\begin{aligned}
 \partial e_{10,2}&=x_5^2-x_3x_7,&
 \partial e_{16,2}&=x_5x_{11}-x_3x_{13},\\
 \partial e_{18,2}&=x_7x_{11}-x_5x_{13},&
 \partial e_{20,2}&=x_7x_{13}-x_3x_{17}.
\end{aligned}
\]
For instance, in $G_{23,3}$ we have
\[
 x_3x_7x_{13}=x_5^2x_{13}=x_5x_7x_{11}=x_3x_7x_{13}.
\]
This computation is visualized by tracing the boundary of the
two-dimensional cell labeled $e_{23,3}$. The differential is a cellular
boundary map with the additional requirement of $R$-linearity; for
example,
\[
 \partial(x_{13}e_{10,2})=x_{13}\partial(e_{10,2})
 =x_{13}(x_5^2-x_3x_7).
\]
In homological degree two, choose
\[
 \partial(e_{23,3})=x_5e_{18,2}-x_7e_{16,2}+x_{13}e_{10,2}.
\]
Applying $\partial$ again makes the three monomial differences
cancel. The graph has four vertices and four edges and is connected,
so this cycle spans its one-dimensional cycle space. Consequently
it supplies the entire second homological term in this strand over
every field.

The picture also illustrates variable truncation. All its labels
belong to $R^{<18}$. Passing to $R^{<14}$ removes the vertex
$x_3^2x_{17}$ and its incident edge $x_3e_{20,2}$, leaving the filled
triangle. Thus the cubic generator $e_{23,3}$ is already visible before
the full strand stabilizes. With these compatible choices of bases,
$F_{23,3}^{<m}=F_{23,3}$ for $m\geq17$. The relation cutoff
$F_{23,3}^{\le14}$ has the same three vertices but only the edge
$x_{13}e_{10,2}$. At $m=16$ it also has $x_7e_{16,2}$, giving a
two-edge path. The vertex $x_3^2x_{17}$ enters at $m=17$. At $m=18$ the
triangle fills in, while that vertex remains isolated until
$x_3e_{20,2}$ joins it to the triangle at $m=20$.
\end{example}

\Needspace{11cm}
\begin{example}[The strand $F_{31,3}$]\label{ex:original31}
There are six monomials in $R_{31,3}$. To specify the quadratic
generators used in this picture, set
\[
\begin{aligned}
 \partial e_{14,2}&=x_7^2-x_3x_{11},&
 \partial e_{24,2}&=x_{11}x_{13}-x_7x_{17},\\
 \partial e_{24,2}'&=x_{11}x_{13}-x_5x_{19},&
 \partial e_{26,2}&=x_{13}^2-x_7x_{19},\\
 \partial e_{26,2}'&=x_{13}^2-x_3x_{23},&
 \partial e_{28,2}&=x_{11}x_{17}-x_5x_{23}.
\end{aligned}
\]
Together with $e_{18,2}$ and $e_{20,2}$ above, these give the eight
edges below. Each edge carries its full generator label $x_pe_{s,2}^{(j)}$,
where $x_p$ is the fixed variable of the binary move. The middle quadrilateral and the two triangles
are the three two-dimensional cells.
\begin{center}
\begin{tikzpicture}[scale=.96,line join=round,line cap=round,
    dot/.style={circle,fill=black,inner sep=1.5pt},
    vertex/.style={font=\normalsize,inner sep=2pt},
    edge label/.style={font=\small,inner sep=3pt},
    face label/.style={font=\normalsize,inner sep=2pt},
    ed/.style={semithick,mid arrow={Stealth[length=1.45mm,width=1mm]}}]
  \coordinate (F) at (0,7.0);
  \coordinate (D) at (-5.05,3.15);
  \coordinate (E) at (-1.8,3.65);
  \coordinate (B) at ( 1.8,3.65);
  \coordinate (C) at ( 5.05,4.35);
  \coordinate (A) at (0,.35);
  \fill[blue!12] (F)--(D)--(E)--(A)--(B)--(C)--cycle;
  \draw[ed] (D)--(F);
  \draw[ed] (D)--(E);
  \draw[ed] (E)--(F);
  \draw[ed] (A)--(E);
  \draw[ed] (A)--(B);
  \draw[ed] (B)--(F);
  \draw[ed] (B)--(C);
  \draw[ed] (C)--(F);
  \path (D)--(F) node[midway,sloped,above,edge label] {$x_7e'_{24,2}$};
  \path (D)--(E) node[midway,sloped,below,edge label] {$x_5e_{26,2}$};
  \path (E)--(F) node[pos=.60,sloped,below,edge label] {$x_{13}e_{18,2}$};
  \path (A)--(E) node[midway,sloped,below,edge label] {$x_5e'_{26,2}$};
  \path (A)--(B) node[midway,sloped,below,edge label] {$x_3e_{28,2}$};
  \path (F)--(B) node[pos=.40,sloped,below,edge label] {$x_{11}e_{20,2}$};
  \path (B)--(C) node[midway,sloped,below,edge label] {$x_{17}e_{14,2}$};
  \path (F)--(C) node[midway,sloped,above,edge label] {$x_7e_{24,2}$};
  \foreach \name in {A,B,C,D,E,F} \node[dot] at (\name) {};
  \node[vertex,above=3pt] at (F) {$x_7x_{11}x_{13}$};
  \node[vertex,left=3pt] at (D) {$x_5x_7x_{19}$};
  \node[vertex,right=4pt] at (E) {$x_5x_{13}^2$};
  \node[vertex,left=4pt] at (B) {$x_3x_{11}x_{17}$};
  \node[vertex,right=3pt] at (C) {$x_7^2x_{17}$};
  \node[vertex,below=3pt] at (A) {$x_3x_5x_{23}$};
  \node[face label] at (-2.35,4.62) {$e'_{31,3}$};
  \node[face label] at (0,2.55) {$e_{31,3}$};
  \node[face label] at (2.35,4.99) {$e''_{31,3}$};
\end{tikzpicture}
\end{center}
The height in this drawing records the direction of the chosen
moves. Its top vertex is $x_7x_{11}x_{13}$, and its two sources are
$x_5x_7x_{19}$ and $x_3x_5x_{23}$. These are minimal elements
of the reachability order for the selected star edges. The numerical
prime-product order is total and has only the latter as its minimum.
The reachability order is not a lattice: the two sources have no common
lower bound. Nor is the displayed graph a Hasse diagram; for example,
its direct edge from $x_5x_7x_{19}$ to $x_7x_{11}x_{13}$ is
accompanied by a two-edge path through $x_5x_{13}^2$.

On the edge from $x_5x_7x_{19}$ to $x_7x_{11}x_{13}$, the fixed
factor $x_7$ moves from the second position to the first in the
sorted list of primes. The same change occurs for $x_5$ on the edge
from $x_3x_5x_{23}$ to $x_5x_{13}^2$. By contrast, on the edge
from $x_5x_{13}^2$ to $x_7x_{11}x_{13}$, the fixed $x_{13}$ can
occupy the third position at both ends.

The signed cellular boundaries are
\[
\begin{aligned}
 \partial e_{31,3}
   &=x_3e_{28,2}+x_{11}e_{20,2}-x_{13}e_{18,2}-x_5e_{26,2}',\\
 \partial e_{31,3}'
   &=x_5e_{26,2}-x_7e_{24,2}'+x_{13}e_{18,2},\\
 \partial e_{31,3}''
   &=x_{17}e_{14,2}+x_7e_{24,2}-x_{11}e_{20,2}.
\end{aligned}
\]
Each expression has boundary zero. They are independent because the
edges $x_3e_{28,2}$, $x_5e_{26,2}$, and $x_{17}e_{14,2}$ occur in just
one of the three respective expressions. The graph is connected,
with six vertices and eight edges, so its cycle space has dimension
three. Hence these are a basis of the cubic syzygies in this
bidegree, and the displayed complex supports $F_{31,3}$ over every
field.

The two computations
\[
\begin{aligned}
 x_3x_5x_{23}&=x_3x_{11}x_{17}=x_7^2x_{17},\\
 x_3x_5x_{23}&=x_5x_{13}^2=x_7x_{11}x_{13}=x_7^2x_{17}
\end{aligned}
\]
trace different paths between the same vertices. Their difference
is the boundary of the middle quadrilateral together with the right-hand triangle,
with compatible orientations. This is how a relation among sequences
of Goldbach moves appears in the next homological degree.
\end{example}

The last picture also realizes the compatible cutoff construction.
Its right-hand triangle, with generator $e_{31,3}''$, is present in
$F_{31,3}^{<17}$. At cutoff $19$ the left-hand triangle, with generator
$e_{31,3}'$, is added. At cutoff $23$ the middle quadrilateral, with generator
$e_{31,3}$, is added. Each stage retains every earlier cell and
attaching map. Thus these particular illustrated cutoffs form actual
polyhedral subcomplexes. This observation concerns the stated examples;
it does not assert such models for all strands or all cutoffs.

The next cubic example, in weight $33$, also connects these
polygonal strands to explicit Koszul classes. We give its picture
and calculation in Example~\ref{ex:weight33}, after defining the
Koszul complex.

\subsubsection{The quartic three-cell}
\begin{example}[The three-cell in $F_{34,4}$]
\label{ex:pyramid34}
The quadratic prime-product anchors and the following compatible
cubic choices give a unique three-cell in this strand. Its closure,
labeled $e_{34,4}$, is a square pyramid.
The diagram below is its top-down boundary view: the central
vertex is the apex, the four interior triangles are its lateral faces,
and the outer face represents its square base. Vertices carry monomial
labels, edges carry their quadratic-generator labels, and faces carry
their cubic-generator labels. The small edge arrows increase prime product.

\begin{center}
\begin{tikzpicture}[
  scale=1.05,line join=round,
  v/.style={circle,fill=black,inner sep=1.6pt},
  vertex label/.style={font=\normalsize,inner sep=3pt},
  edge label/.style={font=\small,inner sep=2pt},
  face label/.style={font=\normalsize,inner sep=2pt},
  ed/.style={line width=.7pt,mid arrow={Stealth[reversed,length=1.45mm,width=1mm]}}]
\coordinate (A) at (-4,-4);
\coordinate (B) at (4,-4);
\coordinate (C) at (4,4);
\coordinate (D) at (-4,4);
\coordinate (U) at (0,0);
\fill[blue!12] (A)--(B)--(C)--(D)--cycle;
\draw[ed] (B)--(A);
\draw[ed] (C)--(B);
\draw[ed] (C)--(D);
\draw[ed] (A)--(D);
\node[edge label,below=7pt] at (0,-4) {$x_5x_{13}e_{16,2}$};
\node[edge label,right=7pt] at (4,0) {$x_5x_{11}e_{18,2}$};
\node[edge label,above=7pt] at (0,4) {$x_5x_7e_{22,2}'$};
\node[edge label,left=7pt] at (-4,0) {$x_3x_5e_{26,2}$};
\draw[ed] (B)--(U)
  node[pos=.53,sloped,below=4pt,edge label] {$x_{11}x_{13}e_{10,2}$};
\draw[ed] (U)--(A)
  node[pos=.53,sloped,above=4pt,edge label] {$x_3x_{13}e_{18,2}$};
\draw[ed] (C)--(U)
  node[pos=.53,sloped,above=4pt,edge label] {$x_7x_{11}e_{16,2}$};
\draw[ed] (U)--(D)
  node[pos=.53,sloped,above=4pt,edge label] {$x_3x_7e'_{24,2}$};
\foreach \P in {A,B,C,D,U}{\node[v] at (\P) {};}
\node[vertex label,below=5pt] at (A) {$x_3x_5x_{13}^{2}$};
\node[vertex label,below=5pt] at (B) {$x_5^{2}x_{11}x_{13}$};
\node[vertex label,above=5pt] at (C) {$x_5x_7x_{11}^{2}$};
\node[vertex label,above=5pt] at (D) {$x_3x_5x_7x_{19}$};
\node[vertex label] (apex label) at (0,1.25) {$x_3x_7x_{11}x_{13}$};
\draw[-{Stealth[length=1.2mm,width=.8mm]},line width=.4pt]
  (apex label.south)--(0,.12);
\node[face label] at (0,-2.85) {$x_{13}e_{21,3}$};
\node[face label] at (2.85,0) {$x_{11}e_{23,3}$};
\node[face label] at (0,2.85) {$x_7e_{27,3}$};
\node[face label] at (-2.85,0) {$x_3e'_{31,3}$};
\node[face label] at (0,5.05) {Outer face: $x_5e'_{29,3}$};
\end{tikzpicture}

\end{center}
To specify the remaining boundaries, set
$\partial e_{22,2}'=x_{11}^2-x_3x_{19}$ and choose
\[
\begin{aligned}
 \partial e_{21,3}
   &=x_{11}e_{10,2}-x_5e_{16,2}+x_3e_{18,2},\\
 \partial e_{27,3}
   &=x_{11}e_{16,2}-x_5e_{22,2}'+x_3e'_{24,2},\\
 \partial e'_{29,3}
   &=x_{13}e_{16,2}+x_{11}e_{18,2}
      -x_7e_{22,2}'+x_3e_{26,2}.
\end{aligned}
\]
The boundaries of $e_{23,3}$ and $e'_{31,3}$ were given above. With compatible
face orientations, the new generator has boundary
\[
 \partial e_{34,4}
 =x_{13}e_{21,3}-x_{11}e_{23,3}-x_7e_{27,3}
   +x_5e'_{29,3}-x_3e'_{31,3}.
\]
The left triangular face is $x_3e'_{31,3}$, obtained by multiplying
the left triangle of Example~\ref{ex:original31} by $x_3$.
The square base is labeled $x_5e'_{29,3}$.

Every edge occurs twice with opposite incidence signs. The resulting
two-sphere uses five vertices, eight edges, and five faces; filling it
gives the three-cell. The full strand has dimensions $(11,20,11,1)$,
so this picture shows the support of its unique three-cell rather
than the entire strand.

For completeness, the finite integral checks for these examples are
summarized below. The dimensions refer to the full resolution strand,
including monomial multiples of lower generators.
\[
\begin{array}{c|c|c}
 (n,d)&(\dim F_{0,n,d},\dim F_{1,n,d},\ldots)
       &(\operatorname{rank}\partial_1,\operatorname{rank}\partial_2,\ldots)\\\hline
 (18,2)&(2,1)&(1)\\
 (42,2)&(4,3)&(3)\\
 (23,3)&(4,4,1)&(3,1)\\
 (31,3)&(6,8,3)&(5,3)\\
 (33,3)&(8,11,4)&(7,4)\\
 (34,4)&(11,20,11,1)&(10,10,1)
\end{array}
\]
Here is an exact way to verify the last line, rather than infer
uniqueness from the pictured sphere. Enumerate the eleven quartic
monomials, all twenty monomial multiples of the quadratic star
relations, and all eleven variable multiples of the chosen cubic
cycles. Their integer boundary matrices have consecutive products
zero. Integer elimination uses only pivots $1$ or $-1$ and leaves
zero residual matrices, with ranks $10$ and $10$ for the first two
maps. The displayed three-cell boundary is a primitive generator
of the remaining one-dimensional kernel. Thus the existing cubic
multiples fill the entire first cycle group, and exactly one new
third generator is required. These checks hold over every field.

The data are in
\path{supplement/computations/illustrated-strands-data.json}, with the
standard-library checker
\path{supplement/computations/illustrated-strands-check.py}. It reconstructs the
matrices from polynomial boundaries, checks all cubic generator
weights through $31$, and verifies the five face boundaries of this
pyramid. For the example in weight $33$, \path{supplement/computations/koszul33-check.py}
checks the four polygon cycles, the swaps, and their Koszul classes over the integers. All polynomial
coefficients in positive differentials lie in $\m$, so these are
minimal generators. The finite choices can be extended to a minimal
resolution by the preceding construction. This verifies the examples
without assuming a global cellular resolution.
\end{example}

The companion paper \cite{CompanionTopology} compares these binary
strands with strands for the ternary relation ideal.

The binary examples above show how minimal binary moves form the first stage of
an inductive construction: polygon boundaries supply cubic syzygies,
and compatible polygon faces can bound a three-cell. In general,
however, the next step chooses minimal homogeneous generators of a
syzygy module. A product order on monomials does not yet choose those
higher generators or ensure that their boundaries are spheres.
Extending the order and the moves to compatible higher cells is an
open problem here; the displayed models establish only the stated
strands and cutoff inclusions.

\section{Koszul homology and translation classes}
\label{sec:koszul}

The finite-resolution strands give geometric representatives of
syzygies. The Koszul complex gives a second description that is
particularly effective for constructing classes and counting them.
We first define it for a Goldbach ring of any fixed arity $k\ge2$.
We then return to binary rings, identify the polygon classes in
weight $33$, and construct the common-shift classes used in the
Betti estimates.

\subsection{The Koszul complex and its gradings}
Let $\kk$ be our coefficient field of characteristic $\ell$, where
$\ell=0$ or $\ell$ is prime, and put
$R_{\mathbb Z}=\mathbb Z[x_a\mid a\in\A]$ and
$R=R_{\mathbb Z}\otimes_{\mathbb Z}\kk=\kk[x_a\mid a\in\A]$.
We reserve $\mathbb K$ for the integral Koszul free resolution of
$\mathbb Z$ over $R_{\mathbb Z}$, with augmentation sending every
$x_a$ to zero. For infinitely many variables, $\mathbb K$ is the
filtered direct limit of the finite Koszul complexes, extended to
$R_{\mathbb Z}$; its homology is $\mathbb Z$ in homological degree
zero and zero otherwise.
Its scalar extension
$\mathbb K_{\kk}:=\mathbb K\otimes_{\mathbb Z}\kk$
is the Koszul free resolution of $\kk$ over $R$.
For an $R$-module $M$, put
\[
 K(M)=\mathbb K_{\kk}\otimes_R M,\qquad
 K=K(G)=\mathbb K_{\kk}\otimes_R G.
\]
Thus
\[
 \mathrm{H}_i(K)_{n,d}=\Tor_i^R(G,\kk)_{n,d}.
\]
As for $F_{n,d}$, we write $K_{n,d}$ for the whole complex in
bidegree $(n,d)$ and $K_{i,n,d}$ for its term in homological degree $i$.
The complex $K$ has free $G$-modules as its terms, but generally has
nonzero positive homology, so it is not a resolution of $\kk$ over $G$.
Our Betti numbers are
$\beta_{i,n,d}=\dim_\kk\mathrm{H}_i(K)_{n,d}$.

For a finite subset $\sigma=\{a_1<\cdots<a_i\}\subseteq\A$, write
$\varepsilon_\sigma=\varepsilon_{a_1}\cdots\varepsilon_{a_i}$, with
$\varepsilon_\varnothing=1$. The $R_{\mathbb Z}$-linear differential
on $\mathbb K$ is
\begin{equation}\label{setup:koszul-differential}
 \mathrm{d}(\varepsilon_\sigma)=
 \sum_{j=1}^i(-1)^{j-1}x_{a_j}\varepsilon_{\sigma\setminus\{a_j\}}.
\end{equation}
The same formula gives the differential after passing to coefficient
classes in $G_{\mathbb Z}$ or $G$.

\phantomsection\label{setup:koszul-grading}
Give the generators bidegrees and homological degrees
\[
 \deg x_a=\deg\varepsilon_a=(a,1),\qquad
 |x_a|=0,\qquad |\varepsilon_a|=1,
\]
and extend both degrees additively to products. Thus, for $|\sigma|=i$,
\[
 \deg(x^\alpha\varepsilon_\sigma)
   =\deg(x^\alpha)+\sum_{a\in\sigma}(a,1),
 \qquad |x^\alpha\varepsilon_\sigma|=i.
\]
In particular, the integral strand has the free abelian terms
\[
 \mathbb K_{i,n,d}
 =\bigoplus_{\substack{\sigma\subseteq\A,\ |\sigma|=i\\
          \sum_a\alpha_a=d-i\\
          \sum_a a\alpha_a+\sum_{a\in\sigma}a=n}}
       \mathbb Z\,x^\alpha\varepsilon_\sigma,
\]
where $\alpha$ ranges over nonnegative exponent vectors of finite
support. The coefficient has bidegree
$(n-\sum_{a\in\sigma}a,d-i)$; $d$ is the monomial degree of the whole
term. The differential lowers the homological degree by one and
preserves $(n,d)$.
For example, $x_3^2\varepsilon_5\varepsilon_7$ has bidegree $(18,4)$
and homological degree $2$, while its coefficient has bidegree $(6,2)$.

To compare different fields, take the same binomial relations over
$\mathbb Z$. Let $I_{\mathbb Z}\subseteq R_{\mathbb Z}$ be the ideal
generated by these relations, and put
\[
 G_{\mathbb Z}=R_{\mathbb Z}/I_{\mathbb Z},\qquad
 K_{\mathbb Z}=\mathbb K\otimes_{R_{\mathbb Z}}G_{\mathbb Z}.
\]
Extension of scalars gives
\[
 G=G_{\mathbb Z}\otimes_{\mathbb Z}\kk,\qquad
 K\cong K_{\mathbb Z}\otimes_{\mathbb Z}\kk.
\]
Since $I_{\mathbb Z}$ is bihomogeneous, tensoring the whole complex
with $G_{\mathbb Z}$ preserves the bidegree decomposition:
\[
 K_{\mathbb Z,i,n,d}
 =\bigoplus_{\substack{\sigma\subseteq\A\\|\sigma|=i}}
   (G_{\mathbb Z})_{\,n-\sum_{a\in\sigma}a,\,d-i}
   \,\varepsilon_\sigma.
\]
Components of negative monomial degree are zero. The generators are
the classes $\overline{x^\alpha}\varepsilon_\sigma$ of the same
total bidegree $(n,d)$; the binomial relations identify coefficient
monomials only within a single bidegree.
The component-basis proof of Proposition~\ref{setup:components} works
over $\mathbb Z$, so these terms remain free abelian. Since $\A$ is
bounded below, each strand $K_{\mathbb Z,n,d}$ is a finite complex of finitely generated
free abelian groups. Consequently these complexes can be compared
before choosing a field. Their integral homology and cell models are
used in \cite{CompanionTopology} to determine characteristic-dependent
Betti numbers.

\subsection{Finite strands and stabilization}
The same finiteness gives a direct stabilization statement for
Koszul complexes. Unlike the compatible choices of minimal
resolutions in Section~\ref{setup:resolution-construction}, the maps
here are the canonical maps on coefficient monomials and exterior
generators. Before taking cutoffs, apply the affine translation of
Proposition~\ref{setup:affine-invariance} if necessary so that
$\A\subseteq\mathbb Z_{>0}$, and put $\A_m=\A\cap[1,m]$.
For a fixed arity $k$, define $G^{<m}$ by restricting all
$k$ labels on each side of every defining relation to $\A_m$; define
$G^{\le m}$ by also restricting the common Goldbach degree to at most
$m$. This extends the binary cutoff notation to the chosen arity.

\begin{proposition}[Stabilization in a fixed bidegree]
\label{setup:stabilization}
Let $\A\subseteq\mathbb Z_{>0}$, fix an arity $k\ge2$, put
$b=\min\A$, and let $\A_m=\A\cap[1,m]$.
Denote by $K_{\mathbb Z}^{<m}$ the Koszul
complex on the variables indexed by $\A_m$, with coefficients in
$G_{\mathbb Z}^{<m}$. Sending each
$\overline{x^\alpha}\varepsilon_\sigma$ to the same expression over
$G_{\mathbb Z}$ induces an isomorphism of $\mathbb Z$-complexes
\[
 K_{\mathbb Z,n,d}^{<m}
 \longrightarrow K_{\mathbb Z,n,d}
 \qquad\text{if }d\ge1\text{ and }m\ge n-(d-1)b.
\]
The same assertion holds with coefficient ring
$G_{\mathbb Z}^{\le m}$ whenever $m\ge n$.
\end{proposition}
\begin{proof}
Every monomial $x^\alpha\varepsilon_\sigma$ of bidegree $(n,d)$
contains $d$ labels counting multiplicities. Each label $a$ therefore
satisfies
\[
 a+(d-1)b\le n.
\]
Thus every label is at most $m$, so all generators of this strand
occur in the stated cutoff.
The relations in each exterior summand are generated over $\mathbb Z$
by
\[
 x^\gamma(x^\mu-x^\nu)\varepsilon_\sigma,
\]
where $x^\mu-x^\nu$ is a defining binomial and both terms have total
bidegree $(n,d)$. The same bound puts every label in each such relation
inside the cutoff. The two strands therefore have identical generators
and relations, and the differential has the same formula.
By positivity, the common weight of each relevant defining
binomial is at most $n$, so $m\ge n$ also includes all required relations
in the relation-weight truncation.
\end{proof}

Before this bound is reached, new relations may identify old coefficient
classes, so the canonical maps need not be injective.

\subsection{From polygonal strands to Koszul classes}
We now take $G$ to be the binary Goldbach ring on the positive odd
primes. The following example makes the connection between the two
descriptions of Tor explicit. In particular, its upper quadrilateral
will produce the common-shift class used in the introduction.

\Needspace{6cm}
\begin{example}[The strand $F_{33,3}$ and its Koszul cycles]\label{ex:weight33}
The vertices are the monomials $x^\alpha\in\mathcal M_{33,3}$, one for each
of the eight unordered partitions of $33$ into three odd primes.
The eleven product-star edges form the graph below. Its three
quadrilaterals and one triangle subdivide a heptagon, so their
oriented boundaries form an integral basis of its cycle group.
They give the four minimal generators in homological degree two
and bidegree $(33,3)$ over every field. We label the triangle
$e_{33,3}$ and the left, right, and upper quadrilaterals
$e'_{33,3}$, $e''_{33,3}$, and $e'''_{33,3}$, respectively.

\begin{center}
\begin{tikzpicture}[line join=round,line cap=round,
    dot/.style={circle,fill=black,inner sep=1.5pt},
    vertex/.style={font=\small,inner sep=2pt},
    edge label/.style={font=\small,inner sep=1pt},
    face label/.style={font=\normalsize},
    ed/.style={line width=.7pt,mid arrow={Stealth[length=1.45mm,width=1mm]}}]
  \coordinate (H) at (0,3.65);
  \coordinate (B) at (-3.75,2.05);  \coordinate (E) at (3.75,2.05);
  \coordinate (F) at (-4.55,-.65);  \coordinate (D) at (4.55,-.65);
  \coordinate (G) at (-2.05,-3.35); \coordinate (A) at (2.05,-3.35);
  \coordinate (C) at (0,.10);
  \fill[blue!12] (H)--(B)--(F)--(G)--(A)--(D)--(E)--cycle;
  \draw[ed] (B)--(H); \draw[ed] (E)--(H);
  \draw[ed] (B)--(F); \draw[ed] (D)--(E);
  \draw[ed] (F)--(G); \draw[ed] (A)--(D); \draw[ed] (A)--(G);
  \draw[ed] (B)--(C); \draw[ed] (C)--(E);
  \draw[ed] (C)--(G); \draw[ed] (A)--(C);
  \foreach \v in {A,B,C,D,E,F,G,H} \node[dot] at (\v) {};
  \node[vertex,above] at (H) {$x_{11}^3$};
  \node[vertex,above left] at (B) {$x_3x_{11}x_{19}$};
  \node[vertex,above right] at (E) {$x_5x_{11}x_{17}$};
  \node[vertex,left] at (F) {$x_7^2x_{19}$};
  \node[vertex,right] at (D) {$x_5^2x_{23}$};
  \node[vertex,below left] at (G) {$x_7x_{13}^2$};
  \node[vertex,below right] at (A) {$x_3x_7x_{23}$};
  \node[vertex] at (0,.69) {$x_3x_{13}x_{17}$};
  \draw[line width=.35pt] (0,.49)--(C);
  \node[face label] at (0,-2.58) {$e_{33,3}$};
  \node[face label] at (-2.78,-.72) {$e'_{33,3}$};
  \node[face label] at (2.78,-.72) {$e''_{33,3}$};
  \node[face label] at (0,2.24) {$e'''_{33,3}$};
  \node[edge label,above] at (-1.88,3.02) {$x_{11}e'_{22,2}$};
  \node[edge label,above] at (1.88,3.02) {$x_{11}e_{22,2}$};
  \node[edge label,left] at (-4.37,.69) {$x_{19}e_{14,2}$};
  \node[edge label,right] at (4.37,.69) {$x_5e_{28,2}$};
  \node[edge label,left] at (-3.51,-2.01) {$x_7e_{26,2}$};
  \node[edge label,right] at (3.51,-2.01) {$x_{23}e_{10,2}$};
  \node[edge label,below] at (0,-3.47) {$x_7e'_{26,2}$};
  \path (B)--(C) node[midway,sloped,below=3pt,edge label] {$x_3e_{30,2}$};
  \path (C)--(E) node[midway,sloped,below=3pt,edge label] {$x_{17}e_{16,2}$};
  \path (G)--(C) node[midway,sloped,above=3pt,edge label] {$x_{13}e_{20,2}$};
  \path (C)--(A) node[midway,sloped,above=3pt,edge label] {$x_3e'_{30,2}$};
\end{tikzpicture}
\end{center}
The two edges incident to $x_{11}^3$ form the translated star $x_{11}F_{22,2}$.
At $x_{11}^3$ every binary move fixes $x_{11}$, and the two other
pairs of weight $22$ are $(3,19)$ and $(5,17)$. Thus every polygon
passing through $x_{11}^3$ uses both edges.

The fixed factor on each edge also describes the corresponding class in
$\Tor_2^R(G,\kk)_{33,3}$. Use $K=\mathbb K_{\kk}\otimes_R G$, with
$\mathbb K_{\kk}$ the field Koszul resolution introduced above.
At a vertex $x^\alpha$ of an oriented
polygon, let $p$ be the fixed label of the incoming edge and $q$ that of
the outgoing edge. When $p\ne q$ the corner contributes
\[
 \frac{x^\alpha}{x_p x_q}\,\varepsilon_p\varepsilon_q.
\]
When $p=q$ the corner contributes zero since $\varepsilon_p^2=0$.
To follow the polygon in the Koszul complex, read off the fixed labels
of its successive edges. The corner above takes us from label $p$ to
label $q$: its differential is
$(x^\alpha/x_q)\varepsilon_q-(x^\alpha/x_p)\varepsilon_p$.
Along an edge $x^\alpha\to x^{\alpha'}$ fixing $p$, the two adjacent
corners therefore contribute opposite terms with exterior factor
$\varepsilon_p$. Their coefficients agree in $G$ by the binary relation
$(x^\alpha/x_p)=(x^{\alpha'}/x_p)$. Summing around the polygon gives a cycle.

To identify this cycle with the face class in Tor, write the oriented
boundary of its generator $e$ as $\partial e=\sum_j s_j x_{p_j}e_j$,
where $e_j$ is the quadratic generator on the $j$th edge and
$s_j=\pm1$. Let $w$ be the sum of the corners lifted to the Koszul
complex $\mathbb K_{\kk}$ over $R$. In the total complex
$F\otimes_R\mathbb K_{\kk}$ the chain
\[
 e\otimes1+\sum_j s_j e_j\otimes\varepsilon_{p_j}+1\otimes w
\]
is a cycle. The Koszul differential of the middle term cancels
$\partial e\otimes1$ with the tensor sign, and its remaining differential
cancels $1\otimes\mathrm{d}w$ edge by edge. The two augmentations give $e$ modulo
$\m$ and $w$ modulo $I$, respectively.

With the indicated face labels, the following cycles $z_T,z_L,z_R,z_U$
represent $e_{33,3},e'_{33,3},e''_{33,3},e'''_{33,3}$, respectively,
in Tor. We write each exterior pair in increasing order.
\[
\begin{aligned}
 z_T&=-x_{23}\varepsilon_3\varepsilon_7+x_{17}\varepsilon_3\varepsilon_{13}-x_{13}\varepsilon_7\varepsilon_{13},\\
 z_L&=-x_{11}\varepsilon_3\varepsilon_{19}+x_{17}\varepsilon_3\varepsilon_{13}-x_{13}\varepsilon_7\varepsilon_{13}+x_7\varepsilon_7\varepsilon_{19},\\
 z_R&=-x_7\varepsilon_3\varepsilon_{23}+x_{13}\varepsilon_3\varepsilon_{17}-x_{11}\varepsilon_5\varepsilon_{17}+x_5\varepsilon_5\varepsilon_{23},\\
 z_U&=-x_{19}\varepsilon_3\varepsilon_{11}+x_{13}\varepsilon_3\varepsilon_{17}-x_5\varepsilon_{11}\varepsilon_{17}.
\end{aligned}
\]
Each term $\pm x_r\varepsilon_p\varepsilon_q$ comes from the vertex
$x_rx_px_q$. In $z_T,z_L,z_R$ this gives one term for every vertex
of the corresponding face. At the shared vertex $x_3x_{13}x_{17}$,
for instance, the left face contributes $x_{17}\varepsilon_3\varepsilon_{13}$
and the right face contributes $x_{13}\varepsilon_3\varepsilon_{17}$.
The partition is the same. The outgoing edge fixes $13$ on the left face
and $17$ on the right. On the upper face the corner at $x_{11}^3$ contributes
$x_{11}\varepsilon_{11}^2=0$, since both incident edges fix $11$.
The three terms of $z_U$ come from its other three vertices.

For the purpose of following these cycles, one may suppress the
degree-two vertex $x_{11}^3$ and join $x_3x_{11}x_{19}$ directly to
$x_5x_{11}x_{17}$. The new edge fixes $11$ and uses the binary relation
$x_3x_{19}=x_5x_{17}$. This preserves the cycle group and turns the
upper quadrilateral into a triangle. It is a simplification of the
monomial graph; the full strand $F_{33,3}$ still contains the basis
monomial $x_{11}^3$.

Traverse the triangle from $x_3x_7x_{23}$ through
$x_3x_{13}x_{17}$ and $x_7x_{13}^2$. The fixed labels are $3,13,7$.
At $x_3x_{13}x_{17}$ the fixed label changes from $3$ to $13$,
giving $x_{17}\varepsilon_3\varepsilon_{13}$. The next edge keeps $13$ and replaces
$3+17$ by $7+13$. At $x_7x_{13}^2$ the corner is
$x_{13}\varepsilon_{13}\varepsilon_7$. The terms of the two corner
differentials still carrying $\varepsilon_{13}$ are therefore
\[
 +x_3x_{17}\varepsilon_{13}
 \qquad\text{and}\qquad
 -x_7x_{13}\varepsilon_{13}.
\]
They cancel in $G$ by the very move we traversed. Writing the second
corner in increasing order gives $-x_{13}\varepsilon_7\varepsilon_{13}$.
Both exterior labels are $1$ modulo $3$, so the ordered residue pair
$(1,1)$ looks the same before and after this transposition. Keeping the
prime labels on the two positions records its sign.

There is also a choice in moving a factor between the coefficient and the
exterior pair at a vertex with three distinct labels. For $a<b<c$,
\[
 \mathrm{d}(\varepsilon_a\varepsilon_b\varepsilon_c)
     =x_a\varepsilon_b\varepsilon_c-x_b\varepsilon_a\varepsilon_c+x_c\varepsilon_a\varepsilon_b.
\]
Thus one corner can be replaced by the two other sides of this
Koszul triangle.
Adding $\mathrm{d}(\varepsilon_3\varepsilon_7\varepsilon_{23})-\mathrm{d}(\varepsilon_3\varepsilon_{13}\varepsilon_{17})$ to $z_T$ gives the five-term representative
\[
 -x_{13}\varepsilon_7\varepsilon_{13}+x_3\varepsilon_7\varepsilon_{23}-x_7\varepsilon_3\varepsilon_{23}+x_{13}\varepsilon_3\varepsilon_{17}-x_3\varepsilon_{13}\varepsilon_{17}.
\]
Making the corresponding replacements at $x_3x_{13}x_{17}$ and
$x_3x_{11}x_{19}$ in $z_L$ amounts to subtracting
$\mathrm{d}(\varepsilon_3\varepsilon_{13}\varepsilon_{17})+\mathrm{d}(\varepsilon_3\varepsilon_{11}\varepsilon_{19})$. It gives
\[
 x_7\varepsilon_7\varepsilon_{19}-x_{13}\varepsilon_7\varepsilon_{13}+x_{13}\varepsilon_3\varepsilon_{17}-x_3\varepsilon_{11}\varepsilon_{19}
       -x_3\varepsilon_{13}\varepsilon_{17}-x_{19}\varepsilon_3\varepsilon_{11}.
\]
The last two terms account for the intermediate steps of these swaps.

In this bidegree $\mathrm{d}:K_{3,33,3}\to K_{2,33,3}$ has a separate column $\mathrm{d}(\varepsilon_a\varepsilon_b\varepsilon_c)$ for each
partition with three distinct parts. Consequently no boundary has a term
$x_p\varepsilon_p\varepsilon_q$ with $p\ne q$. The vertices $x_7^2x_{19}$ and $x_5^2x_{23}$
give such coordinates for $z_L$ and $z_R$. The vertex $x_7x_{13}^2$
gives one for both $z_T$ and $z_L$. The class of $z_U$ is nonzero
because on the triple $(3,11,19)$ it uses only the coordinate
$x_{19}\varepsilon_3\varepsilon_{11}$, whereas the only incoming
boundary at this triple has all three coordinates.

These four classes are independent over every field. Read the
coefficients of $x_7\varepsilon_7\varepsilon_{19}$ and $x_5\varepsilon_5\varepsilon_{23}$ to detect $z_L$ and $z_R$,
and then $x_{13}\varepsilon_7\varepsilon_{13}$ to detect $z_T$. The preceding test on the
triple $(3,11,19)$ detects the remaining $z_U$.
\end{example}

For a cubic weight $n$, the Koszul differential gives a different
graph, which we denote by $\Lambda_n$. Its vertices are the
\emph{prime labels} occurring in cubic partitions of $n$, whereas
the vertices of $\Gamma_{n,3}$ and the chosen minimal-move graph are
\emph{cubic monomials}. Join distinct labels $a,b$ in $\Lambda_n$
when $n-a-b$ is an odd prime. This edge represents the basis term
$x_{n-a-b}\varepsilon_a\varepsilon_b\in K_{2,n,3}$, with $a<b$.
If $Q_w$ denotes the common quadratic class of represented weight
$w$, the vertices represent $\eta_a=Q_{n-a}\varepsilon_a\in K_{1,n,3}$,
and
\[
 \mathrm d(x_{n-a-b}\varepsilon_a\varepsilon_b)=\eta_b-\eta_a.
\]
Thus $K_{2,n,3}\to K_{1,n,3}$ is the incidence map of $\Lambda_n$.
Each partition $n=a+b+c$ into three distinct primes also contributes
$\varepsilon_a\varepsilon_b\varepsilon_c\in K_{3,n,3}$, whose boundary
is the triangle on $a,b,c$. Fill these triangles; other triangles
in $\Lambda_n$ need not be boundaries. The graph cycle space modulo
these filled-triangle boundaries is $\Tor_2^R(G,\kk)_{n,3}$.

Let $s_n$ count the distinct-part partitions and $t_n$ the partitions
$n=2p+q$ with $p\ne q$. Write $v_n,c_n$ for the numbers of vertices
and connected components of $\Lambda_n$, including isolated vertices.
There are $3s_n+t_n$ edges, and the $s_n$ incoming triangle boundaries
have disjoint supports and unit entries. Consequently
\[
 \beta_{2,n,3}=(3s_n+t_n-v_n+c_n)-s_n
             =2s_n+t_n-v_n+c_n.
\]
For $n=33$, the prime vertices are $3,5,7,11,13,17,19,23$.
The graph has $15$ edges and cycle rank $8$. Its four filled triangles
come from $(3,7,23)$, $(3,11,19)$, $(3,13,17)$, and $(5,11,17)$,
leaving four independent homology classes. These are precisely the
four displayed cycles obtained from the faces of $F_{33,3}$ by
recording successive fixed labels. In particular, the upper face
gives the prime-label cycle $3\to17\to11\to3$; this triangle is
unfilled because $3+17+11\ne33$.

A partition $2p+q$ with $p\ne q$ detects a class through $x_p\varepsilon_p\varepsilon_q$
exactly when $\{p,q\}$ lies on a graph cycle, since all incoming boundaries
vanish on this coordinate. In weight $19$ the partitions are $(3,3,13)$,
$(3,5,11)$, and $(5,7,7)$. Both repeated-part edges are pendant, and
$\beta_{2,19,3}=0$.

\subsection{Common shifts and independent classes}
\label{sec:translations}
Let $G$ now be the binary Goldbach ring on a bounded-below support
$\A$, and use $K=K(G)$ for its coefficient Koszul complex.
The upper quadrilateral of Example~\ref{ex:weight33} is one instance
of a construction that works in every monomial degree.
Suppose $\sigma=\{a_1<\cdots<a_d\}$ and $\sigma+h$ both lie in the support
$\A$, for a positive integer $h$. The two lists need not be consecutive
terms of an arithmetic progression. For the binary ring consider
\begin{equation}\label{eq:translated-cycle}
 z_{\sigma,h}=\sum_{j=1}^d(-1)^{j-1}x_{a_j+h}
       \varepsilon_{\sigma\setminus\{a_j\}}.
\end{equation}

\Needspace{8\baselineskip}
\begin{theorem}[Independent translation classes]\label{thm:translation}
Fix $h>0$ and $d\ge2$. As $\sigma$ varies over the $d$-element subsets of
$\A\cap(\A-h)$, the chains $z_{\sigma,h}$ give linearly independent classes
over every field. For each $\sigma$ the class has homological degree $d-1$
and bidegree $(n,d)$, where $n=h+\sum_{a\in \sigma}a$. In particular,
\[
 \beta_{d-1,n,d}\ \ge\
 \#\{\sigma\subseteq\A\cap(\A-h): |\sigma|=d,\ h+\textstyle\sum_{a\in \sigma}a=n\}.
\]
These classes remain independent under enlargement of the support.
\end{theorem}
\begin{proof}
For $j<r$, the terms in the differential omitting both exterior
indices have opposite signs and coefficients $x_{a_j+h}x_{a_r}$
and $x_{a_j}x_{a_r+h}$. Their pair sums agree, so their
difference vanishes in $G$.

Suppose that a nonzero finite combination of the displayed cycles
were a boundary. Let $a$ be the largest label in any $\sigma$ with nonzero
coefficient, and choose such a $\sigma$ containing $a$. The coordinate
$x_{a+h}\varepsilon_{\sigma\setminus\{a\}}$ determines $\sigma$ for the fixed
$h$ and therefore has nonzero coefficient. Any preimage has monomial degree
and homological degree $d$, so its coefficients are scalars. The only
column giving the chosen coordinate is $\mathrm{d}(\varepsilon_\tau)$, where $\tau=(\sigma\setminus\{a\})\cup\{a+h\}$.
Since $d\ge2$, this column has another term with $a+h$ in its
exterior part. No selected cycle has such a term, and no other column
can cancel it: its coefficient label and its exterior labels together
determine $\tau$. The variables of $G$ are independent in monomial degree one,
so this is a contradiction. The proof uses only coefficients $\pm1$
and applies in every larger support.
\end{proof}

\subsection{Consequences for projective dimension}
The independence theorem turns repeated additive differences into
lower bounds for the size and length of a resolution.

\begin{corollary}\label{cor:density}
If the positive differences in $\A$ have unbounded multiplicities,
the binary Goldbach ring on $\A$ has infinite projective dimension over
$R$. In fact each positive homological degree of the linear
strand then has infinitely many generators. Both statements apply to
the positive odd primes.
\end{corollary}
\begin{proof}
An overlap of size $r$ at difference $h$ gives $\binom{r}{i+1}$
independent classes in homological degree $i$. For each fixed $i$ this
becomes arbitrarily large as $r$ increases. On the odd primes, put
$\A_m=\{p\le m:p\text{ odd prime}\}$. Its $\binom{|\A_m|}{2}$
unordered pairs have fewer than $m$ distinct positive differences. Some
difference therefore occurs at least $\binom{|\A_m|}{2}/m$ times,
which tends to infinity by the prime number theorem. Alternatively,
the prime progressions of Green and Tao \cite{GreenTao08} give
arbitrarily large overlaps by shifting each progression one step.
\end{proof}

Appendix~\ref{app:scroll-shift} gives a stronger bound for each fixed
shift using a scroll resolution. In Section~\ref{sec:betti-euler-shifts}
we fix the Goldbach degree instead and select independent classes
whose positive shifts may vary.

\section{Prime representations and shifted Koszul cycles}
\label{sec:betti-euler-shifts}

Theorem~\ref{thm:fixed-weight-betti-euler} gives an algebraic estimate
for the linear Betti number in terms of ordered representation counts,
for any positive support.
Corollary~\ref{cor:prime-betti-singular-series} combines this estimate
with known prime asymptotics. Theorem~\ref{thm:prime-shift-proportion}
then constructs independent common-shift cycles spanning a positive
fraction of the linear strand at each sufficiently large weight of
the appropriate parity. We then deduce a Hilbert-series consequence
from the smaller Betti numbers in lower homological degrees.
For a positive support $\A$, let $u$ mark Goldbach degree and put
\[
 P(u)=\sum_{a\in\A}u^a,\qquad T_d(n)=[u^n]P(u)^d .
\]
Thus $T_d(n)$ counts \emph{ordered} representations of $n$ by $d$
labels in the specified support, with repetitions allowed.

\subsection{The fixed-weight Betti estimate}
Consider any $k$-ary Goldbach ring $G$ on a positive support $\A$.
Translation of labels as in Proposition~\ref{setup:affine-invariance}
reduces a bounded-below support to this case.

\begin{theorem}[Fixed-weight Betti estimate]
\label{thm:fixed-weight-betti-euler}
For every field and every $d\ge k\ge2$, put $i=d-k+1$. Then
\begin{equation}\label{eq:fixed-weight-leading}
 \beta_{i,n,d}
 =\frac{\binom{d-1}{i-1}}{d!}\,T_d(n)
      +O_d\bigl((1+|\A_n|)^{d-2}\bigr).
\end{equation}
\end{theorem}
\begin{proof}
Let $K=K(G)=\mathbb K_\kk\otimes_R G$ be the Koszul complex with
coefficients in $G$, and write $K_{i,n,d}$ for its term in homological
degree $i$ and bidegree $(n,d)$. Its differential is
$\mathrm d_{i,n,d}:K_{i,n,d}\to K_{i-1,n,d}$.
Since $\operatorname{im}\mathrm d_{i+1,n,d}\subseteq\ker\mathrm d_{i,n,d}$,
the rank--nullity theorem gives
\begin{equation}\label{eq:koszul-rank-nullity}
 \beta_{i,n,d}
 =\dim_\kk\bigl(K_{i,n,d}/\operatorname{im}\mathrm d_{i+1,n,d}\bigr)
  -\dim_\kk\operatorname{im}\mathrm d_{i,n,d}.
\end{equation}
We count the first term by grouping monomials modulo incoming
boundaries, then bound the second by the dimension of the target.

Recall that $i=d-k+1$ and $\deg\varepsilon_a=(a,1)$.
A term $f\varepsilon_\sigma$ in $K_{r,n,d}$, with $|\sigma|=r$,
therefore has coefficient bidegree
$\deg f=(n-\sum_{a\in\sigma}a,d-r)$. In the three relevant terms,
the coefficient monomial degrees are
\[
 d-(i+1)=k-2,\qquad d-i=k-1,\qquad d-(i-1)=k.
\]
Since $I$ is generated in monomial degree $k$, its components below
that degree vanish. Thus the coefficient spaces in $K_{i+1,n,d}$
and $K_{i,n,d}$ agree with those in the polynomial-ring Koszul complex.

For $x^\alpha\in\mathcal M_{n,d}$, put
$S=\operatorname{supp}\alpha$ and $s_\alpha=|S|$. For $\tau\subseteq S$,
write $u_\tau=(x^\alpha/x_\tau)\varepsilon_\tau$.
The $\binom{s_\alpha}{i}$ terms $u_\tau$ with $|\tau|=i$ form the
group belonging to $x^\alpha$. For
$\eta=\{a_1<\cdots<a_{i+1}\}\subseteq S$, the incoming column is
\[
 \mathrm d u_\eta
   =\sum_{j=1}^{i+1}(-1)^{j-1}u_{\eta\setminus\{a_j\}},
\]
the oriented simplex boundary on $S$.
To compute its rank, fix $a_0=\min S$ and take the boundaries
$\mathrm d u_{\{a_0\}\cup\tau}$ for the $i$-subsets
$\tau\subseteq S\setminus\{a_0\}$. These are independent:
each has coefficient $1$ at its own $u_\tau$, and all its other
exterior subsets contain $a_0$. They also span. If an $(i+1)$-subset $\eta$
omits $a_0$, the simplex identity $\partial^2=0$ on
$\{a_0\}\cup\eta$ expresses its boundary as a signed sum of
boundaries of faces containing $a_0$.
There are $\binom{s_\alpha-1}{i}$ selected boundaries, so this is
the incoming rank over every field. Consequently $x^\alpha$ contributes
\[
 \binom{s_\alpha}{i}-\binom{s_\alpha-1}{i}
   =\binom{s_\alpha-1}{i-1}
\]
dimensions modulo incoming boundaries. Summing gives
\[
 \dim_\kk(K_{i,n,d}/\operatorname{im}\mathrm{d}_{i+1,n,d})
   =\sum_{x^\alpha\in\mathcal M_{n,d}}
             \binom{s_\alpha-1}{i-1}.
\]
We take out-of-range binomial coefficients to be zero.

Let $T_d^{\mathrm{dist}}(n)$ count ordered representations of $n$
by $d$ pairwise distinct labels in $\A$, and put
$T_d^{\mathrm{rep}}(n)=T_d(n)-T_d^{\mathrm{dist}}(n)$.
Each squarefree monomial has exactly $d!$ orderings and has
$s_\alpha=d$. The combined squarefree contribution is therefore exactly
\[
 \frac{\binom{d-1}{i-1}}{d!}\,T_d^{\mathrm{dist}}(n).
\]
In the binary case $i=d-1$, each squarefree monomial has $d$ terms,
corresponding to the choice of a single coefficient label. One incoming boundary relates them,
leaving $d-1$ dimensions and the factor $(d-1)/d!$.

For a representation with repetitions, choose two equal positions
and fill the other $d-2$ positions with labels from $\A_n$.
The equation $2a+\sum_{\text{other positions}}a_j=n$ determines
the repeated label $a$ uniquely, if it exists. Hence
\[
 T_d^{\mathrm{rep}}(n)
 \le\binom d2(1+|\A_n|)^{d-2}.
\]
Every nonsquarefree monomial contributes a number of dimensions
bounded in terms of $d$, and their number is at most
$T_d^{\mathrm{rep}}(n)$. Thus both their contribution and the
correction from replacing $T_d^{\mathrm{dist}}(n)$ by $T_d(n)$
are absorbed in the stated error:
\[
 \begin{aligned}
 \dim_\kk(K_{i,n,d}/\operatorname{im}\mathrm{d}_{i+1,n,d})
 &=\frac{\binom{d-1}{i-1}}{d!}\,T_d^{\mathrm{dist}}(n)
       +O_d\bigl(T_d^{\mathrm{rep}}(n)\bigr)\\
 &=\frac{\binom{d-1}{i-1}}{d!}\,T_d(n)
       +O_d\bigl((1+|\A_n|)^{d-2}\bigr).
 \end{aligned}
\]

It remains to bound the second term in \eqref{eq:koszul-rank-nullity}.
In the target $K_{i-1,n,d}$, choosing the $i-1$ exterior labels fixes the residual
coefficient weight. Its monomial degree is $k$, and all represented
monomials of that degree and weight agree in $G$. Thus each exterior set gives
at most one target basis term, and
\[
 0\le\dim_\kk\operatorname{im}\mathrm d_{i,n,d}
 \le\dim_\kk K_{i-1,n,d}\le\binom{|\A_n|}{i-1}.
\]
Since $i-1=d-k\le d-2$, this loss fits the same error bound.
This proves \eqref{eq:fixed-weight-leading} for the Betti number.

\end{proof}

\subsection{Prime asymptotics and the cubic correction}

Theorem~\ref{thm:fixed-weight-betti-euler} uses only the Koszul complex
and elementary counting, and applies to finite as well as infinite
positive supports. For the full odd-prime support, the known
asymptotic for $T_d(n)$ makes its error relatively small when $d\ge3$
is fixed and $n\to\infty$ through $n\equiv d\pmod2$.
The following proof obtains the needed prime-counting asymptotic from
Tao and Ter\"av\"ainen's theorem \cite[Theorem~1.6]{TT}, then applies
the algebraic estimate to obtain the Betti asymptotic. The qualitative
linear-forms theorem comes from Green--Tao \cite{GT10,GT12} and
Green--Tao--Ziegler \cite{GTZ12}; the cited quantitative version
provides the uniform estimate used here and in the completion counts
below.

\begin{corollary}[Singular-series asymptotics]\label{cor:prime-betti-singular-series}
Let $G$ be the binary Goldbach ring on all odd primes, and fix $d\ge3$.
Define the $d$-prime singular series by
\begin{equation}\label{eq:d-prime-singular-series}
 \mathfrak S_d(n)=
 \prod_{\substack{s\text{ prime}\\s\mid n}}
       \left(1+\frac{(-1)^d}{(s-1)^{d-1}}\right)
 \prod_{\substack{s\text{ prime}\\s\nmid n}}
       \left(1-\frac{(-1)^d}{(s-1)^d}\right).
\end{equation}
As $n\to\infty$ through $n\equiv d\pmod2$, uniformly over coefficient fields,
\begin{equation}\label{eq:prime-betti-singular-asymptotic}
 \beta_{d-1,n,d}
 \sim\frac{\mathfrak S_d(n)}{d!(d-2)!}
             \frac{n^{d-1}}{(\log n)^d}.
\end{equation}
\end{corollary}
\begin{proof}
Apply the theorem of Tao and Ter\"av\"ainen
\cite[Theorem~1.6]{TT} to
$a_1,\ldots,a_{d-1},n-\sum_j a_j$ on the positive simplex
$a_j>0$, $\sum_j a_j<n$. Their linear parts are pairwise
nonproportional because $d\ge3$, and the affine constant $n$ is
allowed at scale $n$. At a prime $s$, the normalized local factor is
\[
 \frac{s}{(s-1)^d}
 \#\{(u_1,\ldots,u_d)\in(\mathbb F_s^\times)^d:
                                      u_1+\cdots+u_d=n\},
\]
which equals the corresponding factor in
\eqref{eq:d-prime-singular-series}. The factor at $2$ is $2$ on the
stated parity; the odd-prime factors give
$0<c_d\le\mathfrak S_d(n)\le C_d$ uniformly in $n$.
The simplex has volume $n^{d-1}/(d-1)!$, so the theorem yields
\[
 R_d(n):=\sum_{\substack{u_1+\cdots+u_d=n\\u_j\ge1}}
              \prod_{j=1}^d\Lambda(u_j)
 =\frac{\mathfrak S_d(n)}{(d-1)!}\,n^{d-1}+o_d(n^{d-1}).
\]
Tuples involving a proper prime power contribute at most
$O_d(n^{d-3/2}(\log n)^{d+1})=o(n^{d-1})$ to this weighted count.
Discard also tuples with an entry at most $n/(\log n)^{d+1}$.
Their unweighted number is $O_d(n^{d-1}/(\log n)^{d+1})$, and their
total weight is $O_d(n^{d-1}/\log n)$. Every remaining prime entry
is odd for large $n$ and has logarithm
$\log n+O_d(\log\log n)$. Hence
\begin{equation}\label{eq:d-prime-count-asymptotic}
 T_d(n)\sim\frac{\mathfrak S_d(n)}{(d-1)!}
                    \frac{n^{d-1}}{(\log n)^d}.
\end{equation}
Finally, $\pi(n)\ll n/\log n$ gives
\[
 \frac{(1+\pi(n))^{d-2}}{T_d(n)}
       \ll_d\frac{(\log n)^2}{n}\longrightarrow0.
\]
Theorem~\ref{thm:fixed-weight-betti-euler}, with coefficient
$(d-1)/d!$ for the binary ring, proves
\eqref{eq:prime-betti-singular-asymptotic}.
\end{proof}

For the binary ring on a positive support $\A$, at $d=3$ we can read
the correction directly from the cubic graph of
Example~\ref{ex:weight33}. Let $v_n$ and $c_n$ count its vertices and
components, where the vertices are the labels occurring in a triple
and $a\ne b$ are joined when $n-a-b\in\A$. A triple with three, two,
or one distinct labels contributes $2,1,0$ to
$\dim_\kk(K_{2,n,3}/\operatorname{im}\mathrm{d}_{3,n,3})$. Its outgoing
rank loss is the incidence rank $v_n-c_n$. Applying
\eqref{eq:koszul-rank-nullity} gives
\begin{equation}\label{eq:cubic-betti-euler}
 \beta_{2,n,3}
  =\frac{T_3(n)-\mathbf1_{\{\,n/3\in\A\,\}}}{3}-v_n+c_n.
\end{equation}
For example at weight $33$ there are $34$ ordered triples, including
$(11,11,11)$, and the graph is connected on eight vertices. The formula
gives $(34-1)/3-8+1=4$, as in the drawing.

\subsection{How many classes do common shifts give?}
Corollary~\ref{cor:prime-betti-singular-series} estimates the dimension
of the linear strand. We now identify explicit classes accounting for
a positive fraction of it. Theorem~\ref{thm:prime-shift-proportion}
selects an independent family from the common-shift cycles of
Section~\ref{sec:translations}, using estimates for simultaneous prime
patterns to count these classes.

Take $G$ to be the binary ring on the odd primes $\Podd$. Let
$W^+_{n,d}$ be the span in $\Tor^R_{d-1}(G,\kk)_{n,d}$ of the classes
$[z_{\sigma,h}]$ of Theorem~\ref{thm:translation}, for all $d$-element sets
$\sigma$ and all $h>0$ such that $\sigma,\sigma+h\subseteq\Podd$ and
$\sum_{a\in \sigma}a+h=n$. Different shifts can give dependent classes.
For example at weight $33$ the three possible positive shifts arise from
\[
 (\sigma,h)=(\{3,7,13\},10),\quad(\{3,7,19\},4),\quad(\{3,11,17\},2).
\]
These exhaust the possibilities by enumerating the prime triples of
weight $33$ in Example~\ref{ex:weight33}. In the first two cycles the
terms $x_{13}\varepsilon_7\varepsilon_{13}$ and
$x_7\varepsilon_7\varepsilon_{19}$ each occur only in their indicated
cycle, and no incoming Koszul boundary uses either term. The last
cycle is nonzero by Theorem~\ref{thm:translation}. Thus the shifts give
three of the four independent classes at this weight.

\begin{theorem}[Positive shifts at a fixed weight]
\label{thm:prime-shift-proportion}
For each $d\ge3$ there are constants $c_d,C_d>0$ depending only on $d$ such that, over every
field and for every sufficiently large $n\equiv d\pmod 2$,
\[
 c_d\,\frac{n^{d-1}}{(\log n)^d}
 \ \le\ \dim_\kk W^+_{n,d}\ \le\ \beta_{d-1,n,d}
 \ \le\ T_d(n)\ \le\ C_d\,\frac{n^{d-1}}{(\log n)^d}.
\]
In particular the positive shifts span at least a fixed positive fraction
of the linear strand in each such weight.
\end{theorem}
\begin{proof}
Fix $n\equiv d\pmod 2$; all intervals and completions below depend on this $n$.
Divide $[n/(4d),n/(3d)]$ into $d-1$ equal pieces, numbered from left
to right, and let $J_j$ be the middle half of the $j$th piece. Each
$J_j$ has length $n/(24d(d-1))$. We choose one odd prime $p_j\in J_j$
for $1\le j\le d-1$, so $\boldsymbol p=(p_1,\ldots,p_{d-1})$ is
strictly increasing. Set $J_0=[n/(8d),n/(6d)]$, to the left of $J_1$.

Call an odd prime $p_0\in J_0$ a \emph{completion} of $\boldsymbol p$ at weight $n$
if, with $h=n-\sum_{j=0}^{d-1}p_j$, every $p_j+h$ is prime
($0\le j\le d-1$). A completion gives the $d$ primes
$\sigma=\{p_0<p_1<\cdots<p_{d-1}\}$, and shifting them all by $h$ gives
the $d$ primes $\sigma+h$, so $z_{\sigma,h}$ is a translation cycle of
weight $n$. The smallest shifted prime is
\[
 p_0+h=n-\sum_{j=1}^{d-1}p_j .
\]
Thus, with $n$ fixed, $p_0+h$ depends only on $\boldsymbol p$, and it exceeds $2n/3$.
For each tuple admitting a completion, choose its least such $p_0$.
This choice determines $h$; the shift is not chosen independently and
may vary with the tuple.

The selected cycles are linearly independent, although different
tuples may have different shifts.
In the cycle for $\boldsymbol p$, consider the term
$x_{p_0+h}\varepsilon_{p_1}\cdots\varepsilon_{p_{d-1}}$. Every other term
of this cycle has $p_0\in J_0$ in its exterior part. Thus each selected
cycle has just one term whose exterior labels all lie in
$J_1,\ldots,J_{d-1}$, and those labels are the entries of $\boldsymbol p$.
Two tuples with the same sum share the coefficient label $p_0+h$, but
their distinguished terms still differ in the exterior part.
The least-completion rule then determines $p_0$, $h$, and the selected cycle.
The sole possible incoming Koszul column
containing it is $\mathrm{d}(\varepsilon_\tau)$ for
$\tau=\{p_1<\cdots<p_{d-1}<p_0+h\}$, since the incoming coefficients have monomial degree zero. The same column has a term with $p_0+h$ in
the exterior part. No selected cycle has such a term, and the
coefficient label and exterior indices determine this incoming column
uniquely. Thus in any boundary equal to a combination of the selected
cycles the coefficient of this column must vanish. Reading the coefficient
of the distinguished term proves the claimed independence over every field.

Let $\mathcal T_+$ be the set of tuples $\boldsymbol p$ admitting
a completion, and let $f(\boldsymbol p)$ count all their completions.
Write
\[
 \langle f,1\rangle=\sum_{\boldsymbol p\in\mathcal T_+}f(\boldsymbol p),
 \qquad
 \langle f,f\rangle=\sum_{\boldsymbol p\in\mathcal T_+}f(\boldsymbol p)^2.
\]
These count completions and ordered pairs of completions of a common
tuple, respectively, including repeated pairs. By Cauchy--Schwarz,
$|\mathcal T_+|\ge\langle f,1\rangle^2/\langle f,f\rangle$.
We will prove
\begin{equation}\label{eq:shift-moments}
 \langle f,1\rangle\gg_d\frac{n^d}{(\log n)^{2d}},\qquad
 \langle f,f\rangle\ll_d\frac{n^{d+1}}{(\log n)^{3d}}.
\end{equation}
The selected independent classes then give
\[
 \dim_\kk W^+_{n,d}\ge|\mathcal T_+|
 \ge\frac{\langle f,1\rangle^2}{\langle f,f\rangle}
 \gg_d\frac{n^{d-1}}{(\log n)^d}.
\]

Both moments count simultaneous prime values of linear forms on a
box. We use Tao and Ter\"av\"ainen \cite[Thm.~1.6]{TT}: for bounded,
pairwise nonproportional slopes and constant terms $O(n)$, the
von Mangoldt weighted count has a uniform asymptotic whose main term
is the box volume times the product of local factors.
With $D$ free variables and $t$ forms, its error here is
$O_d(n^D(\log\log n)^{-\eta_d})$ for some $\eta_d>0$.
We use $(D,t)=(d,2d)$ and $(d+1,3d)$.
Since all forms on our boxes are comparable to $n$, converting weighted
bounds to prime counts divides both the main-term and error scales by
$(\log n)^t$. Thus the saving in $\log\log n$ suffices.

\emph{The first bound.} Use one integer variable $a_j\in J_j$ for each
$0\le j\le d-1$; the index $j$ specifies the interval in which the variable
ranges. A point $\boldsymbol a=(a_0,\ldots,a_{d-1})$ of
$J_0\times\cdots\times J_{d-1}$ gives a tuple and a completion, with
$p_j=a_j$, precisely when all $2d$ forms
\[
 a_j,\qquad n-\sum_{\substack{0\le r\le d-1\\r\ne j}}a_r
 \qquad (0\le j\le d-1)
\]
take prime values; the second family is $a_j+h$. Thus $\langle f,1\rangle$
is the number of these points.
Their slopes are pairwise nonproportional because $d\ge3$.
At an odd prime $s$ choose nonzero residues for the $a_j$ whose sum
is $n$ modulo $s$. This is possible since every residue is the sum
of two nonzero residues. The shifted forms then also have nonzero
residues. At $2$ take all $a_j$ odd and use $n\equiv d\pmod 2$.
These choices bound the finitely many small local factors below
uniformly in $n$. At all sufficiently large primes $s$ the slopes remain
pairwise independent modulo $s$. Counting intersections of pairs
of their zero hyperplanes gives normalized local factors
$1+O_d(s^{-2})$. Hence the product of the local factors is uniformly
bounded below. Our box has volume $\gg_d n^d$, and all forms there
are comparable to $n$. The contribution of proper prime powers is
$o(n^d)$ in the weighted count: one form takes at most
$O(\sqrt n\log n)$ such values, each on $O_d(n^{d-1})$ points.
Including the $2d$ logarithmic weights gives a total of at most
$O_d(n^{d-1/2}(\log n)^{2d+1})=o_d(n^d)$.
Converting logarithmic weights to an unweighted count gives the first
bound in \eqref{eq:shift-moments}.

\emph{The second bound.} A pair of completions $(p_0,p'_0)$ of a common
tuple is a point of $J_0\times J_0\times J_1\times\cdots\times J_{d-1}$:
two integer variables $a_0,a'_0\in J_0$ and the same $a_j\in J_j$ for
$1\le j\le d-1$. Put $a_d=n-\sum_{j=1}^{d-1}a_j$, which takes the value
$p_0+h$ at $\boldsymbol p$. The point is counted by $\langle f,f\rangle$
precisely when the forms
\[
 a_d,\ a_0,\ a'_0,\quad
 a_j,\ a_d-a_0+a_j,\ a_d-a'_0+a_j \qquad(1\le j\le d-1)
\]
are prime. There are $d+1$ free variables because the two completions share the
$d-1$ exterior labels. Their prime requirements also share those labels
and $a_d$, leaving $4d-d=3d$ distinct forms. Their slopes are pairwise
nonproportional. Each of the last two families uses its own completion
variable, and the form indexed by $j$ omits precisely $a_j$.
Coincident completions $a_0=a'_0$ are included: equality of some values
at such points does not change these linear parts, and the weighted
product still contains all $3d$ logarithmic factors. The same pairwise intersection
estimate bounds the product of the local factors above uniformly.
Each prime value in this box has logarithmic weight comparable to
$\log n$, so the linear forms theorem gives the second bound in
\eqref{eq:shift-moments}.

The entire Koszul module in homological degree $d-1$ has at most $T_d(n)$
basis terms, since placing the coefficient label first and the
increasing exterior labels after it gives an ordered representation.
The upper bound $T_d(n)\ll_d n^{d-1}/(\log n)^d$ follows from
\eqref{eq:d-prime-count-asymptotic} and the uniform upper bound for
$\mathfrak S_d(n)$.
\end{proof}

\subsection{A Hilbert-series consequence}
\label{sec:hilbert-euler-consequence}
The preceding estimates concern individual Betti numbers and explicit
cycles. Their alternating sum connects this growth with the Hilbert
series. For a $k$-ary Goldbach ring on a positive support, define
\[
 \gamma_{n,d}=\sum_{i\ge0}(-1)^i\beta_{i,n,d},\qquad
 E_G(u,v)=\sum_{n,d}\gamma_{n,d}u^nv^d.
\]
Let $H_G(u,v)=\sum_{n,d}(\dim_\kk G_{n,d})u^nv^d$ and let $H_R$
be the corresponding series of $R$. Taking the Euler characteristic
of the Koszul complex in each bidegree gives
\begin{equation}\label{eq:betti-euler-hilbert}
 E_G(u,v)=\frac{H_G(u,v)}{H_R(u,v)}
        =H_G(u,v)\prod_{a\in\A}(1-u^av).
\end{equation}
Each coefficient is a finite sum. By the component basis in
Proposition~\ref{setup:components}, $H_G$ is independent of the
coefficient field, and hence so are $E_G$ and $\gamma_{n,d}$.

\begin{corollary}[Comparison with the Euler coefficient]
\label{cor:hilbert-euler-comparison}
For every field and every $d\ge k\ge2$, put $i=d-k+1$. Then
\[
 \bigl|(-1)^i\gamma_{n,d}-\beta_{i,n,d}\bigr|
 \le\sum_{j<i}\beta_{j,n,d}
 =O_d\bigl((1+|\A_n|)^{i-1}\bigr).
\]
Thus $(-1)^i\gamma_{n,d}$ has the same leading term as
$\beta_{i,n,d}$ in Theorem~\ref{thm:fixed-weight-betti-euler}.
For the binary ring on all odd primes, fixed $d\ge3$, and
$n\to\infty$ through $n\equiv d\pmod2$,
\[
 (-1)^{d-1}\gamma_{n,d}\sim\beta_{d-1,n,d}
 \sim\frac{d-1}{d!}T_d(n).
\]
\end{corollary}
\begin{proof}
Write $K=K(G)$ for the Koszul complex, with the conventions of
Theorem~\ref{thm:fixed-weight-betti-euler}.
In homological degrees $j>i$ the Koszul differential and both adjacent terms
in this bidegree agree with the exact polynomial-ring complex in
positive homological degrees, so $\beta_{j,n,d}=0$. Therefore
\[
 (-1)^i\gamma_{n,d}=\beta_{i,n,d}
       +\sum_{j<i}(-1)^{i+j}\beta_{j,n,d}.
\]
For $j<i$ the coefficients have monomial degree
$d-j\ge k$. Choose the $j$ exterior labels and $d-j-k$ coefficient
labels. The remaining $k$ labels have prescribed sum and contribute
at most one coefficient class. These choices span $K_{j,n,d}$ and require
choosing only $j+(d-j-k)=d-k=i-1$ labels. Hence
\[
 \beta_{j,n,d}\le\dim_\kk K_{j,n,d}
   =O_d\bigl((1+|\A_n|)^{i-1}\bigr).
\]
There are only finitely many such homological degrees, bounded in number by $d$.
It follows that
\[
 |(-1)^i\gamma_{n,d}-\beta_{i,n,d}|
 \le\sum_{j<i}\beta_{j,n,d}
 =O_d\bigl((1+|\A_n|)^{i-1}\bigr).
\]
Since $i-1=d-k\le d-2$, Theorem~\ref{thm:fixed-weight-betti-euler}
gives the asserted leading term. For binary odd-prime support,
Corollary~\ref{cor:prime-betti-singular-series} makes the error
negligible relative to that term.
\end{proof}

For a binary ring, the ideal is generated in monomial degree two
and $\Tor_0^R(G,\kk)$ has monomial degree zero. Consequently
$\gamma_{n,3}=\beta_{2,n,3}$, so \eqref{eq:cubic-betti-euler}
also computes this Euler coefficient exactly. In monomial degree four,
$-\gamma_{n,4}=\beta_{3,n,4}-\beta_{2,n,4}$.
The identity \eqref{eq:betti-euler-hilbert} therefore recovers
the dominant Betti growth from Hilbert-series coefficients, while
retaining the cancellations between homological degrees.
The same Euler-characteristic principle compares binary and ternary
rings through mixed Tor in \cite{CompanionTopology}.

\section{The universal relation group and interpolation}
\label{sec:depth}

Binary Goldbach moves determine both a monoid of monomial classes and
its group completion. We compute the latter on the odd primes and
deduce the interpolation formula of Theorem~\ref{main:prime-interpolation}.
We then prove that the finite relation groups stabilize, and use
localization to bound the depth and projective dimension of the two
truncations.

For the binary Goldbach ring $G=R/I$ on a finite support $\A$, let
$e_a$ denote the unit vector at $a$ in the free abelian group
$\mathbb Z^{\A}$. Set
\[
 L_{\A}=\langle e_a+e_b-e_{a'}-e_{b'}:a+b=a'+b',\ a,b,a',b'\in\A\rangle,
 \qquad U_{\A}=\mathbb Z^{\A}/L_{\A},
 \qquad r_{\A}=\operatorname{rank}_{\mathbb Z}U_{\A}.
\]
This is the universal ambient group for Freiman homomorphisms of
order two, as constructed by Konyagin--Lev \cite{KonyaginLev00}
and treated categorically by Blanco--Haghverdi \cite{BlancoHaghverdi26}.
Explicitly, every map $g:\A\to H$ to an abelian group that preserves
equal pair sums extends uniquely to a group homomorphism
$U_{\A}\to H$. We use this universal group to relate arithmetic
identities to localization and depth.

Write $[a]$ for the image of $e_a$. A monomial
$x_{\boldsymbol a}=x_{a_1}\cdots x_{a_d}$ determines the element
$[a_1]+\cdots+[a_d]$ of $U_{\A}$. The classes of monomials modulo binary moves form a commutative
monoid under multiplication. Its monoid algebra is $G$, by
Proposition~\ref{setup:components}. The group $U_{\A}$ is
its group completion, in which the class of each variable has an inverse.

\subsection{The integral group of prime relations}

For an infinite support $\A$, the same construction uses finitely
supported integer exponent vectors. Write $U_\infty=U_{\Podd}$ for the
group on all odd primes. Its relations are
$[p]+[q]=[p']+[q']$ whenever $p+q=p'+q'$. To compute
this group, first omit $3$. In a pair of primes greater
than three, the number $j$ of entries congruent to one modulo three
is determined by the sum, since the three possibilities $j=0,1,2$ give
distinct values of $2j-2\bmod3$. For example, $5+19=7+17$ replaces one prime in each residue class
by another in the same class. Every binary move among primes greater
than three preserves both counts. To prove the converse in the group, we construct
pair-sum identities among primes and cancel the auxiliary generators.

\begin{lemma}\label{depth:prime-pair-conversion}
Let $\boldsymbol p=(p_1,\ldots,p_d)$ and $\boldsymbol q=(q_1,\ldots,q_d)$
be lists of primes in $\Podd_{>3}=\Podd\setminus\{3\}$ with the same sum
and the same counts in each residue class modulo three. There is a
(possibly empty) list $\boldsymbol r$ of primes in $\Podd_{>3}$ such
that binary moves on this support connect $x_{\boldsymbol r}x_{\boldsymbol p}$
to $x_{\boldsymbol r}x_{\boldsymbol q}$. In particular,
$\sum_i[p_i]=\sum_i[q_i]$ in $U_{\Podd_{>3}}$.
\end{lemma}
\begin{proof}
For $d\le1$ the monomials coincide. For $d=2$ the two monomials
either coincide or are joined by a single binary move.
Assume $d\ge3$, and reorder the lists so that
$p_i\equiv q_i\pmod3$ for each $i$.
Put $\Delta_0=0$ and $\Delta_i=\sum_{j=1}^i(p_j-q_j)$.
Each $\Delta_i$ is divisible by six, and the equality of the sums gives
$\Delta_d=0$.

\emph{The auxiliary forms.}
Choose integer parameters $a_1,\ldots,a_d,c$. With indices read cyclically
on $a_i$ and $b_i$, define
\[
 b_i=a_{i+1}+c+\Delta_i,\qquad r_i=a_i+b_i-p_i .
\]
For every choice of the parameters we have the integer identities
\[
 p_i+r_i=a_i+b_i,\qquad q_i+r_i=a_{i+1}+b_{i-1}.
\]
Indeed,
\[
 q_i+r_i
 =a_i+a_{i+1}+c+\Delta_i-(p_i-q_i)
 =a_{i+1}+(a_i+c+\Delta_{i-1}).
\]
The condition $\Delta_d=\Delta_0=0$ also verifies the cyclic case $i=1$.
Consequently, if all $3d$ values $a_i,b_i,r_i$ are primes greater
than three, the two lists can be converted, after adjoining the same
$r_i$, into the same list of $a_i$ and $b_i$. It remains to produce
such prime values simultaneously.

Let $\Psi$ be this system of $t=3d$ affine forms in $D=d+1$
parameters $\boldsymbol z=(a_1,\ldots,a_d,c)$. Its linear parts are
\[
 e_i,\qquad e_{i+1}+e_c,\qquad e_i+e_{i+1}+e_c
 \qquad(1\le i\le d).
\]
They are distinct nonzero $0$--$1$ vectors, since $d\ge3$.
Two such vectors can be proportional over a field only if they are
equal: a nonzero coordinate forces the proportionality factor to be
$1$. Thus the linear parts are pairwise independent both over
$\mathbb Q$ and over every $\mathbb F_s$.
Keeping $c$ free is essential for this argument. Fixing $c$ in advance
would make $b_i$ and $a_{i+1}$ have the same linear part.
The parameter $c$ need not be prime; at a successful tuple it is even,
because $c=b_i-a_{i+1}-\Delta_i$.

\emph{Local admissibility and the singular series.}
For each prime $s$, define the admissible residue density and its
normalization by
\[
 \alpha_s=s^{-D}
   \#\{\boldsymbol z\in\mathbb F_s^D:
                         \psi(\boldsymbol z)\ne0\text{ for every }\psi\in\Psi\},
 \qquad
 \beta_s=(1-1/s)^{-t}\alpha_s .
\]
The quantity $\alpha_s$ is a probability. For $t$ independent uniform
residues modulo $s$, the probability that all are nonzero is
$(1-1/s)^t$, since each residue is nonzero with probability $1-1/s$.
The ratio $\beta_s$ compares the actual density with this independent
model; unlike $\alpha_s$, the correction factor $\beta_s$ can exceed one.

We first show $\alpha_s>0$ for every prime $s$.
Modulo $2$ take $c=0$ and all $a_i=1$; every form is odd.
Modulo $3$ take $c=1$ and all $a_i=1$. Then $b_i=2$ and
$r_i=-p_i\ne0\bmod3$. At a prime $s\ge5$ set $c=0$ modulo $s$
and choose the $a_i$ successively. First choose $a_1\ne0$;
since $\Delta_d=0$, this also makes $b_d=a_1\ne0$.
For $2\le i\le d$, choose $a_i$ avoiding
\[
 0,\qquad -\Delta_{i-1},\qquad
 p_{i-1}-a_{i-1}-\Delta_{i-1},
\]
so that $a_i,b_{i-1},r_{i-1}$ are nonzero.
At the final choice also avoid $p_d-a_1$, ensuring
$r_d=a_d+a_1-p_d\ne0$. Each choice excludes at most four residues,
so it is possible even at $s=5$.

For the product of the $\beta_s$, we also need a bound as $s$ grows.
Each equation $\psi(\boldsymbol z)=0$ defines a hyperplane with
$s^{D-1}$ points. By pairwise independence of the linear parts, any
two of these hyperplanes intersect in exactly $s^{D-2}$ points.
The union bound and the first two inclusion--exclusion bounds therefore give
\[
 1-\frac ts\ \le\ \alpha_s\
 \le\ 1-\frac ts+\binom t2\frac1{s^2}.
\]
No independence assertion about three or more hyperplanes is needed.
It follows that
\[
 \alpha_s=1-\frac ts+O_d(s^{-2}),\qquad
 \beta_s=(1-1/s)^{-t}\alpha_s=1+O_d(s^{-2}).
\]
Thus the positive local factors have a convergent, positive product
$\mathfrak S_\Psi=\prod_s\beta_s$.
Here local admissibility supplies positivity at the finitely many
small primes, while the summable $O_d(s^{-2})$ error controls the tail.
The product normalization can also be seen directly from the Chinese
remainder theorem. For $W=\prod_{s\le w}s$,
\[
 \frac{\#\{\boldsymbol z\bmod W:
                  \gcd(\psi(\boldsymbol z),W)=1\text{ for all }\psi\in\Psi\}}{W^D}
       \left(\frac W{\varphi(W)}\right)^t
 =\prod_{s\le w}\beta_s .
\]
This finite congruence calculation determines the local factor in the
prime-counting theorem; it does not by itself produce prime values.

\emph{Prime values at scale $N$.}
Keep the two lists, hence all affine constants, fixed, and let
$N\to\infty$. Take
$\Omega_N=[N,2N]^{d+1}$. Once
$|\Delta_i|,|\Delta_i-p_i|\le N$ for every $i$, the forms satisfy
\[
 N\le a_i\le2N,\qquad N\le b_i\le5N,\qquad 2N\le r_i\le7N.
\]
In particular they are positive and comparable to $N$ throughout
this box, whose volume is $N^{d+1}$.
Apply Tao--Ter\"av\"ainen \cite[Theorem~1.6]{TT} at scale $2N$,
since $\Omega_N\subset[-2N,2N]^{d+1}$. The linear coefficients are
bounded, the affine constants are $O(N)$, and the linear parts are
pairwise nonproportional. With $\Lambda$ the von Mangoldt function,
the theorem gives, for some $\eta_d>0$,
\[
 \sum_{\boldsymbol z\in\Omega_N\cap\mathbb Z^{d+1}}
           \prod_{\psi\in\Psi}\Lambda(\psi(\boldsymbol z))
 =\mathfrak S_\Psi N^{d+1}
   +O_{\boldsymbol p,\boldsymbol q}
        \bigl(N^{d+1}(\log\log N)^{-\eta_d}\bigr).
\]
Only the qualitative $o(N^{d+1})$ error is needed here. That version
was already available from Green--Tao and Green--Tao--Ziegler
\cite{GT10,GT12,GTZ12}, as recalled after the cited theorem.

The weight $\Lambda$ also detects proper prime powers, which must
be discarded. A specified form can take only $O(\sqrt N\log N)$
such values in $[1,7N]$. For each value, one parameter is determined
by the other $d$, since the form has a coefficient $1$.
The union over the $3d$ forms, including the logarithmic weights,
therefore contributes at most
\[
 O_d\bigl(N^{d+1/2}(\log N)^{3d+1}\bigr)=o(N^{d+1}).
\]
Since $\mathfrak S_\Psi>0$, a positive weighted contribution remains
on tuples for which every form is prime. In fact their number is
\[
 \bigl(\mathfrak S_\Psi+o_{\boldsymbol p,\boldsymbol q}(1)\bigr)
                 \frac{N^{d+1}}{(\log N)^{3d}},
\]
because every prime value has logarithm $\log N+O(1)$.
Their proportion among all parameter tuples tends to zero; the
positive limiting quantity above is the normalized weighted count.
For sufficiently large $N$, all these primes exceed three.

\emph{The moves and cancellation.}
Choose one successful tuple and put $\boldsymbol a=(a_1,\ldots,a_d)$,
$\boldsymbol b=(b_1,\ldots,b_d)$ and $\boldsymbol r=(r_1,\ldots,r_d)$.
Replacing each pair $(p_i,r_i)$ by $(a_i,b_i)$ gives a path from
$x_{\boldsymbol r}x_{\boldsymbol p}$ to $x_{\boldsymbol a}x_{\boldsymbol b}$.
Replacing $(q_i,r_i)$ by $(a_{i+1},b_{i-1})$ gives a path from
$x_{\boldsymbol r}x_{\boldsymbol q}$ to the same monomial, by cyclic indexing.
Thus
\[
 x_{\boldsymbol r}x_{\boldsymbol p}
 \longrightarrow x_{\boldsymbol a}x_{\boldsymbol b}
 \longleftarrow x_{\boldsymbol r}x_{\boldsymbol q},
\]
where each arrow denotes $d$ binary moves. Passing to the group allows
us to cancel $\sum_i[r_i]$. Cancellation of the monomial
$x_{\boldsymbol r}$ in the Goldbach ring is not asserted.
\end{proof}

\paragraph{A cyclic conversion.}
For $\boldsymbol p=(5,5,23)$ and $\boldsymbol q=(11,11,11)$,
both sums are $33$ and every entry is $2\bmod3$.
Here $(\Delta_1,\Delta_2,\Delta_3)=(-6,-12,0)$.
One successful choice is
\[
 c=12,\qquad
 \boldsymbol a=(29,31,41),\qquad
 \boldsymbol b=(37,41,41),\qquad
 \boldsymbol r=(61,67,59).
\]
The identities used by the two paths are
\[
 \begin{array}{lll}
 5+61=29+37,&5+67=31+41,&23+59=41+41,\\
 11+61=31+41,&11+67=41+37,&11+59=29+41.
 \end{array}
\]
Both paths reach $x_{29}x_{31}x_{37}x_{41}^3$. Thus the construction
realizes the equality $[23]+2[5]=3[11]$ after the common auxiliary
classes $[59]+[61]+[67]$ are canceled in the group.

\begin{theorem}[The integral group of binary prime relations]
\label{depth:prime-pair-group}
The group on the primes greater than three is free of rank three.
On all odd primes the group $U_\infty$ is free of rank two,
with
\[
 [p]=[5]+\frac{p-5}{2}\bigl([7]-[5]\bigr)
 \qquad(p\in\Podd).
\]
\end{theorem}
\begin{proof}
Write $A=[5]$, $B=[7]$, and $D=[11]-[5]$. Every prime $p>3$ has
the form $s+6j$ with $s=5$ or $7$ and $j\ge0$.
Apply Lemma~\ref{depth:prime-pair-conversion} to $p$ together with $j$
copies of $5$, and $s$ together with $j$ copies of $11$. It gives
$[p]=[s]+jD$. For example, $5+17=11+11$ gives $[17]=[5]+2D$.
Thus $A,B,D$ generate the group on primes greater
than three. The residue counts distinguish $A$ and $B$. Both counts
vanish on $D$, which has weight six. These three elements are
therefore integrally independent.

Now adjoin $3$. The relations $3+7=5+5$ and $3+11=7+7$ give
$[3]=2A-B$ and $D=3(B-A)$. Thus the step between primes in the
same residue class is $3(B-A)$ per six units of weight, and
$[7]-[5]=B-A$ accounts for the two-unit step between $5$ and $7$.
Substituting in $[p]=[s]+jD$ gives the asserted formula for every
$p>3$, and $[3]=2A-B$ gives it also at $3$. The generators $A$ and
$B-A$ have independent vectors of Goldbach and monomial degrees
$(5,1)$ and $(2,0)$, so
they freely generate.
\end{proof}

To deduce Theorem~\ref{main:prime-interpolation}, observe that any map
$g:\Podd\to H$ preserving equal pair sums extends to a homomorphism
$U_\infty\to H$ by $[p]\mapsto g(p)$. Since
$[7]-[5]=[5]-[3]$, the formula of Theorem~\ref{depth:prime-pair-group}
becomes $[p]=[3]+((p-3)/2)([5]-[3])$. Applying the homomorphism gives
the interpolation formula in $H$.

For complex-valued maps the formula is equivalent to $g(p)=Ap+B$
on the odd primes. If $f(p+q)=g(p)+g(q)$ for odd primes $p,q$, it
also gives $f(n)=An+2B$ on $\Podd+\Podd$. This proves the odd-prime
form of the binary affine-classification conjecture in
\cite[Conj.~1]{KP15}. The corresponding theorem for sums of at least
three primes is proved there; Park \cite[p.~6, after Thm.~3.7]{Park25}
still presents the binary analogue as conjectural. Our integral
computation additionally allows arbitrary abelian targets, including
groups with torsion.

The odd-prime restriction must be distinguished from the all-primes
wording of \cite[Conj.~1]{KP15}. Every equal-pair-sum relation involving
$2$ is trivial: an odd sum $2+p$ has only the unordered prime
representation $\{2,p\}$, and $4$ has only $2+2$.
Thus, if $U_{\{2\}\cup\Podd}$ denotes the universal group on all primes,
\[
 U_{\{2\}\cup\Podd}=\mathbb Z[2]\oplus U_\infty
                    \cong\mathbb Z^3.
\]
The value $g(2)$ is independent of the values on the odd primes, so
pair-sum preservation alone cannot imply one affine formula on all
primes. For example, $g(2)=1$ and $g(p)=0$ for every odd prime $p$
preserves equal pair sums but is not affine on the full prime set.
The comparison with \cite[Conj.~1]{KP15} is therefore with its
odd-prime form. The shifted conjecture addressed next already asks
for an odd-prime conclusion.

\begin{corollary}[K\'atai--Phong interpolation]\label{depth:KP-interpolation}
Fix positive integers $a,b$. Let $f,g$ be functions on the positive
integers with values in an abelian group, and suppose that
\[
 f(p+q+a+b)=g(p+a)+g(q+b)
\]
for all odd primes $p,q$. Then for every odd prime $p$,
\[
 g(p+a)=\frac{p-3}{2}\,g(5+a)-\frac{p-5}{2}\,g(3+a).
\]
In particular this proves Conjecture~6 of K\'atai and Phong
\cite[p.~235]{KP15Problems}, whose hypothesis imposes the displayed
identity for all primes and whose conclusion concerns the odd primes.
\end{corollary}
\begin{proof}
Put $g_a(p)=g(p+a)$ and $g_b(p)=g(p+b)$. Interchanging $p$ and $q$ gives
$g_a(p)-g_b(p)=g_a(q)-g_b(q)$. Thus $g_a(p)+g_a(q)$ differs from $f(p+q+a+b)$
by a constant and depends only on $p+q$. Apply
Theorem~\ref{main:prime-interpolation} to $g_a$ and collect the
coefficients of $g_a(3)$ and $g_a(5)$.
\end{proof}

\subsection{Finite stabilization}

Put $\A_m=\{p\le m:p\text{ odd prime}\}$ and $U_m=U_{\A_m}$.

The existence of a relation involving only smaller primes was stated by
K\'atai--Phong \cite[Theorem~3]{KP15}: every sufficiently large prime
is a sum of two smaller primes minus a third smaller prime. We use this
arithmetic fact to obtain integral stabilization below. For completeness,
we give a short proof with balanced ranges from the classical almost-all
Goldbach estimate. The balanced form will also distinguish the two
truncations.

Goldbach questions also enter the numerical semigroups generated by
primes above a threshold: Hellus--Rechenauer--Waldi
\cite{HellusRechenauerWaldi20} study their Frobenius numbers and
minimal generators. The relations used here preserve both the sum
and the number of summands.

Write $N$ for the interval scale and $n$ for Goldbach degree. Set
\[
 \begin{aligned}
 T_{2,N}(n)&=\#\{(p,p'):p,p'\text{ odd primes},\ p+p'=n,\ N<p,p'\le2N\},\\
 R_{2,N}(n)&=\sum_{\substack{u+v=n\\N<u,v\le2N}}\Lambda(u)\Lambda(v),
 \end{aligned}
\]
where $u,v$ range over integers and $\Lambda$ is the von Mangoldt function.
Let $V_{2,N}(n)$ count all integer pairs in the same range, without weights,
and let $\mathfrak S(n)$ denote the Goldbach singular series.
The analytic estimate recalled in the introduction of
Bazzanella--Languasco \cite[eq.~(2)]{BL00}, with $H=1$ and the same scale $N$, gives the single-target estimate
$R_{2,N}(n)\sim V_{2,N}(n)\mathfrak S(n)$ for all but
$O_C(N/(\log N)^C)$ even integers $n\in(5N/2,7N/2]$, for every fixed $C>0$.
In this interval the first factor is $\gg N$
and the second is bounded below by a positive constant. The contribution
involving a proper prime power is $O(N^{1/2}(\log N)^3)=o(N)$: there are
$O(N^{1/2}\log N)$ such prime powers up to $2N$, and each weight is
at most $(\log(2N))^2$. Thus $T_{2,N}(n)>0$ for every nonexceptional
even target and sufficiently large $N$, giving a representation by primes
\[
 n=p+p',\qquad N\le p,p'\le2N.
\]
In fact the primes can be chosen strictly greater than $N$. The parameter
$m$ remains the prime cutoff: $2N\le m$ ensures that both primes belong to
$\A_m=\{p\le m:p\text{ odd prime}\}$ for large $N$. We use the recalled
almost-all estimate, not the main short-interval theorem of \cite{BL00}.

\begin{lemma}[Balanced lower-prime relations]\label{depth:moving-anchor}
There is an integer $q_0$ such that every prime $q>q_0$ admits odd primes
$p,p',p''<q$ such that
\[
 p+q=p'+p'',\qquad q/4<p<q/2,\qquad (q-1)/2<p',p''\le q-1.
\]
\end{lemma}
\begin{proof}
Put $N=(q-1)/2$. By the prime number theorem, the interval $(q/4,q/2)$
contains asymptotically $q/(4\log q)$ prime anchors $p$. The distinct
even targets $n=p+q$ lie in $(5N/2,7N/2]$ for large $q$. With a fixed
$C>1$, the exceptional set in the preceding estimate has
$O_C(q/(\log q)^C)=o(q/\log q)$ elements. At least one target
therefore has a genuine prime representation $p'+p''$ with
$N<p',p''\le2N=q-1$, as required.
\end{proof}

Return to the odd-prime supports $\A_m$.
Inclusion of supports induces maps
$U_m\to U_{m'}$, where
$U_m=U_{\A_m}$ and $m\le m'$.

\begin{theorem}\label{depth:integral-stabilization}
Let $U_\infty$ be the group on all odd primes modulo the
equal-pair-sum relations. For every sufficiently large $m$, the natural
map $U_m\to U_\infty$ is an isomorphism.
\end{theorem}
\begin{proof}
Choose $m_0>q_0$, so the prefix $\A_{m_0}$ contains all odd primes
at most $q_0$. For each
subsequent prime $q$, Lemma~\ref{depth:moving-anchor} gives
\[
 [q]=[p']+[p'']-[p]
\]
using only earlier labels. Hence every subsequent transition is
surjective, and all later groups are quotients of one fixed finitely
generated abelian group $U_{m_0}$. The kernels of the maps from $U_{m_0}$ form an
ascending chain of subgroups. It stabilizes because $U_{m_0}$ is a Noetherian
$\mathbb Z$-module. Stabilization of the kernels proves stabilization
of the groups and their transition maps, including their torsion.
The direct limit is $U_\infty$, since every defining relation
uses finitely many prime labels.
\end{proof}

\subsection{Depth and projective dimension}

For the binary odd-prime ring, use the cutoff notation of
Section~\ref{setup:resolution-construction}, with
$\A_m=\{p\le m:p\text{ odd prime}\}$ and $R_m=R^{<m}$.
Depth is taken at the homogeneous maximal ideal unless another
localization is specified. By \eqref{setup:AB-cutoffs}, depth measures
the difference between $|\A_m|$ and projective dimension. We bound
this difference for both truncations.

Both bounds use the same homological observation: localization preserves
a free resolution, although it need not preserve minimality. Thus the
projective dimension at the all-ones point is a lower bound for the
projective dimension over $R$. A dimension bound there, together with
Auslander--Buchsbaum, will give the desired bound.

\begin{lemma}[Triangular elimination]\label{depth:triangular}
Let $\kk$ be a field, let $\A$ be a finite ordered support with least label $a_0$, and put
$R=\kk[x_a:a\in\A]$. Suppose
$\mathfrak a\subseteq\mathfrak q=\langle x_a-1\mid a\in\A\rangle$.
Let $\mathcal B\subseteq\A\setminus\{a_0\}$, and suppose that for each
$b\in\mathcal B$ there are $a',a''\in\A$, depending on $b$, such that
\[
 x_{a_0}x_b-x_{a'}x_{a''}\in\mathfrak a,\qquad a',a''<b.
\]
Then $\pd_R(R/\mathfrak a)\ge|\mathcal B|$. If $\mathfrak a$ is homogeneous in monomial degree,
then $\operatorname{depth}(R/\mathfrak a)\le|\A|-|\mathcal B|$.
\end{lemma}
\begin{proof}
Put $M=R/\mathfrak a$, and let $\mathfrak b$ be generated by the displayed
relations. For recursive substitution, only $x_{a_0}$ needs to be inverted:
each equation gives $x_b=x_{a_0}^{-1}x_{a'}x_{a''}$. In increasing order
of $b\in\mathcal B$, we therefore solve for $x_b$ in terms of earlier
variables, obtaining
\[
 R[x_{a_0}^{-1}]/\mathfrak b\cong
 \kk[x_{a_0},x_{a_0}^{-1},x_c:c\in\A\setminus(\mathcal B\cup\{a_0\})],
\]
whose dimension is $|\A|-|\mathcal B|$.

We now localize further at $\mathfrak q$ to apply Auslander--Buchsbaum
over the regular local ring $R_{\mathfrak q}$, whose depth and dimension
are $|\A|$. Every $x_a$ becomes a unit there, since $x_a\notin\mathfrak q$,
but these additional inverses were not needed for substitution.
The inclusions $\mathfrak b\subseteq\mathfrak a\subseteq\mathfrak q$
ensure that $M_{\mathfrak q}\ne0$ and
$\dim M_{\mathfrak q}\le|\A|-|\mathcal B|$.
Localization preserves a free resolution, so Auslander--Buchsbaum gives
\[
 \pd_R M\ge\pd_{R_{\mathfrak q}}M_{\mathfrak q}
 =|\A|-\operatorname{depth}M_{\mathfrak q}\ge|\mathcal B|.
\]
If $\mathfrak a$ is homogeneous, the graded minimal resolution of $M$
stays minimal at the homogeneous maximal ideal
$\mathfrak m=\langle x_a\mid a\in\A\rangle$. Since $R_{\mathfrak m}$
also has depth $|\A|$, applying Auslander--Buchsbaum there gives the
depth comparison and bound explicitly:
\[
 \begin{aligned}
 \operatorname{depth}M
 &=|\A|-\pd_R M\\
 &\le|\A|-\pd_{R_{\mathfrak q}}M_{\mathfrak q}
   =\operatorname{depth}M_{\mathfrak q}\\
 &\le\dim M_{\mathfrak q}\le|\A|-|\mathcal B|.
 \end{aligned}\qedhere
\]
\end{proof}

For a Goldbach relation the label identity is $a_0+b=a'+a''$;
all three labels $a_0,a',a''$ are smaller than $b$.

\begin{proposition}\label{depth:relation-group-bound}
For every finite nonempty support $\A$ and every field $\kk$,
\[
 \operatorname{depth}G\le r_{\A},\qquad
 \pd_{R}G\ge|\A|-r_{\A}.
\]
\end{proposition}
\begin{proof}
Put $w=x_{\A}=\prod_{a\in\A}x_a$. Inverting $w$ makes each $x_a$
invertible. In the Laurent polynomial ring write
$x^\alpha=\prod_{a\in\A}x_a^{\alpha(a)}$ for $\alpha\in\mathbb Z^{\A}$,
allowing negative exponents. Each binary relation becomes
$x^{e_a+e_b-e_{a'}-e_{b'}}=1$. Since $x^{\alpha+\beta}=x^\alpha x^\beta$
and $x^{-\alpha}=(x^\alpha)^{-1}$, these relations identify exactly
the Laurent monomials whose exponents lie in the same coset of $L_{\A}$.
Thus
\[
 G[w^{-1}]\cong\kk[U_{\A}],
\]
with one basis element for each coset and multiplication given by
addition in $U_{\A}$. This is the Laurent form of the lattice
ideal construction \cite[Theorem~2.1]{ES96}, \cite[\S7.1]{MS05}.
By the structure theorem for finitely generated abelian groups,
$U_{\A}\cong\mathbb Z^{r_{\A}}\oplus T_{\A}$, with $T_{\A}$ finite.
Its group algebra is a Laurent polynomial algebra over the nonzero
finite-dimensional algebra $\kk[T_{\A}]$, and hence has dimension $r_{\A}$.
This remains true when the characteristic divides $|T_{\A}|$.

All defining binomials vanish at the all-ones point, so the quotient
localized at $\mathfrak q=\langle x_a-1\mid a\in\A\rangle$ is nonzero and has dimension
at most $r_{\A}$. Repeating the final argument of
Lemma~\ref{depth:triangular}, local
Auslander--Buchsbaum gives $\pd_{R}G\ge|\A|-r_{\A}$, and graded
Auslander--Buchsbaum gives the asserted depth bound.
\end{proof}

\begin{corollary}[Depth and total Betti numbers]
\label{depth:bounded-depth}
For all sufficiently large $m$, over every field,
\[
 \delta^{<m}\le2,\qquad
 |\A_m|-2\le\rho^{<m}\le |\A_m|.
\]
Moreover the total Betti numbers satisfy
$\beta^{<m}_{j,*,*}\ge\binom{|\A_m|-2}{j}$ for $0\le j\le |\A_m|-2$.
\end{corollary}
\begin{proof}
By Theorems~\ref{depth:integral-stabilization} and
\ref{depth:prime-pair-group}, $U_m\cong\mathbb Z^2$ for all
sufficiently large $m$. Hence $L_{\A_m}$ is a direct summand of
$\mathbb Z^{|\A_m|}$ of rank $|\A_m|-2$. After inverting all variables,
extend a basis of $L_{\A_m}$ to a basis of the exponent lattice.
In the resulting Laurent coordinates the ideal is generated by
$y_1-1,\ldots,y_{|\A_m|-2}-1$, a regular sequence. At the point with
all coordinates one its Koszul resolution is minimal. The graded
minimal resolution over $R_m$ localizes to a free resolution there,
and tensoring with the residue field gives a complex whose $j$th term
has dimension $\beta^{<m}_{j,*,*}$ and whose $j$th homology is the
local Tor group of dimension $\binom{|\A_m|-2}{j}$. Since homology
cannot have larger dimension than its chain group, this gives the
termwise Betti bounds. The
projective dimension bounds now follow from the last nonzero local
rank and Hilbert's syzygy theorem. Graded Auslander--Buchsbaum gives
the depth bound.
\end{proof}

In particular $\rho^{<m}=|\A_m|-O(1)$, with
$|\A_m|=\pi(m)-1$ for integers $m\ge3$. The auxiliary primes in
Lemma~\ref{depth:prime-pair-conversion} may exceed the given cutoff,
so the isomorphism in Theorem~\ref{depth:integral-stabilization} is
essential to this conclusion. Its proof leaves the numerical onset of
stabilization open. Corollary~\ref{depth:verified-range} below verifies
the same rank-two group and depth bound throughout an explicit finite
range, without identifying the onset in the eventual theorem.

\subsection{The relation-weight cutoff}

The relation $p+q=p'+p''$ in Lemma~\ref{depth:moving-anchor} can have weight
larger than $m$ when $q$ is near the cutoff. It therefore does not give
the same uniform depth bound for $G^{\le m}$. A fixed anchor yields a
weaker but still asymptotically sharp projective-dimension estimate.

For an integer $m\ge6$, set
\[
 \begin{aligned}
 E_m&=\{p\in\A_m:p+3\ne q+r\text{ for all primes }q,r>3\},
       &e_m&=|E_m|,\\
 u_m&=\#\{p\in\A_m:p>m-3\}.&&
 \end{aligned}
\]
Repeated summands are allowed, and $u_m\le2$ since $(m-3,m]$ contains at most two odd integers.

\begin{proposition}\label{depth:fixed-anchor}
Over every field,
\[
 \delta^{<m}\le e_m,
 \qquad\delta^{\le m}\le e_m+u_m.
\]
\end{proposition}
\begin{proof}
For $p\notin E_m$, choose $p+3=q+r$ with $q,r>3$ prime. Then
$q,r\le p-2<p$, so $x_3x_p-x_qx_r$ is a triangular relation in the
variable cutoff. It is present in the relation cutoff when $p\le m-3$.
All binomials vanish at the all-ones point. Apply
Lemma~\ref{depth:triangular}, omitting at most $u_m$ relations in the
second case.
\end{proof}

\begin{corollary}[A verified range of variable cutoffs]
\label{depth:verified-range}
Put $B=20\,000\,000$. Over every field, for every integer $3\le m\le B$,
\[
 \depth G^{<m}=\delta^{<m}\le2.
\]
For $5\le m\le B$, moreover, $U_m\cong\mathbb Z^2$ and
\[
 \beta^{<m}_{i,*,*}\ge\binom{|\A_m|-2}{i}
 \qquad(0\le i\le|\A_m|-2).
\]
\end{corollary}
\begin{proof}
The exhaustive sieve in
\path{supplement/computations/prime-pair-depth-check.py} verifies that
each of the $1\,270\,604$ primes $7\le p\le B$ has
$p+3=q+r$ with primes $q,r>3$. Thus $E_m=\{3,5\}$ for
$6\le m\le B$, and Proposition~\ref{depth:fixed-anchor} gives the
depth bound. The cases $m=3,4,5$ are polynomial rings in at most two
variables.

Every verified pair has $q,r\le p-2$, so it lies in the same prime
prefix as $p$. Induction on $p$ in the relation group, using
$[p]=[q]+[r]-[3]$, gives
\[
 [p]=(1-j)[3]+j[5],\qquad j=(p-3)/2.
\]
Hence $[3],[5]$ generate $U_m$. The map
$[p]\mapsto(1-j,j)\in\mathbb Z^2$ respects every equal-pair-sum
relation and sends these two generators to the standard basis.
They are therefore independent, proving $U_m\cong\mathbb Z^2$.
The localization argument in the proof of
Corollary~\ref{depth:bounded-depth} now gives the total Betti bounds
over every field.
\end{proof}

This computation does not determine whether its verified range
overlaps the range of the eventual theorem, and it does not prove
the depth bound at every cutoff. The fixed-anchor condition is
sufficient for the bound; no necessity is asserted. For the
relation-weight cutoff the same verification gives only
$\delta^{\le m}\le2+u_m\le4$ for $6\le m\le B$, since a relation
with anchor $3$ and label $p$ has weight $p+3$.

\begin{theorem}\label{depth:relation-cutoff-pd}
For every fixed $C>0$, over every field,
\[
 0\le\delta^{\le m}
   =|\A_m|-\rho^{\le m}
   \ll_C m(\log m)^{-C}.
\]
Consequently $\rho^{\le m}/|\A_m|\to1$.
\end{theorem}
\begin{proof}
Let $E^*(X)$ count the even integers at most $X$ having no representation
as a sum of two primes greater than three. The balanced almost-all
estimate above implies
\[
 E^*(X)\ll_C X(\log X)^{-C}\qquad(C>0).
\]
To see this with the required restriction on summands, cover
$[\sqrt X,X]$ by the intervals $(5N/2,7N/2]$ at geometric scales
with ratio $4/3<7/5$. Their lengths sum to $O(X)$ and their logarithms
are comparable to $\log X$. At each scale, nonexceptional targets
have both prime summands greater than three. Integer rounding does not
destroy the overlap for large scales. The remaining initial interval
contributes $O(\sqrt X)$, which is absorbed in the stated bound.

The injective map $p\mapsto p+3$ sends $E_m$ into the exceptional set
counted by $E^*(m+3)$. Thus $e_m\ll_C m(\log m)^{-C}$. Apply
Proposition~\ref{depth:fixed-anchor} and $u_m\le2$. Finally,
$|\A_m|\sim m/\log m$, so taking $C>1$ proves the limiting ratio.
\end{proof}

The theorem does not show bounded depth for the relation cutoff or
characteristic-independent exact depth at each finite cutoff. For the
variable cutoff, Corollary~\ref{depth:bounded-depth} is stronger than
this limiting ratio.

The group computation also describes the comparison with the toric
quotient after localization. On the full odd-prime support, two
monomials of the same bidegree become equal after multiplication by
some common monomial. Appendix~\ref{app:saturation} proves the
squarefree-multiplier refinement, studies saturation by a single
prime variable, and gives examples where cancellation fails. It also
develops explicit dual Koszul cocycles and the stronger infinite
conjecture associated with them.

\section{Linear-strand growth for the two truncations}
\label{sec:growth}

Section~\ref{sec:betti-euler-shifts} studies $n\to\infty$ with monomial
degree $d$ fixed. Here we fix the homological degree $i$ and relation
arity $k$, sum over all Goldbach degrees $n$, and let the cutoff
$m\to\infty$.
For the binary ring the first possible syzygies in this homological degree
have monomial degree $i+1$. In the Koszul complex their coefficients
have monomial degree one. Their total number has a different leading constant for
each of our two truncations. Both constants can be computed by comparing
three consecutive Koszul terms over the given field $\kk$.

\subsection{An elementary estimate in the binary case}

Let $G=R/I$ be the binary Goldbach ring on a finite nonempty support
$\A$. In homological degree $i\ge1$, the linear strand contributes
$\beta_{i,*,i+1}$ to the total $\beta_{i,*,*}$.

\begin{proposition}\label{grow:binary-bound}
For every $i\ge1$,
\[
 \begin{aligned}
 \beta_{i,*,i+1}
 &\le |\A|\binom{|\A|}{i}-\binom{|\A|}{i+1},\\
 \beta_{i,*,i+1}
 &\ge |\A|\binom{|\A|}{i}-\binom{|\A|}{i+1}
       -|\A+\A|\binom{|\A|}{i-1}.
 \end{aligned}
\]
Consequently, along any family with $|\A|\to\infty$ and $|\A+\A|=o(|\A|^2)$,
\[
 \beta_{i,*,i+1}\sim\frac{i}{(i+1)!}|\A|^{i+1}
 \qquad\text{for every fixed }i\ge1.
\]
\end{proposition}
\begin{proof}
Write $G_{*,r}=\bigoplus_n G_{n,r}$, where $*$ denotes the direct sum
over all Goldbach degrees. In monomial degree $i+1$, the relevant
Koszul terms are
\[
 G_{*,0}\otimes\bigwedge^{i+1}\kk^{\A}
 \xrightarrow{\mathrm{d}_{i+1}}
 G_{*,1}\otimes\bigwedge^i\kk^{\A}
 \xrightarrow{\mathrm{d}_i}
 G_{*,2}\otimes\bigwedge^{i-1}\kk^{\A}.
\]
Their dimensions are $\binom{|\A|}{i+1}$, $|\A|\binom{|\A|}{i}$, and
$|\A+\A|\binom{|\A|}{i-1}$. Indeed, there are no linear relations, and
all quadratic monomials of a given represented weight agree in $G$.
The columns of $\mathrm{d}_{i+1}$ have disjoint supports and unit entries:
a term $x_a\varepsilon_{\sigma}$, with $a\notin \sigma$, determines its unique
source $\varepsilon_{\sigma\cup\{a\}}$. This map is therefore injective
over every field. Summing \eqref{eq:koszul-rank-nullity} over Goldbach
degrees $n$ therefore gives
\[
 \beta_{i,*,i+1}
 =|\A|\binom{|\A|}{i}-\binom{|\A|}{i+1}
   -\dim_\kk\operatorname{im}\mathrm{d}_i.
\]
Bounding the rank by zero and by the target dimension proves the
inequalities. For fixed $i$, the main term is asymptotic to
$i|\A|^{i+1}/(i+1)!$, while the error is $O_i(|\A+\A|\,|\A|^{i-1})$.
\end{proof}

Let $m$ denote the prime cutoff as in Section~\ref{setup:resolution-construction}.
For $\A_m=\{p\le m:p\text{ odd prime}\}$, the prime number theorem gives
$|\A_m|\sim m/\log m$.
Since $|\A_m+\A_m|=O(m)=o(|\A_m|^2)$, the proposition gives the binary
variable-cutoff asymptotic. More generally, the same proof applies to
any standard graded quotient without linear relations, with $|\A+\A|$
replaced by its degree-two Hilbert function.

\subsection{The linear strand at the two cutoffs}

Already among the cubic second syzygies the distinction is visible.
Summed over all weights, their number is asymptotic to $|\A_m|^3/3$
for $G^{<m}$ and to $|\A_m|^3/8$ for $G^{\le m}$. The same comparison
gives the constants in every fixed homological degree.

\begin{theorem}\label{grow:both-cutoffs}
For the binary ring over any field and every fixed $i\ge1$, as $m\to\infty$,
\[
 \beta^{<m}_{i,*,i+1}\sim\frac{i}{(i+1)!}|\A_m|^{i+1},
 \qquad
 \beta^{\le m}_{i,*,i+1}
       \sim\frac{1-2^{-i}}{(i+1)!}|\A_m|^{i+1}.
\]
More generally, fix $k\ge2$ and $i\ge1$, and let $G$ be the $k$-ary Goldbach
ring on the odd primes. Its variable cutoff $G^{<m}$ uses all
defining relations on $\A_m$, and its relation cutoff $G^{\le m}$
uses those of common weight at most $m$. Then
\begin{align}
 \beta^{<m}_{i,*,i+k-1}
 &\sim\frac{|\A_m|^{k+i-1}}{(k-1)!(i-1)!(k+i-1)},
 \label{grow:variable-constant}\\
 \beta^{\le m}_{i,*,i+k-1}
 &\sim C_{k,i}|\A_m|^{k+i-1},\qquad
 C_{k,i}=\frac{\displaystyle\sum_{j=1}^k
       (-1)^{j+1}\binom kj j^{1-i}}{k!(k+i-1)!}.
 \label{grow:relation-constant}
\end{align}
\end{theorem}

The next lemma uses only the monomial grading. Within its statement
$\beta_{i,d}(M)$ denotes $\dim_\kk\Tor_i^R(M,\kk)_d$.
For a bigraded quotient this is $\beta_{i,*,d}(M)$.
In its proof, $K_j(M)_d$ denotes the monomial-degree-$d$ part of the
$j$th Koszul term.

\begin{lemma}\label{grow:koszul-comparison}
Suppose $\mathfrak a\subseteq\mathfrak b\subseteq R=\kk[x_1,\ldots,x_N]$ are homogeneous
ideals with no elements in degrees below $k\ge2$. If
$s=\dim_\kk(\mathfrak b/\mathfrak a)_k$, then for $i\ge1$,
\[
 0\le\beta_{i,i+k-1}(R/\mathfrak b)-\beta_{i,i+k-1}(R/\mathfrak a)
       \le s\binom N{i-1}.
\]
\end{lemma}
\begin{proof}
Put $E=\kk^N$, $d=i+k-1$, and $K(M)=\mathbb K_{\kk}\otimes_R M$.
Since $K_j(M)_d=M_{d-j}\otimes_\kk\bigwedge^jE$, the coefficient
degrees in homological degrees $i+1,i,i-1$ are $k-2,k-1,k$.
The two quotients agree below degree $k$, so they have the same
incoming map $A:U\to V$, where
\[
 U=R_{k-2}\otimes_\kk\bigwedge^{i+1}E,\qquad
 V=R_{k-1}\otimes_\kk\bigwedge^iE.
\]
The common term $\dim_\kk(V/\operatorname{im}A)$ in
\eqref{eq:koszul-rank-nullity} cancels when the two Betti numbers are
subtracted. Only the outgoing targets differ:
\[
 W_{\mathfrak a}=(R/\mathfrak a)_k\otimes_\kk\bigwedge^{i-1}E,
 \qquad
 W_{\mathfrak b}=(R/\mathfrak b)_k\otimes_\kk\bigwedge^{i-1}E.
\]
The quotient $\rho:W_{\mathfrak a}\twoheadrightarrow W_{\mathfrak b}$
satisfies $B_{\mathfrak b}=\rho B_{\mathfrak a}$ for the outgoing maps,
and $\dim_\kk\ker\rho=s\binom N{i-1}$. Thus
$\ker B_{\mathfrak a}\subseteq\ker B_{\mathfrak b}$, while
$\operatorname{im}A$ is unchanged. The map $v\mapsto B_{\mathfrak a}v$
identifies
\[
 \ker B_{\mathfrak b}/\ker B_{\mathfrak a}
 \cong\operatorname{im}B_{\mathfrak a}\cap\ker\rho.
\]
Consequently
\[
 0\le\beta_{i,d}(R/\mathfrak b)-\beta_{i,d}(R/\mathfrak a)
 =\dim_\kk(\operatorname{im}B_{\mathfrak a}\cap\ker\rho)
 \le s\binom N{i-1}.
\]
\end{proof}

\begin{proof}[Proof of Theorem~\ref{grow:both-cutoffs}]
Within this proof, write $I$ for the binomial ideal of the chosen
truncation and $I'$ for its monomial comparison ideal. For the relation
cutoff, take $I=I^{\le m}$ and
\[
 I'=\bigl(x_{p_1}\cdots x_{p_k}\mid p_1,\ldots,p_k\in\A_m,
       \ \textstyle\sum_{j=1}^k p_j\le m\bigr).
\]
Both ideals have no elements below monomial degree $k$, and
$I\subseteq I'$ because each generator of $I$ is a difference of two
listed monomials.

In monomial degree $k$, the binomial ideal $I$ is the coefficient-sum-zero
subspace of each represented fiber of weight at most $m$, whereas $I'$
contains the whole fiber. Therefore $s=\dim_\kk(I'/I)_{*,k}$ counts the
represented weights and is $O_k(m)$. For the variable cutoff, take instead
$I=I^{<m}$ and $I'=\m^k$, where $\m=(x_p\mid p\in\A_m)$;
then $s=|k\A_m|=O_k(m)$. In either case the error from
Lemma~\ref{grow:koszul-comparison} is
\[
 O_{k,i}(m |\A_m|^{i-1})=o(|\A_m|^{k+i-1}).
\]
Indeed, after division by $|\A_m|^{k+i-1}$, the error is
$O_{k,i}(m/|\A_m|^k)$, and the prime number theorem gives
\[
 \frac{m}{|\A_m|^k}
 \sim\frac{(\log m)^k}{m^{k-1}}\longrightarrow0.
\]

List the variables in increasing order of prime label. Both comparison
ideals are strongly stable: replacing a variable by one with smaller
index does not increase the additive weight. They are generated in
degree $k$. If $q(x^\alpha)$ is the largest index of a variable dividing
a generator $x^\alpha$, the Eliahou--Kervaire formula \cite{EK90} gives
\begin{equation}\label{grow:EK}
 \beta_{i,*,i+k-1}(R_m/I')
 =\sum_{x^\alpha\in\mathcal G(I')}\binom{q(x^\alpha)-1}{i-1}.
\end{equation}
All positive Betti numbers of the comparison quotient $R_m/I'$ lie on
this strand. The formula is valid
over every field. In the case $I'=\m^k$, the sum becomes
\[
 \sum_{j=1}^{|\A_m|}\binom{j+k-2}{k-1}\binom{j-1}{i-1}
 \sim\frac{|\A_m|^{k+i-1}}{(k-1)!(i-1)!(k+i-1)}.
\]

For the relation cutoff, the prime number theorem gives
$|\A_{tm}|/|\A_m|\to t$ for each fixed $t\in(0,1]$.
Monotonicity and continuity of the limiting distribution function make
this convergence uniform on $[0,1]$.
Thus the measures $|\A_m|^{-1}\sum_{p\in\A_m}\delta_{p/m}$ converge weakly
to uniform measure, and $q(x^\alpha)/|\A_m|$ differs uniformly by $o(1)$ from
the largest scaled label in $x^\alpha$. Generators with a repeated label
contribute only $O_{k,i}(|\A_m|^{k+i-2})$ to \eqref{grow:EK}. Counting the
remaining ordered tuples and dividing by $k!$ therefore gives
\begin{equation}\label{grow:simplex-moment}
 \frac{\sum_n\beta_{i,n,i+k-1}(R_m/I')}{|\A_m|^{k+i-1}}
 \longrightarrow\frac{1}{k!(i-1)!}
 \int_{\substack{u_j\ge0\\\sum_ju_j\le1}}
       (\max_j u_j)^{i-1}\,du.
\end{equation}
The simplex boundary has measure zero, so the cutoff indicator presents
no obstruction to this limit.

Put $r=i-1$. For $r\ge1$, inclusion--exclusion yields
\[
 \operatorname{vol}\{u\ge0:\textstyle\sum_ju_j\le1,\ \max_j u_j>s\}
 =\frac{1}{k!}\sum_{j=1}^k(-1)^{j+1}\binom kj(1-js)_+^k.
\]
Integrating against $r s^{r-1}\,ds$ and evaluating the beta integral gives
\[
 \int_{\sum_ju_j\le1}(\max_j u_j)^r\,du
 =\frac{r!}{(k+r)!}\sum_{j=1}^k(-1)^{j+1}\binom kj j^{-r}.
\]
For $r=0$ the same formula holds, because the simplex has volume $1/k!$
and the alternating binomial sum is one. Substitution into
\eqref{grow:simplex-moment}, followed by the rank comparison, proves
\eqref{grow:relation-constant}. At $k=2$, the sum simplifies to the
stated binary constant.
\end{proof}

The ratio of the binary
relation-cutoff linear strand to the variable-cutoff linear strand
tends to $(1-2^{-i})/i$. These leading constants are independent of
the characteristic, although individual Betti numbers need not be.

\subsection{The first nonlinear strand}

The first nonlinear contribution to the binary variable cutoff is
smaller. For $\rho\ge0$, write uniquely
$\rho=\binom a2+b$ with $a\ge1$ and $0\le b<a$, and set
\[
 \rho^{\langle2\rangle}=\binom{a+1}{3}+\binom{b+1}{2}.
\]
This is the degree-two Macaulay successor, and
$\rho^{\langle2\rangle}=O(\rho^{3/2})$ for $\rho\ge1$.

\begin{proposition}\label{grow:nonlinear-bound}
Let $G=R/I$ be the binary ring on a finite nonempty support $\A$, and
put $\rho=|\A+\A|$.
Over every field and for every $i\ge1$,
\[
 \beta_{i,*,i+2}
 \le\binom{|\A|-1}{i-1}
       \bigl(\rho^{\langle2\rangle}-\dim_\kk G_{*,3}\bigr)
 \le\binom{|\A|-1}{i-1}\rho^{\langle2\rangle}.
\]
\end{proposition}
\begin{proof}
Write $H(d)=\dim_\kk G_{*,d}$. Order the support as
$a_1<\cdots<a_{|\A|}$, and let $L\subseteq R$ be the lex-segment ideal
whose quotient $R/L$ has Hilbert function $H$, using
$x_{a_1}>\cdots>x_{a_{|\A|}}$.
Macaulay's theorem supplies this ideal; see the discussion preceding
\cite[Lemmas~2--3]{JKV21}. It has no linear generators since $H(1)=|\A|$.

Let $Q$ be the ideal generated by the quadratic monomials of $L$.
Since $H(2)=\rho$ and the quadratic lex segment attains Macaulay's
growth bound,
\[
 \dim_\kk(R/Q)_{*,3}=\rho^{\langle2\rangle}.
\]
The cubic monomials of $L$ not divisible by its quadratic generators
form a basis of $(L/Q)_{*,3}$ and are precisely its minimal cubic
generators, denoted $\mathcal G(L)_3$. Hence
\[
 |\mathcal G(L)_3|=\dim_\kk(L/Q)_{*,3}
 =\rho^{\langle2\rangle}-H(3).
\]
These are generators of the comparison ideal $L$, even though $I$
itself is generated by quadrics.

For a monomial $x^\alpha$, put
$q(x^\alpha)=\max\{j:x_{a_j}\text{ divides }x^\alpha\}$.
The Bigatti--Hulett--Pardue theorem and the Eliahou--Kervaire formula,
valid over arbitrary fields, give
\[
 \beta_{i,*,i+2}\le\beta_{i,*,i+2}(R/L)
 =\sum_{x^\alpha\in\mathcal G(L)_3}\binom{q(x^\alpha)-1}{i-1}
 \le\binom{|\A|-1}{i-1}\bigl(\rho^{\langle2\rangle}-H(3)\bigr);
\]
see \cite[Lemmas~2--3]{JKV21} and \cite{EK90}.
The remaining bound follows from $H(3)\ge0$. Finally,
$a=O(\sqrt\rho)$ and $b<a$ in the Macaulay representation give
$\rho^{\langle2\rangle}=O(\rho^{3/2})$.
\end{proof}

For the odd primes at most $m$, the support has size $|\A_m|$ and
$\rho=|\A_m+\A_m|=O(m)$. As $m^{3/2}/|\A_m|^2\sim(\log m)^2/\sqrt m\to0$,
the proposition gives
\[
 \beta^{<m}_{i,*,i+2}
 =O_i(|\A_m|^{i-1}m^{3/2})=o(|\A_m|^{i+1}).
\]
Thus for the variable cutoff the sum of the first two strands has
the same leading asymptotic as the linear strand alone. Determining
the asymptotics of the total Betti number $\beta^{<m}_{i,*,*}$ requires
control of monomial degrees $i+3,i+4,\ldots$. All estimates here are for fixed $i$.

\section{Further questions}
\label{sec:questions}

The first question asks whether binary moves connect every represented
fiber on the full odd-prime support.

\begin{question}\label{question:full-support-connectivity}
For the full binary Goldbach ring $G$ on the odd primes, is every
nonempty fiber $\mathcal M_{n,d}$ connected by binary moves?
Equivalently, is the natural map
\[
 G\longrightarrow\kk[u^pv:p\text{ an odd prime}]
\]
an isomorphism, and is $G$ a domain?
\end{question}
The equivalence with the toric presentation follows from the component
basis in Proposition~\ref{setup:components}. For the full support,
Theorem~\ref{depth:prime-pair-group} identifies any two monomials of
the same bidegree after multiplication by a common monomial.
Such a multiplier is nonzero in $G$. Thus failure of connectivity
would give a nonzero zero divisor, while the toric quotient is a
domain. The localization argument in Appendix~\ref{app:saturation}
does not by itself permit cancellation before localization.

For each fixed monomial degree $d\ge3$, Theorem~\ref{main:weightwise} gives a
positive lower bound for the fraction of the linear strand spanned by
cycles from positive common shifts. Already in monomial degree three there is
a finer question.

\begin{question}
For the full binary ring on the odd primes, let $W^+_{n,3}$ be the span
of the positive common-shift classes at weight $n$. For odd $n$ tending
to infinity, does
\[
 \frac{\dim_\kk W^+_{n,3}}{\beta_{2,n,3}}\longrightarrow1?
\]
How does the answer depend on the coefficient field?
\end{question}

\begin{question}
For fixed $i$, what is the growth of
$\beta^{<m}_{i,*,*}$ and $\beta^{\le m}_{i,*,*}$
when all monomial degrees are included?
\end{question}
Theorem~\ref{grow:both-cutoffs} determines the linear contribution.
For the binary variable cutoff the next strand has lower order.
Estimating the total requires control of syzygies in all later
monomial degrees. One may also let $i$ vary with $m$, a regime
outside the asymptotic formula for fixed homological degree.

\begin{question}
Does $\depth G^{<m}\le2$ hold at every variable cutoff? Which
cutoffs have depth zero, one, or two? Is the depth of the binary
relation-weight truncations uniformly bounded?
\end{question}
Theorem~\ref{depth:prime-pair-group} identifies the stable integral
group with $\mathbb Z^2$. Corollary~\ref{depth:bounded-depth} gives
the bound two for all sufficiently large variable cutoffs, while
Corollary~\ref{depth:verified-range} verifies it for
$3\le m\le20\,000\,000$. We do not know that these two ranges
overlap. The exact depth at each cutoff and the universal bound
remain open. The prime-pair condition at $3$ is sufficient for the
bound, but the depth question does not assert that this condition
is necessary.
For $G^{\le m}$ the bound in Theorem~\ref{depth:relation-cutoff-pd}
uses only pair relations whose weight is at most $m$.
Conjecture~\ref{depth:dual-conjecture} proposes a nonzero infinite
dual Koszul class. Its finite representatives give explicit depth
bounds at the checked cutoffs; Proposition~\ref{depth:dual-prime-pair}
relates the conjecture to the prime-pair condition for $x_3$-saturation.

\subsection*{Acknowledgments}
The author thanks his Clemson University Ph.D. advisor,
Keri Sather-Wagstaff, for frequent conversations and meetings and for
encouraging him to pursue this research.

\subsection*{Computational data}
The arithmetic supplement accompanying this paper contains the
illustrated resolution strands, the weight-$33$ Koszul cycles, the
linear-strand constants, and the saturation, variable-quotient, and
dual-cocycle calculations in Appendix~\ref{app:saturation}. It also
contains the exhaustive prime-pair check for
Corollary~\ref{depth:verified-range}. The supplement is supplied in
\path{anc/supplement} with the arXiv sources; paths in the text are
relative to its parent directory.
Its \path{README.txt} identifies the assertions checked, dependencies,
and expected output. The \path{verify.py} runner verifies the manifest
hashes and executes six standard-library Python checks, including the
prime-pair sieve through $20\,000\,000$. The ideal, exact-depth, and
Tor calculations in Appendix~\ref{app:saturation} require running the
three supplied \textsc{Singular} scripts separately; the Python runner
does not execute those scripts. The finite depth corollary uses the
prime-pair check and the localization argument in its proof, not a
computed minimal resolution. The arithmetic asymptotic and
interpolation theorems are proved in the text.

\Needspace{17\baselineskip}
\subsection*{Use of generative AI}
The Goldbach-ring framework, the construction and existence proof of its
minimal bigraded free resolution, the study of cellular models for its
fixed-bidegree complexes, and the shifted-Koszul-cycle proof of
infinite projective dimension originated with the author. The initial
computational identification of characteristic-two dependence, treated
in the companion paper, also originated with the author.

During the subsequent development of this work, the author used OpenAI's
GPT-6 Astra for substantive mathematical assistance. The model proposed
the results stated as Theorems~\ref{main:weightwise},
\ref{main:cutoff-growth}, and~\ref{main:prime-interpolation} and supplied
crucial additive-number-theoretic arguments and initial proof drafts,
building on the author's constructions. The author directed this
development, checked and verified the arguments, and revised their
presentation. AI assistance was also used to develop computational code
and verification routines and to organize and edit the manuscript.

The author takes responsibility for the mathematical claims, proofs,
computations, references, and final presentation.

\appendix
\section{The scroll bound for a fixed shift}
\label{app:scroll-shift}

We refine the fixed-shift bound of Theorem~\ref{thm:translation} using
the full set of linear syzygies of a rational normal scroll.
Castryck--Cools--Demeyer--Lemmens
use the linear strand of an ambient scroll for a toric surface
\cite[Sec.~4.1]{CCDL19}. Here the Goldbach grading records the weight
of each scroll generator.

\begin{proposition}[The scroll bound for a fixed shift]\label{thm:scroll-shift}
Fix $h>0$, let $B\subseteq\A\cap(\A-h)$ be finite, and put $r=|B|$.
For the binary Goldbach ring on $\A$, over every field and for $i\ge1$,
\[
 \beta_{i,n,i+1}\ \ge\
 \sum_{j=1}^{i}\#\left\{\sigma\subseteq B:\
       |\sigma|=i+1,\quad \sum_{b\in\sigma} b+jh=n\right\}.
\]
In particular $\beta_{i,*,i+1}\ge i\binom{r}{i+1}$.
\end{proposition}
\begin{proof}
We may assume $r\ge i+1$. In the polynomial ring on the labels of
$B\cup(B+h)$, form the matrix $M$ with columns $(x_b,x_{b+h})^{\mathsf T}$
for $b\in B$. Its minors belong to $I$, since $b+(b'+h)=(b+h)+b'$.
The edges $b\to b+h$ form disjoint paths. Each vertex has at most one
predecessor and one successor, and labels increase along each edge.
Thus the columns of a path with vertices $a_0,\ldots,a_s$, in order,
form a block
\[
 \begin{pmatrix}x_{a_0}&\cdots&x_{a_{s-1}}\\
                 x_{a_1}&\cdots&x_{a_s}\end{pmatrix}.
\]
Different blocks use disjoint variables. This is a rational normal
scroll matrix. If there are $t$ paths, there are $r+t$ variables. Over
an algebraic closure, on the locus of the minor equations where the top row
is nonzero, the bottom row
equals a scalar multiple of it. That scalar and the initial variable
of each path parametrize this locus, of dimension $t+1$. Where the
top row vanishes only the $t$ terminal variables are free. Hence the
minor ideal has height $r-1$, and the Eagon--Northcott complex is a
minimal resolution of its quotient \cite{EN62}. Extending it to $R$ by the remaining
variables preserves exactness and minimality, also for infinite $\A$.

To read its bidegrees, assign weight $b$ and monomial degree one to the column
indexed by $b$, and give the two row basis elements weights $0,-h$
and monomial degree zero. This makes $M$ homogeneous. In homological degree $i$ an
Eagon--Northcott generator uses $i+1$ columns indexed by $\sigma$,
the dual row determinant, and a divided power of degree $i-1$ in
the dual rows. The determinant contributes weight $h$. Taking $j-1$
copies of the second dual row contributes another $(j-1)h$, for
$j=1,\ldots,i$. Thus the generators have precisely the weights and
monomial degrees in the formula.

The quotient map $R/I_2(M)\to G$ induces an injection on the linear
strand. Indeed, in monomial degree $i+1$, the Koszul terms in homological degrees
$i+1$ and $i$ have coefficient monomial degrees zero and one. The quotient map
identifies both terms and their incoming differential. It takes each
cycle to a cycle, and both complexes have the same boundaries in
homological degree $i$. This proves the weightwise bound. Summing over $n$
counts each $(i+1)$-subset $i$ times.
\end{proof}

The classes $z_{\sigma,h}$ give the $j=1$ contribution. To see this,
first use independent entries $X_b,Y_b$ in the two rows of $M$.
Give both entries their column degree and count each lower-row entry
once in a further grading. The chain \eqref{eq:translated-cycle}, with
coefficient $Y_b$ and exterior generator indexed by $X_b$, is a cycle
modulo the minors. Its column multidegree is $\sigma$ and its lower-row
degree is one. The scroll resolution has just one generator in this
fine degree. Substituting $x_b,x_{b+h}$ for $X_b,Y_b$ gives
$z_{\sigma,h}$, whose image in $G$ is nonzero by
Theorem~\ref{thm:translation}.

For example, $B=\{3,11,17\}$ and $h=2$ give two cubic classes, at
weights $33$ and $35$. The first is the introductory cycle. The second
can be written
\[
 x_3\varepsilon_{13}\varepsilon_{19}
 -x_{11}\varepsilon_5\varepsilon_{19}
 +x_{17}\varepsilon_5\varepsilon_{13}.
\]
It comes from interchanging the rows of $M$. The proposition counts
these classes for one shift at a time. For instance the shifts
$h=4$ on $\{3,13\}$ and $h=10$ on $\{3,7\}$ give cycles whose
difference is $\mathrm{d}(\varepsilon_7\varepsilon_{13})$.
Thus the fixed-shift bounds need not add across different shifts.

\section{Saturation, variable quotients, and dual Koszul cocycles}
\label{app:saturation}

The relation-group theorem in Section~\ref{sec:depth} identifies
monomials after all prime variables are inverted. This appendix
examines what remains before localization: common multipliers,
quotients by one variable, and explicit dual Koszul classes.
Write $G=R/I$ for the binary ring and $S=R/J$ for its toric quotient.
At a finite variable cutoff put $\A_m=\{p\le m:p\text{ odd prime}\}$
and $w_m=\prod_{p\in\A_m}x_p$.

\subsection{Localization and saturation}

We ask whether the binary relations generate every
relation between monomials of the same bidegree. The toric map
\[
 G=R/I\twoheadrightarrow S=R/J
       =\kk[u^pv:p\text{ odd prime}],\qquad x_p\longmapsto u^pv,
\]
identifies all such monomials. By Proposition~\ref{setup:components},
$I=J$ precisely when every nonempty fiber is connected by binary
Goldbach moves, or equivalently when $\dim_\kk G_{n,d}\le1$ for every
$(n,d)$. The finite toric quotient $S^{<m}$ has dimension two for
$m>5$, so its depth is at most two. The full $S$ also has dimension
two; the associated projective curve has dimension one, as explained
in Section~\ref{sec:setup}. The surjection $G^{<m}\to S^{<m}$ alone,
however, does not give a depth bound for $G^{<m}$.

For a general finite support $\A$, put $w=\prod_{a\in\A}x_a$.
The toric map $G\twoheadrightarrow S$ from Section~\ref{sec:setup}
records only Goldbach degree and monomial degree. The group $U_{\A}$ retains the
binary relations and maps onto the subgroup of $\mathbb Z^2$ generated
by the bidegrees $(a,1)$, by $[a]\mapsto(a,1)$. In polynomial
coordinates the kernel of $G\to\kk[U_{\A}]$ is
$(I:w^\infty)/I$. In particular, two monomials $x^\alpha,x^\beta$
have the same image in this group algebra precisely when
$w^j x^\alpha=w^j x^\beta$ in $G$ for some $j\ge0$.
Here $I:w^\infty=\bigcup_{j\ge0}(I:w^j)$. In a finite primary
decomposition this saturation removes exactly the components whose
radicals contain at least one variable. Such a component may be
minimal or embedded. Since $S$ is a domain and $\phi(w)\ne0$,
\[
 I\subseteq (I:w)\subseteq (I:w^\infty)\subseteq J.
\]

Appendix~\ref{depth:prime-variable} computes these ideals at cutoff $42$
and explains how the colon enters the resolution of $G/x_p$.

For the full odd-prime support, Theorem~\ref{depth:prime-pair-group}
identifies $U_\infty$ with the subgroup of $\mathbb Z^2$
generated by the bidegrees $(p,1)$. Thus any two monomials
$x_{\boldsymbol p},x_{\boldsymbol q}$ of the same bidegree have the same
image after all variables are inverted. Clearing denominators gives a
finite monomial $x_{\boldsymbol r}$ such that
\[
 x_{\boldsymbol r}x_{\boldsymbol p}
     =x_{\boldsymbol r}x_{\boldsymbol q}\quad\text{in }G .
\]
Equivalently, $J=\bigcup_{\boldsymbol r}(I:x_{\boldsymbol r})$.
If the two monomials represent distinct classes in $G$, their
difference is a zero divisor, since $x_{\boldsymbol r}$ has nonzero
image in $S$.

In fact, in the full polynomial ring $R$ we have
\[
 J=\bigcup_m(I:w_m),\qquad w_m=\prod_{p\in\A_m}x_p .
\]
To see why a squarefree multiplier suffices, take $f\in J$ and a
monomial $x_{\boldsymbol r}$ with $x_{\boldsymbol r}f\in I$. For each
occurrence of a prime $p$ in $\boldsymbol r$ we can find distinct,
arbitrarily large primes $q,q'$ such that $q+q'-p$ is also prime.
Indeed, the three forms $a,b,a+b-p$ have pairwise nonproportional
linear parts. Modulo $2$ take $a=b=1$. At each odd prime $s$ take
$a=1$ and choose $b$ outside the two residues $0,p-1$. Thus the
linear forms theorem \cite[Theorem~1.6]{TT}, applied at scale $4N$ on
$[N,2N]\times[3N,4N]$ for large $N$, gives the desired primes, with proper prime
powers discarded as in Lemma~\ref{depth:prime-pair-conversion}.
Choose each pair above all the preceding pairs. Since
$x_qx_{q'}=x_px_{q+q'-p}$ in $G$, the product of these distinct
variables is a multiple of $x_{\boldsymbol r}$ in $G$ and annihilates
the class of $f$. It divides $w_m$ for sufficiently large $m$. The
reverse inclusion follows from $I\subseteq J$ and $w_m\notin J$,
since $J$ is prime.

Here $m$ may depend on $f$. For the ideal of a fixed finite cutoff,
Theorem~\ref{depth:integral-stabilization} instead identifies the
relation group with the bidegree group once $m$ is sufficiently large.
Consequently $J^{<m}=I^{<m}:w_m^\infty$ for every such $m$.

\subsection{Quotients by a prime variable}\label{depth:prime-variable}
We use the variable cutoff on the odd primes. Throughout this
subsection $R=R^{<m}$, $I=I^{<m}$, $G=R/I$, and $S=R/J$ is the toric
quotient at the same cutoff. For $p\in\A_m$ write $G/x_p$ for
$G/x_pG$. All Betti numbers below are taken over $R$ unless another
polynomial ring is specified.

\paragraph{The cutoff $m=42$.}
Write $P_s=(x_p\mid p\le s,\ p\text{ odd prime})$ in $R^{<42}$.
Over $\mathbb Q$ and $\mathbb F_2$, a primary decomposition gives
\[
 I^{<42}=J^{<42}\cap P_{29}\cap Q,\qquad \sqrt Q=P_{31},
\]
where $J^{<42}$ is the kernel of $x_p\mapsto u^pv$ on $R^{<42}$.
One choice of the embedded component is
$Q=(I^{<42}+(x_{31}^3)):x_{37}^\infty$.
The two minimal primes are $J^{<42}$ and $P_{29}$.
Thus the variable zero divisors are exactly the $x_p$ with $p\le31$.
Both $x_{37}$ and $x_{41}$ are regular. The supplement verifies the
decomposition and the variable colon ideals in \textsc{Singular}.
Here one colon suffices. The calculation over $\mathbb Q$ gives
$I^{<42}:x_3=J^{<42}$. This equality holds over every field: for each
binomial of $J^{<42}$, membership of its $x_3$-multiple in $I^{<42}$
means that its two monomials after multiplication are connected by
binary moves (Proposition~\ref{setup:components}). Such binomials
span $J^{<42}$ over each field. Hence also $I^{<42}:w_{42}=J^{<42}$.
The zero-divisor variables need not be an initial segment at other
cutoffs. Over the same fields, $x_{19}$ is regular in $G^{<30}$
whereas $x_{23}$ is a zero divisor.

One can see the loss of information under localization directly over
every field. Take $f=x_{13}x_{37}^2-x_5x_{41}^2$. Its monomials both
map to $u^{87}v^3$ in the toric ring. At cutoff $42$ the component of
$x_{13}x_{37}^2$ consists of this monomial and $x_{19}x_{31}x_{37}$,
so $f\ne0$ in $G^{<42}$. Over the two fields above the colon ideal
$(I^{<42}:f)$ is exactly $P_{31}$. One can check $x_3 f=0$ directly
by the three moves
\[
 x_3x_{13}x_{37}^2=x_5x_{11}x_{37}^2
 =x_5x_7x_{37}x_{41}=x_3x_5x_{41}^2 .
\]
At cutoff $43$ the new prime $43$ connects the original monomials by the path
\[
 x_{13}x_{37}^2=x_7x_{37}x_{43}=x_3x_{41}x_{43}=x_5x_{41}^2 .
\]
Thus $f$ dies on passing to $G^{<43}$, the first variable cutoff
at which $43$ is available.

\begin{proposition}[Betti numbers after setting a variable to zero]
\label{depth:variable-betti}
Fix $p\in\A_m$ and let $I'=I:x_p$. For $i\ge1$,
\[
 \beta_{i,n,i+1}(G/x_p)
 =\beta_{i,n,i+1}(G)+\beta_{i-1,n-p,i}(S).
\]
If $\operatorname{reg}_R G\le2$, then in every bidegree and for $i\ge0$,
\begin{equation}\label{depth:quotient-betti}
 \beta_{i,n,d}(G/x_p)
 =\beta_{i,n,d}(G)+\beta_{i-1,n-p,d-1}(R/I'),
\end{equation}
where a Betti number in negative homological degree is zero.
In particular the last term is a Betti number of $S$ when $I'=J$.
\end{proposition}
\begin{proof}
Since $I$ contains all quadratic equal-sum binomials on $\A_m$,
the map $G\to S$ is an isomorphism in monomial degrees $0,1,2$.
Also $I\subseteq I'\subseteq J$, because the image of $x_p$ in $S$ is
nonzero. Thus $G\to R/I'\to S$ identifies the linear strands in Tor:
in monomial degree $i+1$, the three Koszul terms computing $\Tor_i$ have
coefficient monomial degrees $0,1,2$. The same comparison applies to
$G/x_p\to S/x_pS$, which again agrees through monomial degree two.
Multiplication by $x_p$ is injective on $S$. The mapping cone of
multiplication by $x_p$ on its minimal resolution is therefore a
minimal resolution of $S/x_pS$, giving the first formula.

For the second formula use the exact sequence
\[
 0\longrightarrow (R/I')(-p,-1)
 \xrightarrow{\ \cdot x_p\ } G\longrightarrow G/x_p\longrightarrow0.
\]
We show that the first map induces zero on Tor. In homological degree
zero this follows from $x_p\in\m$. For $i\ge1$, the first possible
monomial degree in $\Tor_i^R(R/I',\kk)$ is $i+1$, since $I'$ has no linear forms.
In that monomial degree the classes lift from $\Tor_i^R(G,\kk)$ by the linear
strand comparison. The composite $G(-p,-1)\to(R/I')(-p,-1)\to G$
is multiplication by $x_p$, which acts as zero on Tor with $\kk$.
Classes of greater monomial degree have shifted target monomial degree at least $i+3$,
where $\Tor_i^R(G,\kk)$ vanishes by the regularity bound.
The long exact sequence now gives \eqref{depth:quotient-betti}.
\end{proof}

At cutoff $42$, over both $\mathbb Q$ and $\mathbb F_2$, the supplement
verifies $\operatorname{reg}_R G\le2$ and gives
\[
\begin{array}{c|ccc}
 &G&S&G/x_3\\ \hline
 \operatorname{depth}_R &2&1&0\\
 \pd_R &10&11&12
\end{array}
\]
Thus \eqref{depth:quotient-betti} applies with $p=3$ and $R/I'=S$. For
example, $x_5x_7$ is a nonzero class in $G/x_3$ killed by every
variable. It is one of six socle classes, represented by
\[
 x_5x_7,\quad x_5x_{13},\quad x_5x_{19},\quad
 x_5x_{23},\quad x_7x_{23},\quad x_5x_{31}.
\]
Since $\pd_R G=10$, the connecting map identifies these six top
Koszul classes in $\Tor_{12}^R(G/x_3,\kk)_{n,d}$ with
$\Tor_{11}^R(S,\kk)_{n-3,d-1}$, bidegree by bidegree.

\paragraph{Powers of \texorpdfstring{$x_3$}{x3}.}
For a fixed finite cutoff the ascending chain $(I:x_3^j)$ stabilizes
by Noetherianity. There is a simple criterion for its limit to be $J$.
For $m>5$ this happens precisely when every prime $7\le p\le m$ admits
a representation $p+3=q+q'$ with $q,q'>3$ prime. Indeed, $q,q'<p$,
so after inverting $x_3$ the relations $x_3x_p=x_qx_{q'}$ express
each $x_p$ recursively in $x_3^{\pm1},x_5$. The toric images
$u^3v,u^5v$ are algebraically independent, so the resulting localized
map to the localized toric ring is injective. Contracting to $R$ gives
$I:x_3^\infty=J$.

Conversely, if $p$ has no such representation, put $j=(p-3)/2$.
The binomial $x_3^{j-1}x_p-x_5^j$ belongs to $J$. After multiplication
by any power of $x_3$, the first monomial remains isolated under binary
moves: each available pair has sum $6$ or $3+p$, both uniquely represented
on this support. Consequently no power of $x_3$ times this binomial belongs to $I$.
A failure of the pair condition at $p$ therefore persists at every
cutoff containing $p$. It follows that $I^{<m}:x_3^\infty=J^{<m}$ for
arbitrarily large $m$ precisely when every odd prime $p\ge7$ meets the
pair condition. We have no general bound on the power needed when
equality holds, and the results here leave open whether
$I^{<m}:x_3=J^{<m}$ for infinitely many $m$. The equality
$J=\bigcup_m(I:w_m)$ for the full ring allows the multiplier to
depend on the binomial.

\begin{example}[The first failure of a single colon]\label{depth:cutoff132}
At every nonempty variable cutoff before $m=131$ we have $I:x_3=J$.
At $m=131$, set
\[
 f=x_{127}^6-x_{107}x_{131}^5,\qquad I'=I:x_3.
\]
The two monomials have bidegree $(762,6)$, and
\[
 f\notin I',\qquad J=I'+(f)=I:x_3^2.
\]

To see why one copy of $x_3$ fails, consider the two partitions of
$765$ obtained after multiplying by $x_3$. Starting on the left, two legal moves give
\[
 x_3x_{127}^6 \longrightarrow x_{23}x_{107}x_{127}^5
 \longrightarrow x_{19}x_{107}x_{127}^4x_{131}.
\]
Sort the primes in increasing order and then reduce each entry modulo $3$.
The two moves give the residue words
\[
 (0,1,1,1,1,1,1)\ \longrightarrow\
 (2,2,1,1,1,1,1)\ \longrightarrow\
 (1,2,1,1,1,1,2).
\]
The first move, $3+127=23+107$, replaces $(0,1)$ by $(2,2)$ and
uses the only $3$. The second, $23+127=19+131$, preserves the
residue counts even though it changes the primes. To turn two entries
congruent to $1$ back into a $3$ and an entry congruent to $2$ at
this vertex, the only candidate sums are $19+127=3+143$ and
$127+127=3+251$. The first requires the composite $143$ and the
second requires $251>131$. The two entries congruent to $2$ would need
$107+131=3+235$, with $235=5\cdot47$. Thus no next move from this
vertex can reintroduce $3$. Here the residue counts constrain which moves are possible. To follow a
loop as a Koszul cycle, Example~\ref{ex:weight33} also keeps track of
which prime is fixed on each edge and the order of the two fixed primes
at each corner.

Every partition in this fiber has seven entries and at most one $3$:
two copies of $3$ would give a sum at most $6+5\cdot131=661$.
Write $h$ for the count of entries congruent to $1$ and $z$ for the
number of threes. Reduction of the sum modulo $3$ gives
$2h-(7-z)\equiv0\pmod3$. Thus
\[
 \begin{array}{c|c}
 z & h\text{ allowed by the congruence}\\ \hline
 0 & 2,\ 5\\
 1 & 0,\ 3,\ 6
 \end{array}
\]
A move replaces only two entries, so it changes $h$ by at most two.
The two endpoints have $h=6$ and $h=0$, so any path between them must pass through $h=3$. Such a partition
would contain one $3$ and \emph{six other entries}, three from each
nonzero residue class. These six entries must sum to $765-3=762$.
Their largest possible sum is $3\cdot127+3\cdot131=774$. The
next-largest primes in the two classes are $109$ and $113$, respectively.
Any reduction therefore costs at least $18$, whereas the required
reduction is $12$. This excludes the intermediate partition and
proves that $x_3f\notin I$.

A second $3$ allows us to continue after the same two moves from the
left. Use the remaining $3$ in $3+19=11+11$, and then use
$11+127=7+131$ and $7+127=3+131$ twice in succession. This gives
$x_3^2x_{107}x_{131}^5$, hence $x_3^2f\in I$. The obstruction is
specific to the cutoff. When $137$ enters, a path with just one $3$
can pass through $(3,113,127,127,127,131,137)$. Its residue word
$(0,2,1,1,1,2,2)$ has $h=3$, with all three entries congruent to
$1$ coming from $127$.

For completeness, the supplement checks each odd-prime prefix through
$131$. It verifies the pair condition above and uses \textsc{Singular}
over $\mathbb Q$ to compare consecutive colons. They agree at the first
power for each earlier prefix. At the prefix ending at $131$, the second and third colons agree,
and the second is $I'+(f)$. These equalities hold over every field.
Indeed, membership of a binomial in $I:x_3^j$ is connectivity of its
endpoints after multiplication by $x_3^j$. These colons are spanned by
such binomials, by the component basis. Adjoining $f$ identifies the
components joined by its monomial multiples, again independently of
the field.

The calculation also gives $I':f=(x_p\mid p\le113)$ at this cutoff.
Each of those variables kills the class of $f$ over every field by
the same connectivity test. Every $x_{127}^a x_{131}^b f$ survives,
and distinct pairs $(a,b)$ give different bidegrees. It follows that
\[
 J/I'\cong\kk[x_{127},x_{131}](-762,-6).
\]
This describes the difference between the actual term $R/I'$ of the
colon sequence and the toric quotient $S$. The Koszul complex on the
$29$ variables through $113$ resolves $J/I'$, with the indicated shift.
Consequently $0\to J/I'\to R/I'\to S\to0$ induces isomorphisms
$\Tor_i^R(R/I',\kk)_{n,d}\cong\Tor_i^R(S,\kk)_{n,d}$ outside the two
monomial degrees $i+5$ and $i+6$.
\end{example}

\subsection{Dual Koszul cocycles and the prime-pair condition}
\label{depth:dual-cocycles}

The depth bound can also be approached by constructing classes in the
dual Koszul complex. At a variable cutoff write $R_m=R^{<m}$ and
$G_m=G^{<m}$, and let $\mathbb K_{\kk}^{<m}$ be the Koszul resolution of
$\kk$ over $R_m$. Put
\[
 C_m=\operatorname{Hom}_{R_m}(\mathbb K_{\kk}^{<m},G_m).
\]
Write $C_{m,n,d}$ for its complex in internal bidegree $(n,d)$ and
$C^i_{m,n,d}$ for the term in cohomological degree $i$.
In this subsection we use $\varepsilon_q$ for the dual exterior
generators, of cohomological degree one and internal bidegree
$(-q,-1)$. With the transpose convention for the codifferential,
\[
 \mathrm{d}^\star\omega=\theta_m\wedge\omega,
 \qquad \theta_m=\sum_{q\in\A_m}x_q\varepsilon_q.
\]
Thus $\mathrm{d}^\star$ preserves internal bidegree, and
$\mathrm H^i(C_m)=\operatorname{Ext}_{R_m}^i(\kk,G_m)$.
The least $i$ with nonzero cohomology is $\delta^{<m}$.
If $s_m=\sum_{q\in\A_m}q$, finite Koszul self-duality gives
\begin{equation}\label{depth:dual-betti}
 \dim_\kk\mathrm H^i(C_m)_{n,d}
 =\beta^{<m}_{|\A_m|-i,\,s_m+n,\,|\A_m|+d}.
\end{equation}

Suppose $m\ge11$ and every prime $5\le q\le m$ has a representation
$q+1=p_1+p_2$ by odd primes. Both summands are at most $q-2$, so they
belong to $\A_m$. Choose $p_1\le p_2$ with $p_1p_2$ maximal, using
the product order of Proposition~\ref{moves:order}, and set
\[
 y_q=x_{p_1}x_{p_2}\quad\text{in }G_m.
\]
The maximum-product choice fixes a representative; the class $y_q$
is independent of this choice by the defining quadratic relations.
In particular,
\[
 y_5=x_3^2,\qquad y_7=x_3x_5,\qquad y_{11}=x_5x_7.
\]
Consider the cochain
\begin{equation}\label{depth:dual-family}
 \omega_m=\varepsilon_3\wedge
       \sum_{\substack{5\le q\le m\\q\text{ prime}}}y_q\varepsilon_q
       \in C^2_{m,-2,0}.
\end{equation}

\begin{proposition}\label{depth:dual-family-criterion}
The cochain $\omega_m$ is closed if and only if
\begin{equation}\label{depth:dual-cubic-identities}
 x_qy_{q'}=x_{q'}y_q\quad\text{in }G_m
 \qquad(5\le q<q'\le m).
\end{equation}
Equivalently, if $q+1=p_1+p_2$ and $q'+1=p'_1+p'_2$ are the
chosen representations, then
\[
 x_q(x_{p'_1}x_{p'_2})=x_{q'}(x_{p_1}x_{p_2})\quad\text{in }G_m.
\]
Whenever it is closed, its class is nonzero over every field. In
particular,
\[
 \beta^{<m}_{|\A_m|-2,\,s_m-2,\,|\A_m|}\ge1,
 \qquad \delta^{<m}\le2.
\]
\end{proposition}
\begin{proof}
Exterior multiplication gives
\[
 \mathrm{d}^\star\omega_m=
 \sum_{5\le q<q'\le m}(x_{q'}y_q-x_qy_{q'})
        \varepsilon_3\varepsilon_q\varepsilon_{q'},
\]
which proves the criterion. A possible primitive in internal
bidegree $(-2,0)$ has the form $\omega'=\sum_q f_q\varepsilon_q$, with
$f_q\in(G_m)_{q-2,1}$. Since $1$ and $9$ are not prime labels,
$f_3=f_{11}=0$. The coefficient of
$\varepsilon_3\varepsilon_{11}$ in every such $\mathrm{d}^\star\omega'$
is therefore $x_3f_{11}-x_{11}f_3=0$. In $\omega_m$ that coefficient is the nonzero
monomial class $x_5x_7$. The component basis of
Proposition~\ref{setup:components} proves its nonvanishing over every
field. Equation~\eqref{depth:dual-betti} and the depth characterization
give the final assertions.
\end{proof}

\begin{example}[Finite cocycles]\label{depth:dual-finite}
At cutoff $12$ the cocycle is
\[
 \omega_{12}=x_3^2\varepsilon_3\varepsilon_5
       +x_3x_5\varepsilon_3\varepsilon_7
       +x_5x_7\varepsilon_3\varepsilon_{11}.
\]
Its closure follows from $x_3x_7=x_5^2$ and $x_3x_{11}=x_7^2$.
For example,
$x_5y_{11}=x_5^2x_7=x_3x_7^2=x_3^2x_{11}=x_{11}y_5$.
The following choices extend the displayed expression by one term
whenever a new prime is adjoined:
\[
\begin{array}{c|ccccc}
 q&13&17&19&23&29\\ \hline
 y_q&x_7^2&x_7x_{11}&x_7x_{13}&x_{11}x_{13}&x_{13}x_{17}
\end{array}
\]
Thus $\omega_{14}=\omega_{12}+x_7^2\varepsilon_3\varepsilon_{13}$ and
$\omega_{18}=\omega_{14}+x_7x_{11}\varepsilon_3\varepsilon_{17}$, and
similarly at cutoffs $20$, $24$, and $30$.

The supplement verifies~\eqref{depth:dual-cubic-identities} at every
odd-prime prefix through $101$. Hence $\omega_m$ is a nonboundary cocycle
over every field for every integer $11\le m\le102$.
The verification enumerates cubic monomials and their binary-move
components at each support; it does not replace $G_m$ by its toric
quotient. The coefficient test in the preceding proof certifies
nontriviality independently of the field.

At the smaller cutoff $8$ there is instead a linear-coefficient class
\[
 \omega_8=x_3\varepsilon_3\varepsilon_5+x_5\varepsilon_3\varepsilon_7
       \in C^2_{8,-5,-1}.
\]
Its codifferential is
$\mathrm{d}^\star\omega_8=(x_3x_7-x_5^2)\varepsilon_3\varepsilon_5\varepsilon_7=0$.
The only possible primitive in this bidegree is a scalar multiple
of $\varepsilon_5$; the codifferential of every nonzero such multiple
has a nonzero $\varepsilon_5\varepsilon_7$ coefficient, so the displayed class is
nonzero.
\end{example}

\begin{conjecture}[Infinite dual Koszul cocycle]\label{depth:dual-conjecture}
For every odd prime $q\ge5$, there are odd primes $p_1,p_2$ with
$q+1=p_1+p_2$. Let $y_q=x_{p_1}x_{p_2}\in G_{q+1,2}$ be the
resulting quadratic class in the full binary Goldbach ring. Then
\[
 x_qy_{q'}=x_{q'}y_q\quad\text{in }G\qquad(5\le q<q').
\]
\end{conjecture}

For infinitely many variables,
$C^i=\operatorname{Hom}_R(\mathbb K_{\kk,i},G)$ is a product over
the $i$-element exterior supports. Each coordinate of $\mathrm{d}^\star$ uses
only finitely many input coordinates. Under the conjecture,
\[
 \omega=\varepsilon_3\wedge\sum_{q\ge5}y_q\varepsilon_q
\]
therefore defines a cocycle in $C^2$. It also belongs to the
completion of the direct-sum cochain module in the Goldbach weight
of its coefficients: the coefficient $y_q$ has weight $q+1$.
This is a coefficient filtration; every whole term has dual internal
bidegree $(-2,0)$. The same $\varepsilon_3\varepsilon_{11}$
coefficient test proves that $\omega$ is not a boundary, even in this
completion. In particular, the conjecture gives
$\operatorname{Ext}_R^2(\kk,G)_{-2,0}\ne0$.

The conjecture concerns identities in the full ring. Requiring every
$\omega_m$ to be closed already in $G^{<m}$ is stronger; the finite
examples have this additional property. A cubic identity of weight
$q+q'+1$ may in general use primes as large as $q+q'-5$ along a move
path, beyond the original labels $q,q'$.

\begin{proposition}\label{depth:dual-prime-pair}
Conjecture~\ref{depth:dual-conjecture} implies that every prime
$q\ge7$ has a representation
\[
 q+3=p_1+p_2,\qquad p_1,p_2>3\text{ prime}.
\]
Consequently it implies $I^{<m}:x_3^\infty=J^{<m}$ at every
variable cutoff $m>5$.
\end{proposition}
\begin{proof}
The identity with the label $5$ is $x_5y_q=x_3^2x_q$.
If $q+3$ has no representation by primes greater than $3$, the
monomial $x_3^2x_q$ is isolated under binary moves: its only pair
sums are $6$ and $q+3$, each uniquely represented. For $q\ge7$ the
monomial $x_5y_q$ is different, since it contains $x_5$.
The component basis therefore contradicts the required equality.
Moreover, every prime occurring along a cubic path of weight $q+6$
is at most $q$, so larger cutoffs cannot repair a failure of this
particular identity. The saturation conclusion follows from the
prime-pair criterion in Appendix~\ref{depth:prime-variable}.
\end{proof}

The relation-group theorem does not supply the cubic identities in
Conjecture~\ref{depth:dual-conjecture}. When $y_q,y_{q'}$ exist,
$x_qy_{q'}-x_{q'}y_q$ belongs to $J$, so that theorem supplies a common
monomial multiplier annihilating the difference in $G$.
It does not permit cancellation of that multiplier. Even saturation
by powers of $x_3$ only gives annihilation by some power, rather
than equality in $G$. Thus the conjecture asks for specific cubic
binary-move paths beyond the stable connectivity already proved.

\end{document}